\documentclass[11pt,a4paper]{article}
\usepackage{amsfonts,amsmath,amssymb,amscd}
\usepackage{fancybox}
\usepackage[latin1]{inputenc}
\usepackage[usenames]{color}
\usepackage{graphicx}
\usepackage{enumerate}
\usepackage{float}
\usepackage{fullpage}
\usepackage{amssymb}
\usepackage{subcaption}

\usepackage{euscript,color}

\definecolor{red}{rgb}{1,0,0}

\newtheorem{thm}{Theorem}[section]

\newtheorem{lem}[thm]{Lemma}
\newtheorem{prop}[thm]{Proposition}
\newtheorem{defn}[thm]{Definition}
\newtheorem{rem}[thm]{Remark}
\newtheorem{ex}[thm]{Example}

\def \r{\mbox{${\mathbb R}$}}

\begin{document}
	
\title{\textbf{Closed $1/2$-Elasticae in the $2$-Sphere}}

\author{E. Musso and A. P\'ampano}
\date{\today}

\maketitle

\begin{abstract} 
We study critical trajectories in the sphere for the $1/2$-Bernoulli's bending functional with length constraint. For every Lagrange multiplier encoding the conservation of the length during the variation, we show the existence of infinitely many closed trajectories which depend on a pair of relatively prime natural numbers. A geometric description of these numbers and the relation with the shape of the corresponding critical trajectories is also given.\\


\end{abstract}

\section{Introduction}

Functionals on curves depending on the curvatures are ubiquitous in differential geometry, analysis, mathematical physics and biomathematics. Their study dates back to the days of the Bernoulli family and Euler. Indeed, the problem of determining the bending deformation of rods was first formulated by J. Bernoulli in 1691 (\cite{L}) and the possible qualitative types for untwisted plane configurations were completely described by L. Euler (\cite{E}), although some particular cases were already known to J. Bernoulli (\cite{Be}). More generally, in a letter to L. Euler in 1738 (\cite{T}), D. Bernoulli proposed to investigate the extrema of the functionals
$$\mathcal{B}_{p,\lambda}:\gamma\longmapsto\int_\gamma \left(\lvert\kappa\rvert^{p}+\lambda\right)$$
where $\kappa$ is the curvature of the curve $\gamma$ and $\lambda$ is a Lagrange multiplier, encoding the conservation of the length during the variation.$^1$\footnote{$^1$Actually, D. Bernoulli considered the unconstrained case, that is, $\lambda=0$.} Despite the ancient origin, these variational problems are still a very vital field of research. In the current literature, critical curves of these functionals are usually called {\em p-elasticae} (or, {\em free p-elasticae} when $\lambda=0$).

Leaving aside many technical aspects, we now describe some general heuristic ideas which will guide us in dealing with the problems studied in this paper. Consider a functional for immersed curves in a 2-dimensional Riemannian space form ${\rm M}_c^2$ of constant curvature $c$, with Lagrangian $P(\kappa)+\lambda$, where $P$ is a function of class $\mathcal{C}^2$, $\lambda$ is a Lagrange multiplier and $\kappa$ is the geodesic curvature of the curve. The functional is acting on curves of class $\mathcal{C}^4$ parameterized by the arc-length, belonging to the open set in the Whitney's topology, defined by the condition that $\kappa$ takes values in the domain of $P$. We denote by ${\rm G}_c$ the isometry group of ${\rm M}^2_c$. From standard computations for functionals depending on the curvature (see, for instance, \cite{AGM} and references therein) the critical curves, for compactly supported smooth variations, are characterized by the associated Euler-Lagrange equation which, in turn, admits the conservation law $\Phi_{\lambda}(\kappa,\dot{\kappa})=d$, where $\dot{\kappa}$ represents the derivative of $\kappa$ with respect to the arc-length, $d$ is a constant of integration and 
$$\Phi_{\lambda}(x,y)=(P_{xx})^2 y^2+\left(x P_x-P-\lambda\right)^2+c P_x^2\,.$$ 

Geometric methods can be successfully applied if one can find a contact transformation ${\mathtt f}:(x,y)\to (f(x),\dot{f}(x)y)$ such that the equations of the level curves $\Phi_{\lambda}\circ {\mathtt f}=d$ can be reduced to the form $y^2+x^nQ_{\lambda,d}(x)=0$, where $n\in\{0,2\}$ and $Q_{\lambda,d}$ is, for generic values of $\lambda$ and $d$, a polynomial with simple roots such that $Q_{\lambda,d}(0)\neq 0$. The algebraic curves $y^2+x^nQ_{\lambda,d}(x)=0$ are called the \emph{phase curves}, while their connected components are the \emph{reduced phase curves}.
The phase curves are the pre-images under the contact transformation ${\mathtt f}$ of the usual phase portraits (\cite{G}) of the variational problem,  that is the level sets of the function  $\Phi_{\lambda}$. Their main use is in the construction of solutions of the Euler-Lagrange equation.  Indeed,  solutions of the ordinary differential equation $\dot{\mu}^2=-\mu^nQ_{\lambda,d}(\mu)$ can be built inverting incomplete hyperelliptic integrals$^2$\footnote{$^2$The term hyperelliptic is used here in a broad sense, including also the rational and elliptic cases.} and the function $ \kappa=f\circ \mu$ solves the Euler-Lagrange equation. Note that $c$ can modify the algebro-geometric properties of the phase curves with consequences on the analytical behavior of the  integrals. This is a first place where the geometry of the ambient space enters into play.

Once that we know $\kappa_{\lambda,d}$,  the corresponding critical curves $\gamma_{\lambda,d}$ can be found by quadratures. This is trivial in the Euclidean case, i.e., when $c=0$. For $c\neq 0$ the problem can be faced using the momentum map for the Hamiltonian action of ${\rm G}_c$ on the phase space. Clearly, the curvature $c$ also affects the geometry of the critical curves at this point. 

To find closed critical curves, the first step is to describe the domain ${\mathcal U}\subset \r^2$ of the parameters $\lambda,d$ such that $\kappa_{\lambda,d}$ is periodic. If $(\lambda,d)\in {\mathcal U}$,  the least period  $\omega_{\lambda,d}$  of $\kappa_{\lambda,d}$ can be evaluated in terms of complete hyperelliptic integrals. In a second step we need to extract from critical curves with periodic curvatures those that are closed. The behavior of the adjoint representation of ${\rm G}_c$ forces to follow different procedures depending on whether $c=0$, $c<0$ or $c>0$. This is a third point where the value of $c$ comes into play in an essential way. 

Since the Lagrangian of the unconstrained $p$-Bernoulli's bending functionals is positively homogeneous, only the sign of the Gaussian curvature of ${\rm M}_c$ plays a role. Similarly, for the constrained case, possibly modifying the value of $\lambda$, we reach the same conclusion. Thus, without loss of generality, one can assume $c=0$ or $c=\pm1$.

We now focus on the case $c=1$, which is the relevant one for the purposes of the paper. In order to formulate the closure conditions one has to analyze the {\em period map}, which is defined as follows. Let ${\mathcal F}_{\lambda,d}$ be the spherical Frenet frame along $\gamma_{\lambda,d}$.  The monodromy of $\gamma_{\lambda,d}$ is the element of ${\rm G}_c\equiv SO(3)$ given by ${\mathfrak m}(\gamma_{\lambda,d}):={\mathcal F}_{\lambda,d}(\omega_{\lambda,d}){\mathcal F}_{\lambda,d}(0)^{-1}$.  Choosing carefully the initial data, we may assume that ${\mathfrak m}(\gamma_{\lambda,d})$ belongs to $SO (2)$.  This gives rise to a differentiable map ${\mathfrak m}:{\mathcal U}\to SO (2)$. The  period map is then a continuous lift $\widehat{\Psi}$ of ${\mathfrak m}$ to ${\mathbb R}$. Possibly, $\widehat{\Psi}$ can be evaluated in terms of complete hyperelliptic integrals. Closed solutions arise when $\widehat{\Psi} (\lambda,d)=2\pi q$ where $q\in {\mathbb Q}$ is a rational number. In the most favorable cases, $\widehat{\Psi}$ has maximal rank and its fibers do intersect transversely the horizontal lines.  This means that, for every admissible $\lambda$, there exist countably many distinct equivalence classes of closed critical curves with multiplier $\lambda$. In addition, for every $q\in \mathbb{Q}$ such that $2\pi q\in \widehat{\Psi}({\mathcal U})$,  the level set  $\widehat{\Psi}^{-1}(2\pi q)$ is an embedded curve. The rational number $q=m/n$ encodes two pieces of relevant geometric information: $n$ is the order of the stabilizer of the trajectory while, generically, $m$ is the homotopy class of the critical curve in the sphere punctured at the north and south poles (or, equivalently, the linking number with the vertical axis).

Now we  discuss (without claiming to be exhaustive) some known results about $p$-elasticae fitting in or related to the above heuristic scheme.

The case $p=2$ corresponds with the classical elastic curves, which have received a considerable interest in the last decades because of their applications in the following research topics: Willmore surfaces (\cite{LS2,Pin1,VDM}), constrained Willmore surfaces (\cite{BPP,SCDPS}), the Canham-Helfrich-Evans model for lipid bilayers (\cite{Can,CGR1,Ev,He,JMN,TY}), and surfaces with spherical curvature lines (\cite{CPS}), among others. The existence and the geometric properties of closed elasticae have been extensively treated in several works published in the 1980s (see, for instance, \cite{LS2}). Another source of interest is due to the interrelationships with integrable flows of curves governed by the mKdV hierarchy (\cite{CQ1,GP1,GP2,LPe,M3,NSW}). In fact, elasticae do evolve maintaining their shape under the first non-trivial Goldstein-Petrich flow. A similar relation between elastic curves and the nonlinear Schr\"{o}dinger equation can be obtained applying the Hasimoto transformation (\cite{Ha1,Ha2,K,LangerPerline}). These phenomena occur in other contexts such as Lorentzian, centro-affine, equi-affine, projective and pseudo-conformal geometries (\cite{CI,CIMB1,CIMB2,CQ2,CQ3,M2,M4,MN4,MS,Pin2,V}). 

For natural values of $p>2$, free $p$-elasticae have been considered in \cite{AGM}. They have been used to construct Willmore-Chen submanifolds in spaces with Riemannian and pseudo-Riemannian warped product metrics (\cite{ABG,BFLM}) and have been applied to analyze conformal tensions in string theories (\cite{BFL}). In the case of spherical curves, the only closed critical trajectories are geodesics. Of course, this assertion which was proven in \cite{AGM} necessitates the assumption that critical curves are $\mathcal{C}^4$, which is the setting we are interested on in this paper. With a suitable choice of the contact transformation and for generic values of $\lambda$ and $d$, the reduced phase curves are real cycles of smooth hyperelliptic curves of genus $p-1$.

In a recent paper (\cite{New}) Miura-Yoshizawa considered the functionals $\mathcal{B}_{p,\lambda}$ for curves in the Euclidean plane, for every real number $p>1$. Since the contact transformation they used leads, in general, to non-algebraic phase curves,$^3$\footnote{$^3$In \cite{New} there is no explicit mention to phase curves, but they are implicitly defined on page $27$.} they developed a new technical tool to handle the problem, namely, the concept of $p$-elliptic functions. As a consequence, they obtain a complete and remarkable classification of critical curves in the plane. As a corollary of their beautiful results, they proved that planar closed critical curves are either circles or lemniscates, as in the classical case $p=2$. However, their approach cannot be directly applied to the case $p\in (0,1)$ due to the appearance of a singularity at the origin in the phase curves. 

Indeed, when $p\in(0,1)$, less is known about $p$-elasticae. From an analytical viewpoint this is understandable since the Lagrangian is merely continuous at the origin; while, from a geometric viewpoint, the main difference with the case $p>1$ is the above mentioned singularity at the origin in the phase curves. Thus, a word of caution must be spent to clarify the use of the term ``critical curve" in our context. A critical curve is a curve of class at least ${\mathcal C}^4$, which is stationary with respect to compactly supported variations whose support does not contain inflection points.$^4$\footnote{$^4$Borrowing the terminology from the classical Euclidean geometry, $s$ is an inflection point if $\kappa(s)=0$ and $\gamma$ is said convex if $\kappa>0$ everywhere.} In particular, the variational problem can be faced for convex curves and for appropriate choices of $p\in {\mathbb Q}$. As suggested in \cite{AGM} finding closed critical trajectories for $p\in(0,1)$ also merits investigation. In \cite{AGM,AGP,MOP,MP}, for $p=(n-2)/(n+1)$, $n\in {\mathbb N}$, $n>2$ and $\lambda=0$, the existence of infinitely many closed free $p$-elasticae in ${\mathbb S}^2$ was shown. These curves arise in the theory of biconservative hypersurfaces as the generating curves of rotational ones (\cite{MOP,MP}). The existence of countable infinitely many trajectories can be seen as a consequence of the fact that, with a suitable choice of the contact transformation and for a generic $d$, the reduced phase curves are real cycles of singular hyperelliptic curves.

Another class of Bernoulli's functionals fitting in our scheme is when $p=1/n$ is the reciprocal of a natural number $n>1$. The unconstrained cases with $p=1/2$, or $p=1/3$, and $c=0$ were considered by W. Blaschke (\cite{B}, Vol I, 1921, and Vol II, 1923) who showed that the critical curves are catenaries ($p=1/2$) or parabolas ($p=1/3$). On one hand, the case $p=1/3$ and $\lambda=0$ corresponds with the equi-affine length for convex curves. After the seminal paper \cite{COT}, equi-affine geometry of convex curves has been consistently used in recent studies on human curvilinear $2$-dimensional drawing movements and recognition for non-rigid planar shapes (see for instance \cite{FH, RK} and the literature therein). On the other hand, all closed free $1/2$-elasticae in ${\mathbb S}^2$ were recently found in \cite{AGM, AGP1}. Critical curves for an extension of this case in Riemannian and Lorentzian $3$-space forms have been characterized as the profile curves of invariant constant mean curvature surfaces (\cite{AGP}). 

To analyze the family of Bernoulli's functionals with $p=1/n$, one can resort to the contact transformation $(x,y)\to (x^n,nx^{n-1}y)$. The equation of the corresponding phase curves can then be written as
$$y^2+\frac{x^2}{\left(n-1\right)^2}\left((n-1)^2x^{2n}+2n(n-1)\lambda x^{2n-1}+n^2(\lambda^2-d)x^{2(n-1)}+c\right)=0\,.$$
Thus, in general, the reduced phase curves are real cycles of singular hyperelliptic curves of genus $n-1$ (if $c\neq 0$). Among all the possible values of $p=1/n$, the only one that makes the phase curve a singular elliptic curve is $p=1/2$. For this reason, $p=1/2$ can be considered to play the role of the classical $p=2$ among the possible values of $p\in(0,1)$. And, as in the classical case $p=2$, their study can be faced resorting to elliptic functions and integrals. Although in the constrained case there exist $1/2$-elasticae in $\mathbb{R}^2$ with non-constant periodic curvature, it is quite easy to see that none of them are closed. The reason lies in the fact that the phase curves are rational curves with an isolated singularity at the origin (see Appendix B). These critical curves were geometrically described in Theorem 6.1 of \cite{LP}. Note that the present case corresponds to $n=-2$ in the notation of \cite{LP}. Moreover, their shapes are similar to the one represented in Figure 3 of \cite{LP}.

Motivated by the above mentioned results of Arroyo, Garay, Menc\'ia and P\'ampano about the existence of closed free $1/2$-elasticae in ${\mathbb S}^2$ (\cite{AGM,AGP,AGP1}), this paper aims to investigate the constrained case. In perspective, it can be seen as a first step towards a general analysis of closed $1/n$-elasticae. One of the nice features of the spherical case is that the singular points of the phase curves are isolated. This implies that critical curves, with the exception of geodesics, do not have inflection points so, possibly after inverting the orientation, they are convex. Consequently, we may restrict to the space of convex curves without any loss of generality and bypass the problem of having a Lagrangian which is not of class $\mathcal{C}^1$ at the origin. Our goal is, beside proving the existence for every value of $\lambda$ of countably many smooth closed critical curves, implement effective methods to identify and analyze the global geometric properties of closed $1/2$-elasticae.

For the sake of brevity, a critical curve for $\mathcal{B}_{1/2,\lambda}$ with positive, non-constant periodic curvature is said to be a \emph{B-curve}. If, in addition, the curve is periodic, it is called a \emph{B-string}. We next state the main results of the paper. The first one is the following theorem.

\begin{thm}\label{A} Let $\gamma:I\subseteq \mathbb{R}\longrightarrow\mathbb{S}^2$ be a critical curve for $\mathcal{B}_{1/2,\lambda}$ with non-constant curvature $\kappa$ defined on its maximal domain $I\subseteq\mathbb{R}$. Then, $I=\mathbb{R}$ and, possibly reversing the orientation, $\gamma$ is a B-curve  
that  can be parameterized in terms of its arc-length parameter $s\in\mathbb{R}$, up to rigid motions, as
$$\gamma(s)\equiv\gamma_\xi(s)=\frac{1}{2\xi\mu}\left(1,-\sqrt{4\xi^2\mu^2-1}\cos\theta(s),\sqrt{4\xi^2\mu^2-1}\sin\theta(s)\right),$$
where
$$\theta(s):=2\xi\int \frac{\mu^2\left(\mu+2\lambda\right)}{1-4\xi^2\mu^2}\,ds\,,$$
and $\mu\equiv \mu(s)=\sqrt{\kappa(s)}$ is a solution of 
$$\dot{\mu}^2=-\mu^2\left(\mu^4+4\lambda\mu^3+4\left[\lambda^2-\xi^2\right]\mu^2+1\right),$$
for suitable constant $\xi>0$.
\end{thm}

We will say that two B-curves are \emph{equivalent} if there is a rigid motion taking one into another. This theorem implies that $\gamma_\xi$, and so all the curves in the equivalence class of $\gamma_\xi$, are B-strings if and only if the elliptic integral
$$\Psi_\lambda(\xi):=2\xi\int_0^\omega\frac{\mu^2\left(\mu+2\lambda\right)}{1-4\xi^2\mu^2}\,ds\,,$$
is a rational multiple of $2\pi$. Here, $\omega$ denotes the least period of $\mu$. Although the function $\Psi_\lambda$ has a jump discontinuity, it can be regularized to a continuous function $\widehat{\Psi}_\lambda$ such that $\widehat{\Psi}_\lambda\equiv \Psi_\lambda ({\rm mod}\, 2\pi {\mathbb Q})$.  Let $\gamma_{\xi}$ be a B-curve with multiplier $\lambda$ such that $\widehat{\Psi}_\lambda(\xi)= 2\pi q$,  $q=m/n\in {\mathbb Q}$, ${\rm gcd}(m,n)=1$. The function $\widehat{\Psi}_\lambda$ plays the role of the period map. We say that $q$ is the \emph{characteristic number} of  $\gamma_{\xi}$.  The number $n$ is said to be the \emph{wave number}.

Our second main result deals with the existence of B-strings.

\begin{thm}\label{C} For every $\lambda\in\mathbb{R}$ the image of $\widehat{\Psi}_\lambda$ is an open interval $I_\lambda$. Thus, for any $q=m/n$ such that $2\pi q\in I_\lambda$ and every $\xi \in \widehat{\Psi}_\lambda^{-1}(2\pi q)$, $\gamma_{\xi}$ is a B-string with multiplier $\lambda$ and characteristic number $q$.
\end{thm}

The proof of this theorem is based on the analysis of the asymptotic behavior of the complete elliptic integral $\widehat{\Psi}_\lambda$. The numerical experiments strongly support the \emph{ansatz} that $\widehat{\Psi}_\lambda$ is a strictly decreasing function of $\xi$ and that $I_\lambda\subset (0,\pi)$. The validity of the ansatz would imply that for each pair of relatively prime natural numbers $(m,n)$ such that $2\pi m/n\in I_\lambda$ there exists a \emph{unique} equivalence class of B-strings with multiplier $\lambda$ and characteristic number $q=m/n$. Moreover, this would also show that all B-strings have self-intersections and wave number $n\geq 3$.

Finally, the third result is about basic geometric features of B-strings.

\begin{thm}\label{B} Let $\gamma_\xi$, $\xi>0$, be a suitable representative of a equivalence class of B-strings with multiplier $\lambda$ and characteristic number $q=m/n$. Then, the following conclusions hold true:
\begin{enumerate}
\item The trajectory of $\gamma_\xi$ is invariant by the group generated by rotation of $2\pi/n$ around the $Ox$-axis and it is contained in a region of the upper hemisphere $\mathbb{S}^2_+=\{(x,y,z)\in\mathbb{S}^2\,\lvert\,x>0\}$ bounded by two horizontal planes.
\item If $4\lambda\xi+1\neq 0$, then $\gamma_\xi$ does not intersect the $Ox$-axis. Moreover:
\begin{enumerate}
\item If $4\lambda\xi+1>0$, $n-m$ is the linking number with the $Ox$-axis (equipped with the upward orientation) and $\gamma_\xi$ possesses, exactly, $n(n-m-1)$ ordinary double points.
\item If $4\lambda\xi+1<0$, $-m$ is the linking number with the $Ox$-axis (equipped with the upward orientation) and $\gamma_\xi$ possesses, at least, $n m$ points of self-intersection.
\end{enumerate}
\item If $4\lambda\xi+1=0$ (necessarily, $\lambda<0$), then $\gamma_\xi$ intersects the $Ox$-axis $n$ times and the moving point $\gamma_\xi(s)$ travels counter-clockwise around the $Ox$-axis (equipped with the upward orientation). In this case $n-m$ is the turning number of the plane projection of $\gamma_\xi$ into the plane $x=0$.
\end{enumerate}
\end{thm}

\begin{figure}[h!]
\centering
\begin{subfigure}[b]{0.3\linewidth}
\includegraphics[height=5cm]{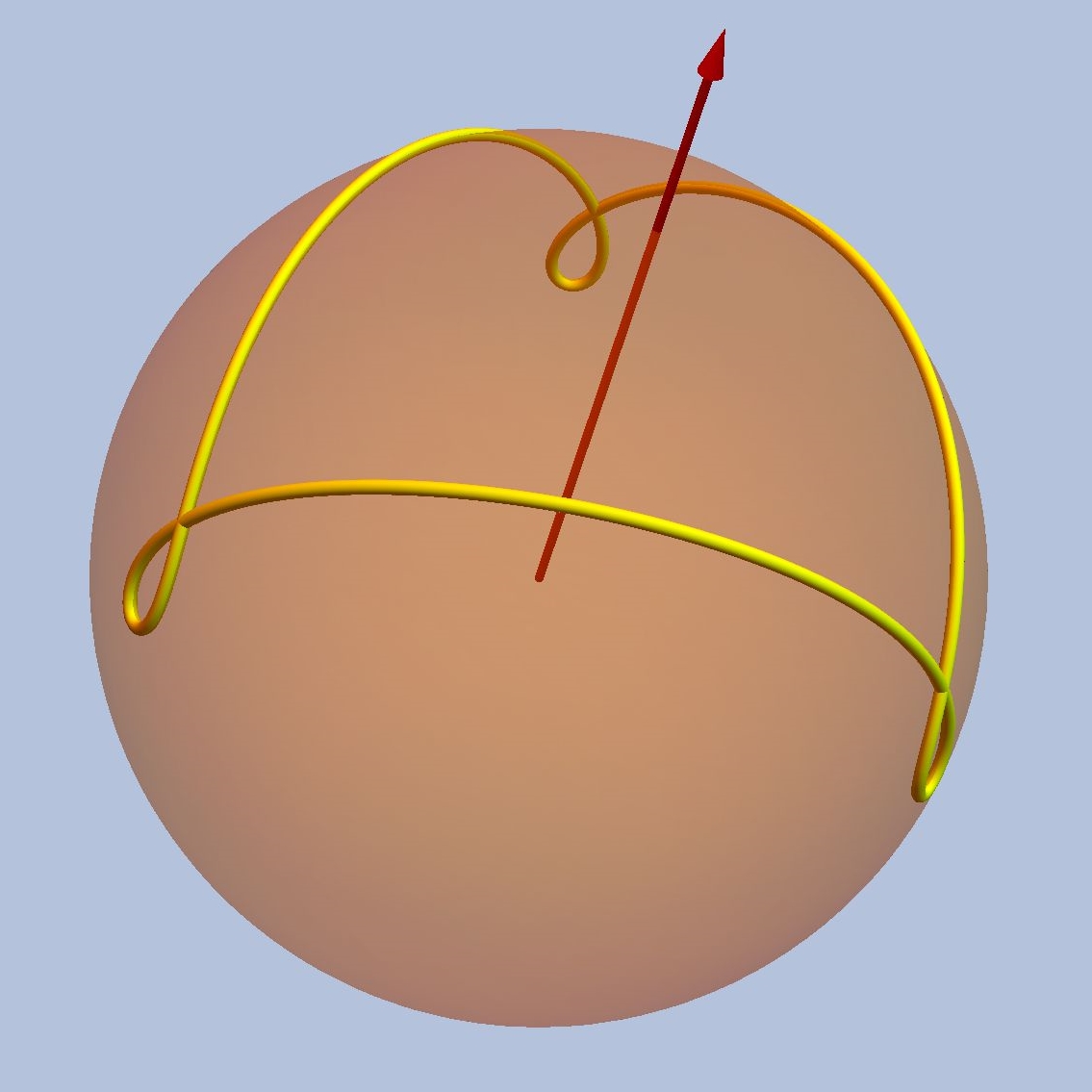}
\caption{$(-1.1,1,3)$}
\end{subfigure}
\quad
\begin{subfigure}[b]{0.3\linewidth}
\includegraphics[height=5cm]{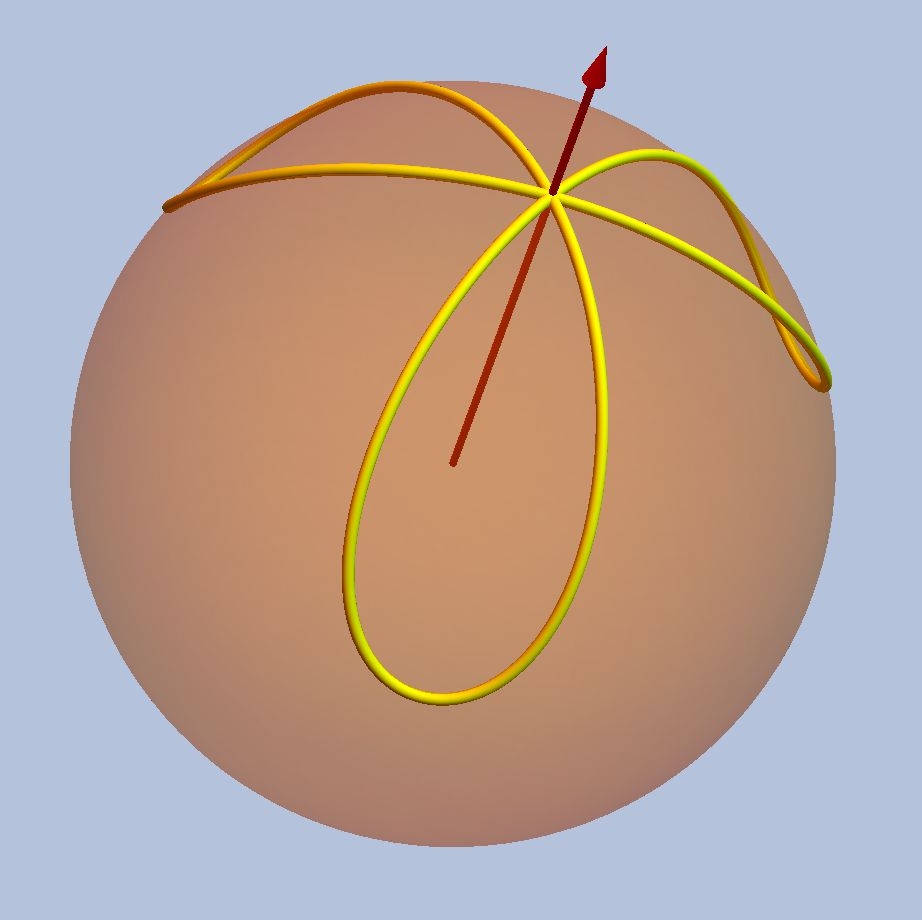}
\caption{$(-0.27,2,3)$}
\end{subfigure}
\quad
\begin{subfigure}[b]{0.3\linewidth}
\includegraphics[height=5cm]{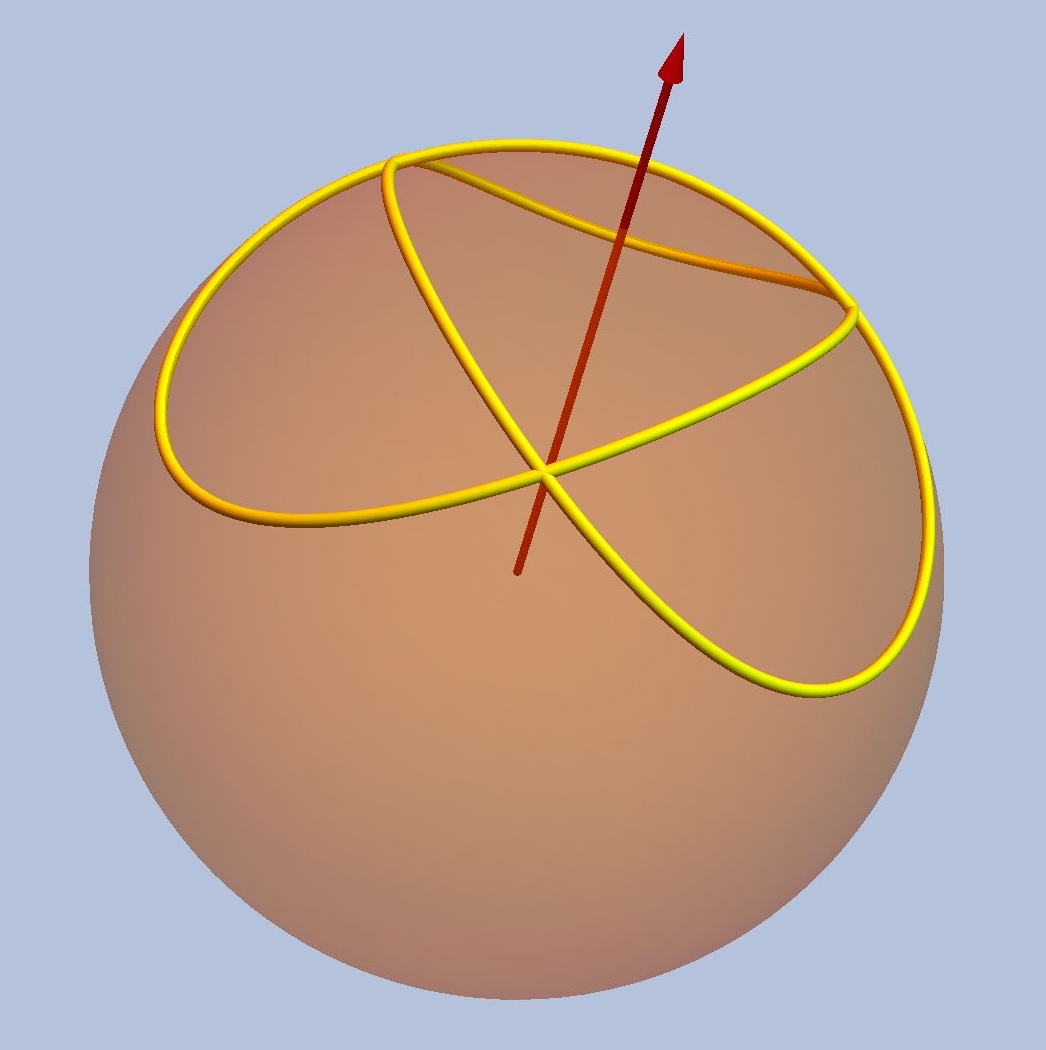}
\caption{$(0.1,1,3)$}
\end{subfigure}
\caption{Three B-strings with three-fold symmetry with different multipliers $\lambda$ and characteristic numbers $q=m/n$. From left to right: $4\lambda\xi+1<0$, $4\lambda\xi+1=0$ and $4\lambda\xi+1>0$. For each of them we show the corresponding parameters $(\lambda,m,n)$.}
\label{B-strings}
\end{figure}

In Figure \ref{B-strings} we show three B-strings with three-fold symmetry for different values of the Lagrange multiplier $\lambda$ and for suitable characteristic numbers $q=m/n$. These cases cover all the possible options for the sign of $4\lambda\xi+1$ discussed in Theorem \ref{B}. More examples will be discussed in detail in Section 5.

The material of this paper is organized into four sections and two appendices. In Section 2 we write the Euler-Lagrange equation and its associated conservation law. Subsequently, we consider the monodromy map and we formulate the closure condition  in terms of the monodromy. In Section 3 we prove the three main theorems (which are stated in a more technical form than the one presented in this introduction).  In Section 4 we focus on the theoretical aspects, analyzing the phase space (\cite{GM,G,Hs})  and the momentum map for the Hamiltonian $SO(3)$-action. Beside its theoretical relevance, the Hamiltonian approach provides a general framework for concretely implement the integration by quadratures. Finally, Section 5 is devoted to the discussion of explicit examples which illustrate the theoretical properties shown in previous sections. In the first part of the Appendix A we construct $\mu$ inverting a complete elliptic integral of the third kind and we compute its least period. In the second part we write $\Psi_\lambda$ in terms of complete elliptic integrals and we prove two relevant limits that have been used in the proof of Theorem \ref{C}. In Appendix B, we briefly discuss the case of $1/2$-elasticae in the plane.

The graphics, symbolic computations and numerical evaluations have been performed with the software \emph{Mathematica 13}. To evaluate the least period of $\mu$ we used its explicit expression given in Appendix A and the library for elliptic integrals implemented in Mathematica, while to evaluate $\mu$ itself we solved numerically (\ref{EL}) on the interval $[0,n\omega]$, where $n$ is the wave number and $\omega$ is the least period of $\mu$, with initial conditions $\mu(0)=e_2$ and $\dot{\mu}(0)=0$. Here, $e_2$ is the lowest positive real root of the polynomial defined in Proposition \ref{conlaw}. We then solved numerically the linear system (\ref{frenet}) on the interval $[0,n\omega]$ to compute the critical curves. The initial conditions are those specified in (\ref{frenetinitial}). This few comments should give a rough idea on how the plots for the graphics were computed. 

Complete elliptic integrals involving the square root of a quartic polynomial do appear several times. We used the monograph \cite{BF} as our reference on this technical topic. The bibliographic references, while relatively consistent, do not completely reflect the vast literature devoted to functionals depending on curvatures and their interrelations with analysis, symplectic geometry and applied mathematics. It is the result of a very partial selection, aimed at the specific themes considered in this work.

\section{Critical Curves}

Let $(x,y,z)$ be the standard coordinates of the Euclidean space $\mathbb{R}^3$ and ${\mathbb S}^2$ be the $2$-sphere of radius one centered at the origin endowed with the induced metric of constant curvature $c=1$.

Let $\gamma:I\subseteq\mathbb{R}\longrightarrow\mathbb{S}^2$ be a smooth immersed curve parameterized by the arc-length $s\in I$. Denote by $T(s):=\dot{\gamma}(s)$ the unit tangent vector field along the curve $\gamma(s)$, where the upper dot represents the derivative with respect to the arc-length parameter, and define the unit normal vector field $N(s)$ along $\gamma(s)$ to be the counter-clockwise rotation of $T(s)$ through an angle $\pi/2$ in the tangent bundle of $\mathbb{S}^2$. In this setting, the (signed) geodesic \emph{curvature} $\kappa(s)$ of $\gamma(s)$ is defined by the Frenet-Serret equation
$$\nabla_T T(s)=\kappa(s)N(s)\,,$$
where $\nabla$ denotes the Levi-Civita connection on $\mathbb{S}^2$.  We will say that a curve is \emph{convex} if $\kappa(s)>0$ for all $s\in I$. For convex curves we introduce the following geometric invariant, referred as the $\mu$-invariant,
$$\mu(s):=\sqrt{\kappa(s)}\,.$$
It follows from the Fundamental Theorem for Spherical Curves that $\mu(s)$ completely determines the (convex) curve, up to rigid motions.

Let $\mathcal{C}^4(\mathbb{R},\mathbb{S}^2)$ be the space of  immersed  curves $\gamma:\mathbb{R}\longrightarrow\mathbb{S}^2$ of class $\mathcal{C}^4$ parameterized by the arc-length $s\in \mathbb{R}$ and let $\lambda\in\mathbb{R}$ be a constant. The \emph{$1/2$-Bernoulli's bending functional} with Lagrange multiplier $\lambda$ is defined by
\begin{equation}\label{energy}
\mathcal{B}_\lambda:\gamma\in\mathcal{C}^4(\mathbb{R},\mathbb{S}^2)\longmapsto\int_0^{L_\gamma}\left(\sqrt{|\kappa|}+\lambda\right)ds\,,
\end{equation}
where $L_\gamma$ stands for the length of $\gamma$.

Using a standard formula for the variational derivative of functionals depending on the curvature $\kappa$ (see for instance \cite{AGM,AGP,LS2,P}) we obtain that the $\mu$-invariant of a {\em convex} critical curve with respect to compactly supported smooth variations satisfies the \emph{Euler-Lagrange equation} 
\begin{equation}\label{EL}
\frac{d^2}{ds^2}\left(\frac{1}{\mu}\right)-\frac{1}{\mu}\left(\mu^4-1\right)-2\lambda\mu^2=0\,.
\end{equation}

\begin{rem}\label{circles} For every $\lambda\in\mathbb{R}$, with the exception of geodesics, there exists a unique circle critical for $\mathcal{B}_\lambda$. Its constant $\mu$-invariant is the (unique) positive solution of
\begin{equation}\label{cst}
\mu^4+2\lambda\mu^3-1=0\,,
\end{equation}
which we denote by $\eta_\lambda$. Consequently, $\eta:\lambda\in\mathbb{R}\longmapsto \eta_\lambda\in\mathbb{R}^+$ is a real-analytic function.
\end{rem}

From now on we assume that $\mu(s)$ is non-constant. Then \eqref{EL} admits a first integral from standard arguments (see \cite{AGM,AGP,LS2,P} again for details). We describe this conservation law in the following result.

\begin{prop}\label{conlaw} Let $\mu(s)$ be a non-constant positive solution of \eqref{EL}.  Then $\mu(s)$ satisfies the first order differential equation
\begin{equation}\label{ODE}
\dot{\mu}^2(s)=-\mu^2(s)\, Q\left(\mu(s)\right),
\end{equation}
where $Q$ is the quartic polynomial defined by 
\begin{equation}\label{polynomial}
Q(t):=t^4+4\lambda t^3+4\left(\lambda^2-\xi^2\right)t^2+1\,,
\end{equation}
and $\xi>0$ is a constant of integration.
\end{prop}

From Proposition \ref{conlaw} it follows that, for fixed $\lambda\in\mathbb{R}$, solutions of the Euler-Lagrange equation \eqref{EL} belong to a two parameter family. Nevertheless, by translating the origin of the arc-length parameter $s$ if necessary, we can assume that non-constant solutions $\mu(s)$ belong to a one parameter family, since the constant of integration arising from integrating \eqref{ODE} may be assumed to be zero. Consequently, the one parameter family of solutions depends on the constant of integration $\xi>0$, whose physical meaning will be clarified in Section 4.

\begin{defn} A \emph{B-curve} (with multiplier $\lambda$) is an arc-length parameterized convex curve $\gamma:\mathbb{R}\longrightarrow\mathbb{S}^2$ with non-constant positive periodic curvature satisfying \eqref{ODE}. For brevity, a periodic B-curve is said a {\em B-string}.
Two B-curves $\gamma$ and $\widetilde{\gamma}$ are said to be \emph{equivalent} if there exists $A\in SO(3)$ and $a\in\mathbb{R}$ such that $\widetilde{\gamma}(s)=A\cdot\gamma(s+a)$, i.e., if there exists an isometry transforming one into another and a translation of the arc-length parameter. We denote the equivalence class of $\gamma$ by $\left[\gamma\right]$ and the set of the equivalence classes of B-curves with multiplier $\lambda$ by $\mathcal{M}_\lambda$.
\end{defn}

In the following result we show that, with the exclusion of geodesics, $1/2$-elasticae do not have inflection points. Thus, possibly reversing the orientation, $\kappa>0$. As a consequence, non-trivial $1/2$-elasticae are B-curves whose $\mu$-invariant is a non-constant periodic solution of \eqref{ODE}.

\begin{prop}\label{periodic} Let $\gamma:I\subseteq\mathbb{R}\longrightarrow\mathbb{S}^2$  be an arc-length parameterized curve with non-constant curvature $\kappa$ and assume that $I$ is its maximal domain of definition. If $\gamma$ is a critical curve$^5$\footnote{$^5$ In the sense specified in the Introduction.} for ${\mathcal B}_{\lambda}$ with respect to compactly supported smooth variations, then $\kappa(s)\neq 0$ for every $s\in I$ (possibly reversing the orientation we may assume $\kappa>0$). Moreover, $I={\mathbb R}$ and $\gamma$ is a B-curve. Conversely, every B-curve is a $1/2$-elastic curve.
\end{prop}
\textit{Proof.} Without loss of generality we may assume that $\kappa(s)>0$ for some $s\in I$. By contradiction, suppose that $\kappa$ vanishes at some point. Then there exist $s_*$ such that $\kappa(s_*)=0$ and $\kappa > 0$ on $(s_*,s_*+\epsilon)$, or $\kappa > 0$ on $(s_*-\epsilon,s_*)$. Consequently, the $\mu$-invariant of $\gamma$ restricted to this interval is a positive solution of (\ref{ODE}). This implies that the polynomial \eqref{polynomial} must be negative for some $t>0$.  We first observe that the limit when $t\to\infty$ of $Q(t)$ is $\infty$, while $Q(0)=1>0$. Regardless of the values of $\lambda$ and $\lambda^2-\xi^2$, $Q(t)$ has either zero or two changes of signs among its coefficients. It then follows from Descartes' rule of signs that $Q(t)$ has either zero, one (double) or two (distinct) positive roots. The case of zero positive roots can be discarded since from above limits one would conclude that $Q(t)>0$ for all $t>0$. The case of the double root corresponds to a circle, which has been considered in Remark \ref{circles}.  On the other hand, since $\kappa$ is not constant, this case is excluded.  Therefore, it only remains the case of two distinct positive roots.  In this case  the algebraic curve $y^2=-x^2 Q(x)$ has a unique $1$-dimensional, closed connected component $\mathcal{C}^*$ contained in the half-plane $x>0$ (see Figures \ref{Example1-1}-\ref{Example5-1}). In addition,  ${\mathcal C}^*$ is smooth and intersects the $y=0$ axis at two distinct points $(e_2,0)$, $(e_1,0)$, where $0<e_2<e_1$ are the simple positive roots of $Q$. Then, $0<e_2^2\le \kappa(s)\le e_1^2$, for every $s\in (s_*,s_*+\epsilon)$, or $s\in(s_*-\epsilon,s_*)$.  A contradiction. This proves the first part of the statement.

\begin{figure}[h]
	\begin{center}
		\includegraphics[height=5cm]{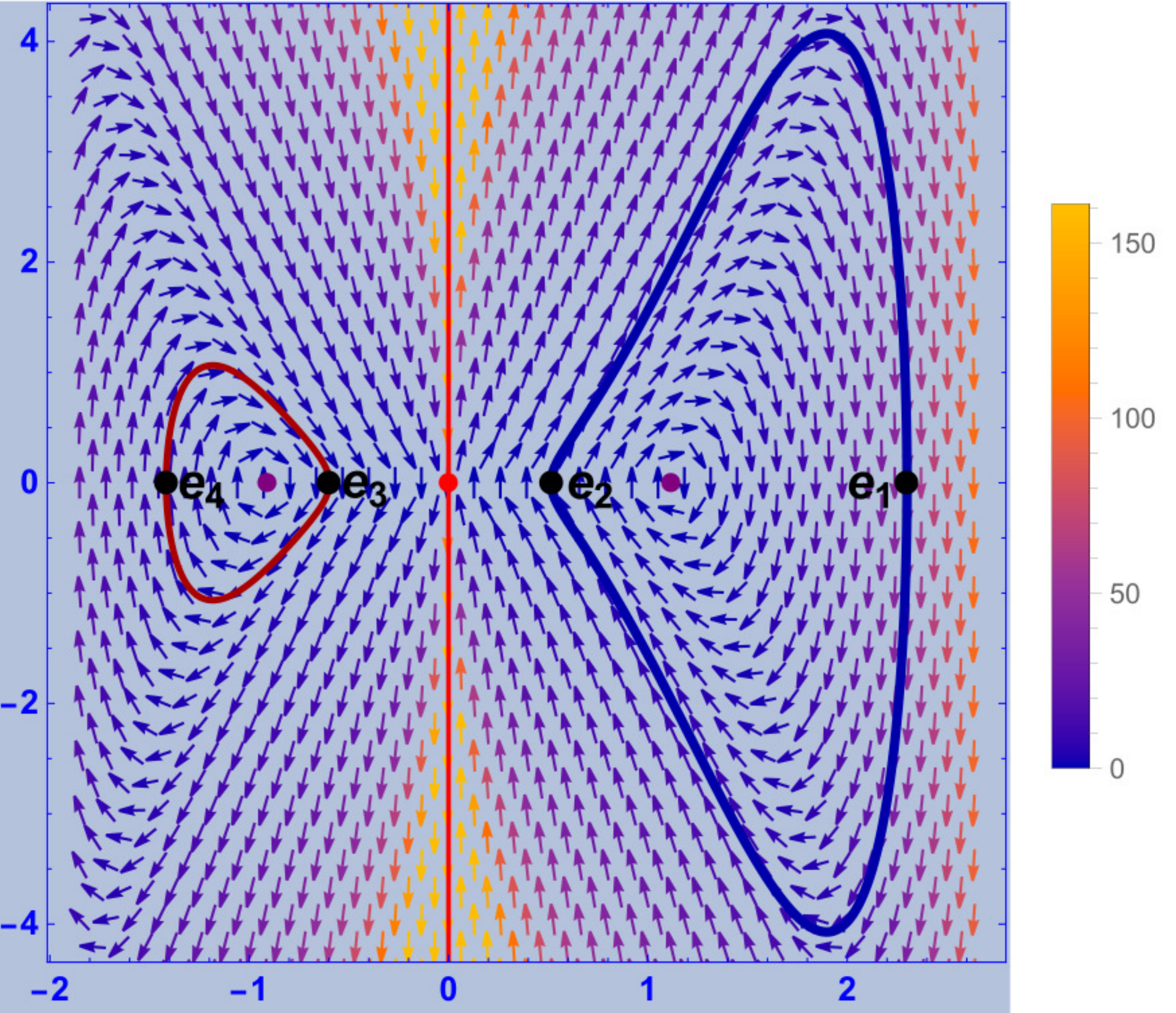}\quad\quad\quad
		\includegraphics[height=5cm]{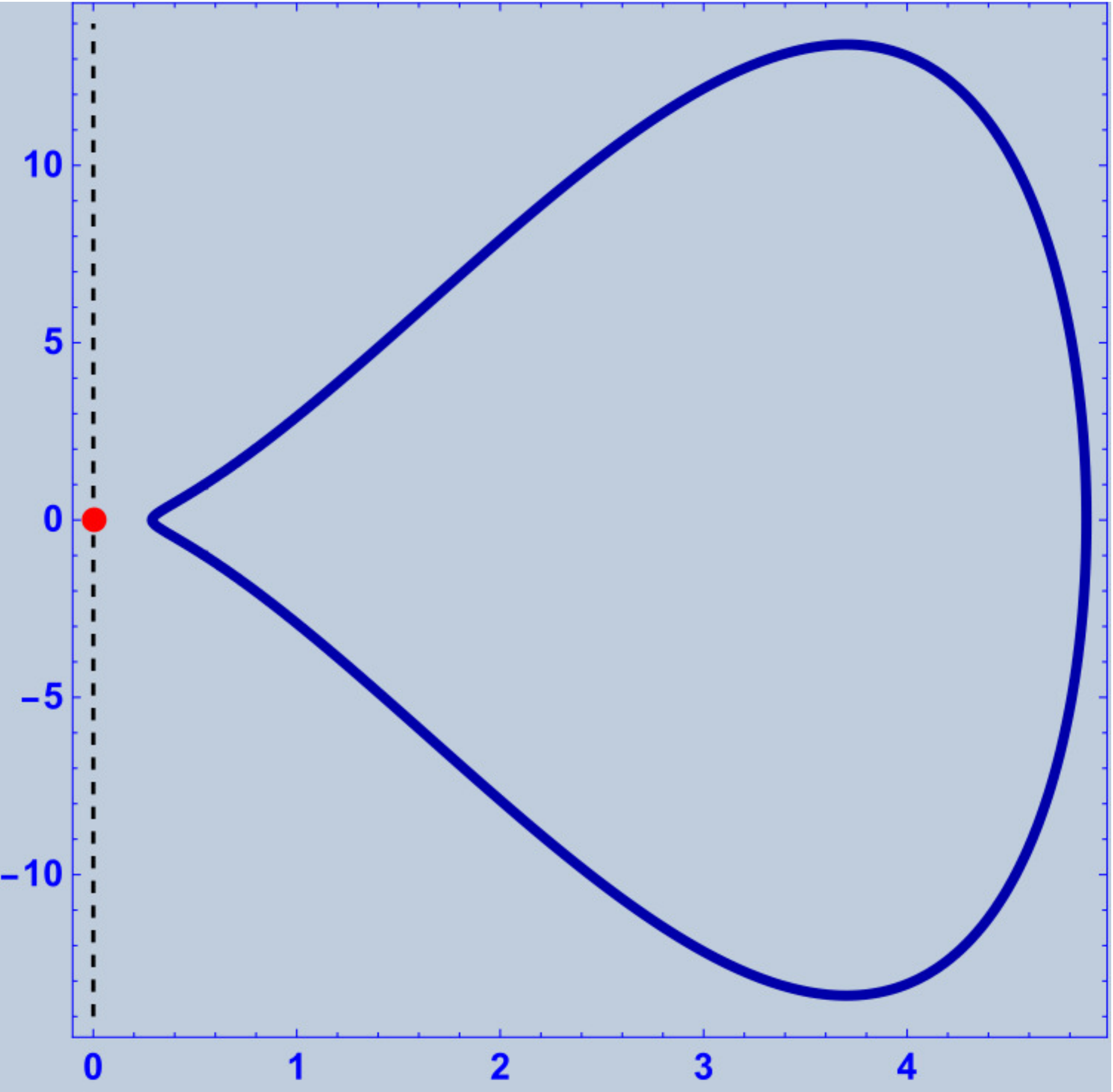}
		\caption{On the left: the plot of the vector field  $\vec{X}_{\lambda}$, $\lambda = -0.2$, and the phase curve with $\lambda=-0.2$ and  $\xi^2= 0.869759$. The curve colored in blue is the reduced phase curve contained in the half-plane ${\mathbb H}=\{(x,y)\,\lvert\, x>0\}$,  parameterized by $(\mu,\dot{\mu})$, where $\mu$ is a positive solution of (\ref{ODE}).  On the right: the corresponding  ``standard" phase portrait, which consists of the origin (singular point) and the image of the blue curve on the left by the contact transformation ${\mathtt f}:(x,y)\in {\mathbb H}\to (x^2,  2xy)\in {\mathbb H}$.}\label{php}
	\end{center}
\end{figure}

Now we prove that $\mu=\sqrt{\kappa}$ is defined on the whole real axis and that it is a periodic function. Let ${\mathbb H}$ be the half-plane $\{(x,y)\,\lvert\, x>0\}$ and $\vec{X}_{\lambda}$ be the vector field defined on ${\mathbb H}$ by
$$\vec{X}_{\lambda}\lvert_{(x,y)}=y\,\partial_x+\frac{1}{x}\left(2y^2+x^2-2\lambda x^5-x^6\right)\partial_y\,.$$
See Figure \ref{php} for an illustration of this vector field on $\mathbb{H}$ and its natural extension to ${\mathbb R}^2$ minus the $Oy$-axis. Define $\psi$ by $\psi(x,y)=y^2+x^2Q(x)$. Then
\begin{equation}\label{XL}
\vec{X}_{\lambda}\lvert_{(x,y)}=\frac{1}{2}\left(\partial_y \psi\, \partial_x - \partial_x\psi\, \partial_y\right)\lvert_{(x,y)}\,,\quad\quad\quad \forall (x,y)\in {\mathcal C}^*.
\end{equation}
From the Euler-Lagrange equation it follows that $\zeta=(\mu,\dot{\mu})$ is an integral curve of $\vec{X}_{\lambda}$.  Moreover, the trajectory of $\zeta$ is contained in ${\mathcal C}^*$, because $\mu$ is a positive solution of  \eqref{ODE}.  Since ${\mathcal C}^*$ is compact there is a compactly supported smooth function 
$\varrho : {\mathbb H}\to {\mathbb R}$ such that $\varrho\lvert_{{\mathcal C}^*}=1$.  Then, $\zeta$ is an integral curve of the complete vector field $\varrho \vec{X}_{\lambda}$. This proves that ${\mathbb R}$ is the maximal interval of definition of the function $\mu$. From (\ref{XL}) we have 
$\zeta({\mathbb R})= {\mathcal C}^*$ and $\varrho \vec{X}_{\lambda}\lvert_{(x,y)}\neq \vec{0}$,   for every $(x,y)\in {\mathcal C}^*$. Taking into account that ${\mathcal C}^*$ is a compact embedded curve and using the Poincar\'e-Bendixson Theorem we may conclude that $\zeta$ is periodic. A fortiori, we also conclude that $\mu$ is periodic. Consequently, also the maximal domain of definition of the associated curve is ${\mathbb R}$.  This shows that $\gamma$ is a B-curve. Since the last assertion is trivial, this finishes the proof.
\hfill$\square$

\begin{rem} Note that only the reduced phase curves lying in the half-plane ${\mathbb H}$ are relevant in the study of B-strings. We emphasize here that the phase curves differ from the standard phase portraits (\cite{G}). In fact, putting $\mu=\sqrt{\kappa}$ in (\ref{ODE}) we obtain
$$\dot{\kappa}^2+4\kappa^2\left(1+4[\lambda^2-\xi^2]\kappa+4\lambda\kappa^{3/2}\right)=0\,.$$
Therefore, the phase portraits, in the standard sense of the term,  are the (non-algebraic) curves defined by equations of the following type
$$y^2+4x^2\left(1+4[\lambda^2-\xi^2]x+4\lambda x^{3/2}\right)=0\,.$$
The reduced phase curves lying in the half-plane $x>0$ are the pre-images of the standard phase portraits under the contact transformation 
${\mathtt f}:(x,y)\in {\mathbb H}\to (x^2,  2xy)\in {\mathbb H}$.
\end{rem}

\begin{figure}[h]
\begin{center}
\includegraphics[height=5cm]{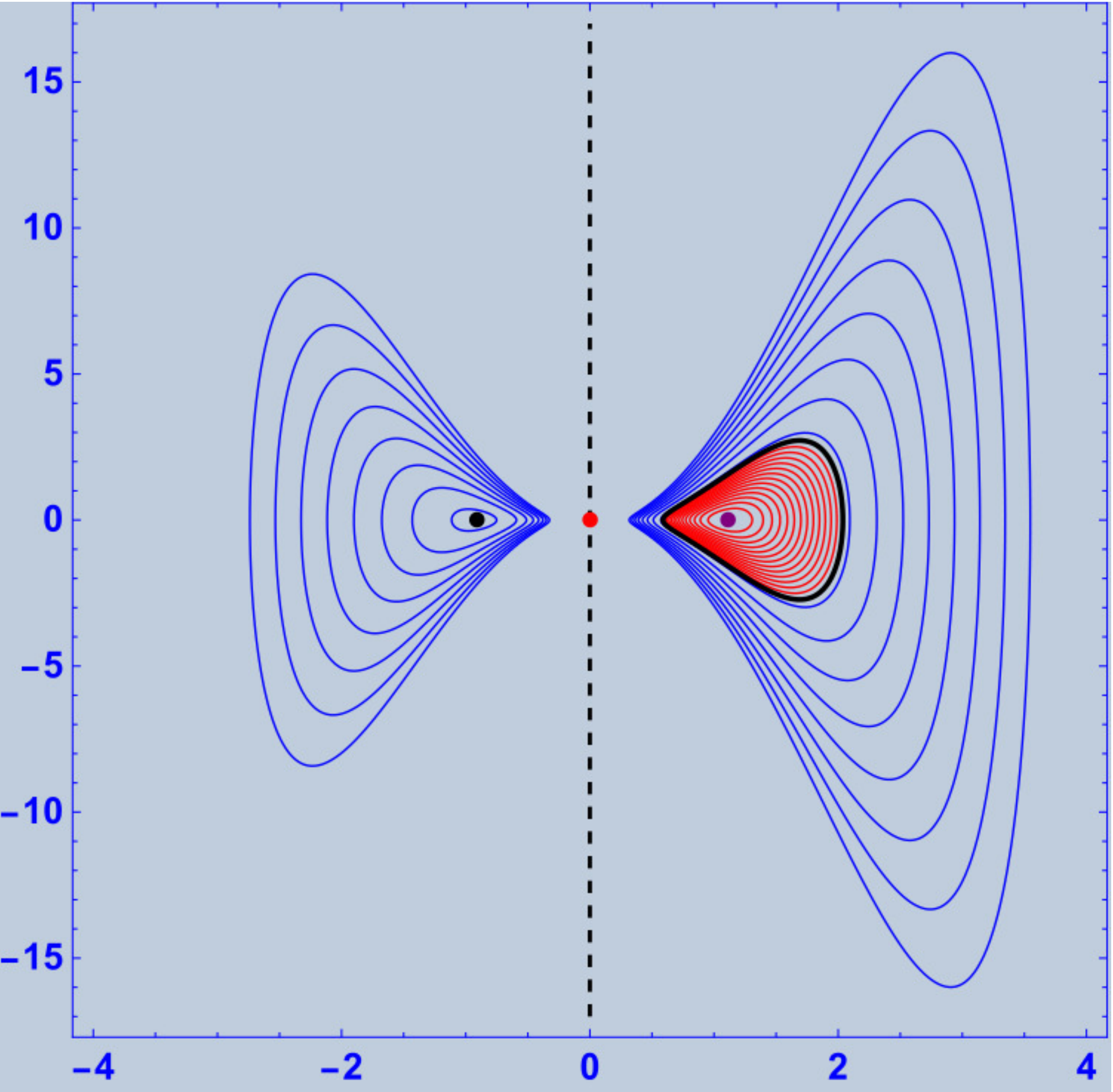}\quad\quad\quad
\includegraphics[height=5cm]{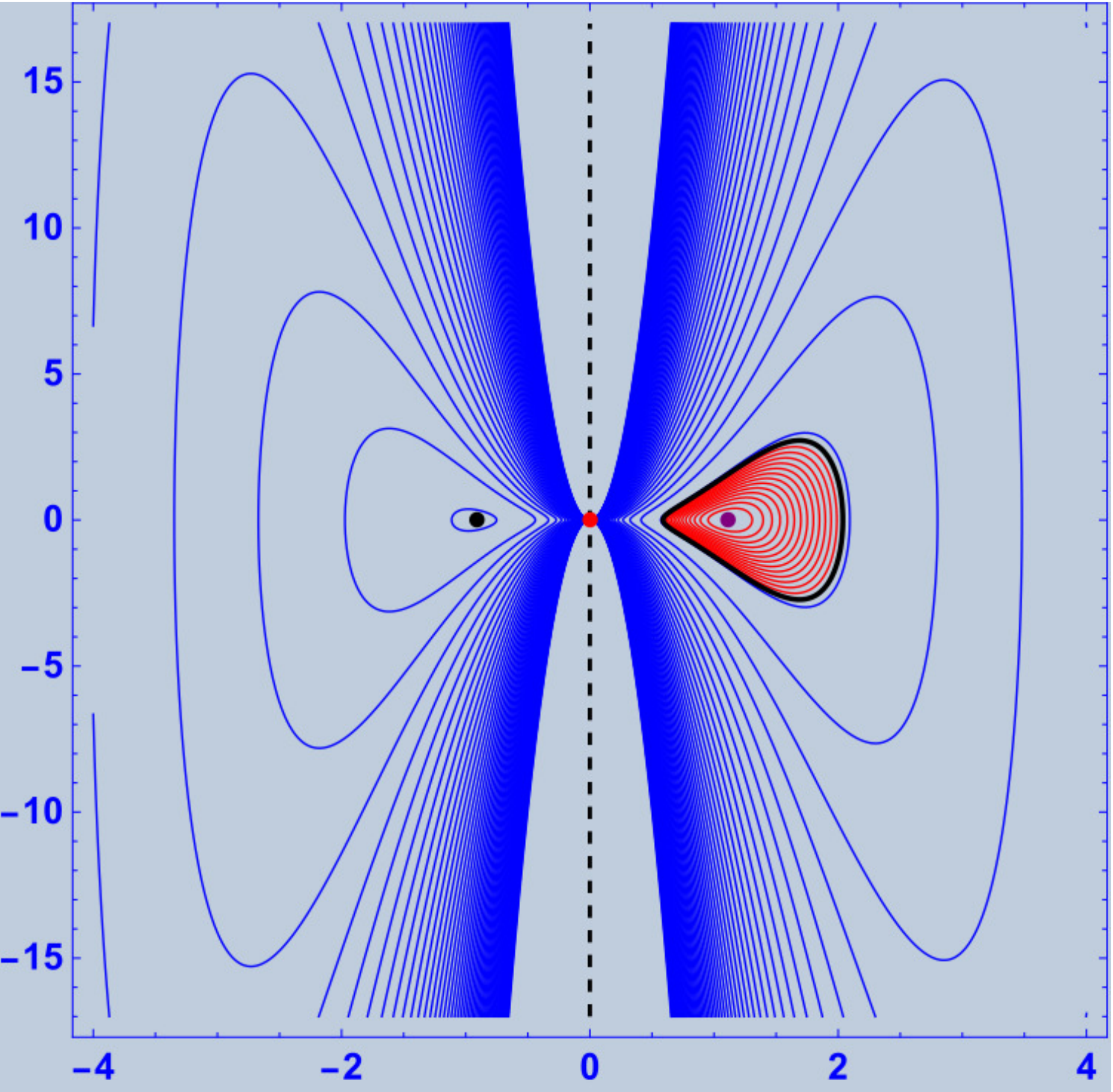}
\caption{The phase curves for $\lambda=-0.2$. On the left are depicted the phase curves when $\xi^2$ varies in the interval $[0,1.675]$. On the right are depicted the phase curves when $\xi^2$ varies in the interval $[0,20]$.}\label{php2}
\end{center}
\end{figure}

\begin{ex}The picture on the left of Figure \ref{php} depicts the plot of the vector field $\vec{X}_{\lambda}$, $\lambda = -0.2$,  thought of as a vector field defined on ${\mathbb R}^2$ minus the $Oy$-axis.  The plot reproduces only the direction of the vector field and the legend specifies the magnitude. The blue and red curves are the reduced phase curves contained in the negative and positive half-planes $x<0$ (dark red) and $x>0$ (blue). The value of the constant $\xi^2$ is $ 0.869759$. The origin (the red point) is the isolated singular point of the phase curve.  The two points colored in purple are the two zeroes of $\vec{X}_{\lambda}$.   The picture on the right of Figure \ref{php} reproduces the corresponding  ``standard" phase portrait.  Figure \ref{php2} reproduces phase curves for $\lambda=-0.2$ and $\xi^2\in [0,1.675]$ (left) and $\xi^2\in [0,20]$ (right). When $0\le \xi^2<0.328893$, the phase curve consists only of a singular point, the origin (in red). When $\xi^2=0.328893$ the phase curve consists of two singular points, the origin and $(1.117181339509767,0)$, one of the two zeroes of the vector field $\vec{X}_{\lambda}$ (in purple). These phase curves are ``virtual" in the sense that are not originated from B-strings. When  $0.328893<\xi^2<0.730907$ the phase curve  has an isolated singular point (the origin) and a smooth reduced phase curve contained in the half-plane $x>0$ (colored in red). When $\xi^2=0.730907$, the phase curve  has two singular points, the origin and the other zero  $(-0.9131837363949253, 0)$ of the vector field $\vec{X}_{\lambda}$ (colored in black) and a smooth connected component (colored in black) contained in the half-plane $x>0$. When  $\xi^2>0.730907$ the phase curve has an isolated singular point (the origin) and two smooth connected components (colored in blue). One contained in the half-plane $x<0$ and the other one in the half-plane $x>0$.
\end{ex}

\begin{rem} Proposition \ref{periodic} proves that, up to equivalence, there exists a one parameter family of B-curves, depending on the constant of integration $\xi>0$. It turns out that this constant of integration $\xi>0$ can be described in terms of the Lagrange multiplier $\lambda$ and the largest root of the polynomial $Q$. Indeed, as shown above, this polynomial has two positive roots, which we denote by $e_1>e_2>0$. Then, it follows from $Q(e_1)=0$ that
\begin{equation}\label{xi}
\xi=\frac{1}{2e_1}\sqrt{1+\left(e_1\left[e_1+2\lambda\right]\right)^2\,}\,.
\end{equation}
From this relation we can assume that B-curves depend on $\lambda$ and $e_1$ and we can consider these parameters as the fundamental ones. In addition, the $\mu$-invariant of a critical curve can be built inverting a function involving incomplete elliptic integrals of the third kind, the Jacobi's amplitude and the Jacobi's sn-function, with parameters that depend on  $\lambda$ and $e_1$.  Consequently,  also the least period of $\mu$ can be written in terms of complete elliptic integrals of the first and third kind (we refer to the first part of Appendix A for the explicit formulae). 
\end{rem}

The set of equivalence classes of B-curves is in one-to-one correspondence with the plane domain $\mathcal{P}=\{(\lambda,e_1)\in\mathbb{R}^2\,\lvert\,e_1>\eta_\lambda\}$, where $\eta_\lambda$ was defined in Remark \ref{circles} as the only positive solution of $\mu^4+2\lambda\mu^3-1=0$.

\begin{thm} For every $\lambda\in\mathbb{R}$, the map $e_1\in\left(\eta_\lambda,\infty\right)\longmapsto\left[\gamma_{\lambda,e_1}\right]\in\mathcal{M}_\lambda$ is bijective.
\end{thm}
\textit{Proof.} The proof is a straightforward consequence of the existence and uniqueness of solutions for ordinary differential equations. For each $\lambda\in\mathbb{R}$ and $e_1\in(\eta_\lambda,\infty)$ fixed, we get a unique $\xi>0$ from \eqref{xi}. Note that for $\lambda\in\mathbb{R}$ fixed the relation \eqref{xi} between $\xi$ and $e_1$ is bijective. Then, there exists a unique, up to translation of the arc-length parameter, solution of \eqref{ODE}, denoted by $\mu_{\lambda,e_1}(s)$. From the Fundamental Theorem for Spherical Curves, this geometric invariant uniquely determines a convex spherical curve $\gamma_{\lambda,e_1}$, up to rigid motions. In conclusion, the curve $\gamma_{\lambda,e_1}$ is unique, up to equivalence. \hfill$\square$
\\

In order to treat the B-curves, it is convenient to fix a suitable representative for each equivalence class. Let $(\lambda,e_1)\in\mathcal{P}$ and $\mu\equiv\mu_{\lambda,e_1}$ be the unique solution of \eqref{ODE} such that $\mu_{\lambda,e_1}(\omega_{\lambda,e_1}/2)=e_1$, where $\omega_{\lambda,e_1}$ denotes the least period of $\mu_{\lambda,e_1}$. Consider $\gamma_{\lambda,e_1}$ to be the unique B-curve with curvature $\mu_{\lambda,e_1}^2$ such that
\begin{equation}\label{frenetinitial}
\left\{\begin{split} &\gamma_{\lambda,e_1}\left(\frac{\omega_{\lambda,e_1}}{2}\right)=\left(\frac{1}{\sqrt{1+\left(e_1\left[e_1+2\lambda\right]\right)^2}},\frac{-e_1\left(e_1+2\lambda\right)}{\sqrt{1+\left(e_1\left[e_1+2\lambda\right]\right)^2}},0\right) \\
& \dot{\gamma}_{\lambda,e_1}\left(\frac{\omega_{\lambda,e_1}}{2}\right)=(0,0,-1) \end{split}\right.\,.
\end{equation}
We call $\gamma_{\lambda,e_1}$, obtained as above, the \emph{standard B-curve} with parameters $\lambda$ and $e_1$.

\begin{rem} It is not restrictive to focus exclusively on standard B-curves and, hence, from now on we implicitly assume that the B-curves in consideration are in their standard form.
\end{rem}

We finish this section by formulating the closure condition for a B-curve in terms of the monodromy map. For a B-curve $\gamma_{\lambda,e_1}$ we define the (spherical) \emph{Frenet frame field} along $\gamma_{\lambda,e_1}$ as the map $\mathcal{F}_{\lambda,e_1}\equiv \left(\gamma_{\lambda,e_1},\dot{\gamma}_{\lambda,e_1},\gamma_{\lambda,e_1}\times\dot{\gamma}_{\lambda,e_1}\right):\mathbb{R}\longrightarrow SO(3)$ where $\times$ denotes the usual vector cross product of $\mathbb{R}^3$. The map $\mathfrak{m}:\left(\lambda,e_1\right)\in\mathcal{P}\longmapsto\mathcal{F}_{\lambda,e_1}(\omega_{\lambda,e_1})\cdot\left[\mathcal{F}_{\lambda,e_1}(0)\right]^{-1}\in SO(3)$ is called the \emph{monodromy}. 

\begin{thm} The monodromy $\mathfrak{m}$ is a continuous function of $\lambda$ and $e_1$. Moreover, $\gamma_{\lambda,e_1}$ is a B-string if and only if $\mathfrak{m}_{\lambda,e_1}$ has finite order.
\end{thm}
\textit{Proof.} Let $(\lambda,e_1)\in\mathcal{P}$ and consider a B-curve $\gamma_{\lambda,e_1}$ (in its standard form). Since $\mu_{\lambda,e_1}$ is a solution of \eqref{ODE}, it also satisfies \eqref{EL}, together with the initial conditions $\mu_{\lambda,e_1}(\omega_{\lambda,e_1}/2)=e_1$ and $\dot{\mu}_{\lambda,e_1}(\omega_{\lambda,e_1}/2)=0$. Therefore, $\mu_{\lambda,e_1}(s)$ is a real-analytic function of $s\in\mathbb{R}$ and $(\lambda,e_1)\in\mathcal{P}$, and so is $\kappa_{\lambda,e_1}(s)=\sqrt{\mu_{\lambda,e_1}(s)}$. On the other hand, $\mathcal{F}_{\lambda,e_1}$ satisfies
\begin{equation}\label{frenet}
\left\{\begin{split} &\dot{\mathcal{F}}_{\lambda,e_1}=\mathcal{F}_{\lambda,e_1}\cdot\mathcal{K}_{\lambda,e_1}\\
&\mathcal{F}_{\lambda,e_1}(\omega_{\lambda,e_1}/2)=F_{\lambda,e_1}\end{split}\right.\,,
\end{equation}
where
\begin{equation*}
\mathcal{K}_{\lambda,e_1}=\begin{pmatrix} 0 & -1 & 0 \\ 1 & 0 & -\kappa_{\lambda,e_1} \\ 0 & \kappa_{\lambda,e_1} & 0 \end{pmatrix}
\end{equation*}
and $F_{\lambda,e_1}=\left([E_1]_{\lambda,e_1},[E_2]_{\lambda,e_1},[E_3]_{\lambda,e_1}\right)$ for
\begin{equation*}
\left\{\begin{split} &(E_1)_{\lambda,e_1}=\left(\frac{1}{\sqrt{1+(e_1[e_1+2\lambda])^2}},\frac{-e_1(e_1+2\lambda)}{\sqrt{1+(e_1[e_1+2\lambda])^2}},0\right)\\
&(E_2)_{\lambda,e_1}=\left(0,0,-1\right)\\
&(E_3)_{\lambda,e_1}=(E_1)_{\lambda,e_1}\times (E_2)_{\lambda,e_1}\end{split}\right.\,.
\end{equation*}
Consequently, $\mathcal{F}_{\lambda,e_1}(s)$ is also a real-analytic function of $s\in\mathbb{R}$ and $(\lambda,e_1)\in\mathcal{P}$.

The least period of $\omega_{\lambda,e_1}$ can be written in terms of complete elliptic integrals of the first and third kind (see Part I of Appendix A). It then follows from the properties of these integrals that $\omega:\left(\lambda,e_1\right)\in\mathcal{P}\longmapsto \omega_{\lambda,e_1}\in\mathbb{R}$ is a continuous function which is real-analytic on $\widehat{\mathcal{P}}$ i.e., the complement of the zero locus of the real-analytic function $(\lambda,e_1)\in\mathcal{P}\longmapsto \left(e_1+[e_2]_{\lambda,e_1}\right)^2-4e_1^3[e_2]_{\lambda,e_1}^3$. Thus, $\mathfrak{m}$ is also a continuous function and real-analytic on $\widehat{\mathcal{P}}$.

Finally, since $\mu_{\lambda,e_1}$ is periodic with least period $\omega_{\lambda,e_1}$, we conclude from \eqref{EL} that for every $k\in\mathbb{Z}$,
$$\mathcal{F}_{\lambda,e_1}(s+k\omega_{\lambda,e_1})=\mathfrak{m}_{\lambda,e_1}^k\cdot\mathcal{F}_{\lambda,e_1}(s)\,.$$
This finishes the proof. \hfill$\square$

\begin{rem} Observe that the order of the monodromy $\mathfrak{m}_{\lambda,e_1}$ is, precisely, the wave number $n$ of $\gamma_{\lambda,e_1}$.
\end{rem}

\section{Integrability by Quadratures and Existence of B-Strings}

In this section we give the parameterization of B-curves in terms of just one quadrature and prove the main theorems of the paper. The parameterization we will obtain in this section is, essentially, the one of Theorem \ref{A}. Observe that, locally, the ordinary differential equation determining $\mu(s)$ can be used to make a change of variable in $\theta(s)$ and so this parameterization depends on just one quadrature.

We begin by defining a curve in the plane domain $\mathcal{P}=\{\left(\lambda,e_1\right)\in\mathbb{R}^2\,\lvert\,e_1>\eta_\lambda\}$. The \emph{exceptional locus} is the smooth curve $\mathcal{P}_*\subset\mathcal{P}$ defined by the equation
\begin{equation}\label{cubic}
4\lambda^2e_1^3+8\lambda^3e_1^2-e_1+2\lambda=0\,.
\end{equation}
It is easy to check that this exceptional curve is contained in $\{\left(\lambda,e_1\right)\in\mathcal{P}\,\lvert\,\lambda<0\}$. We say that the parameters $(\lambda,e_1)$ are \emph{exceptional} if they belong to $\mathcal{P}_*$. Moreover, for any $\lambda<0$, the cubic equation \eqref{cubic} has a unique positive root $u_\lambda$, which may be explicitly computed. The function $u:\lambda\in\mathbb{R}^-\longmapsto u_\lambda\in\mathbb{R}^+$ is real-analytic and $\mathcal{P}_*$ is the graph of $u$ (see Figure \ref{regions}). For convenience, if $\lambda\geq 0$ we will define $u_\lambda=\infty$.

We next introduce some functions which will play an essential role on the parameterization of B-curves. Let $\sigma:(s,\lambda,e_1)\in\mathbb{R}\times\mathcal{P}\longmapsto \sigma_{\lambda,e_1}(s)\in\mathbb{Z}_2$ be defined by
\begin{equation*}
\sigma_{\lambda,e_1}(s):=\left\{\begin{split} &1\,, &&\text{if } s\in\left[2k\omega_{\lambda,e_1},(2k+1)\omega_{\lambda,e_1}\right),\, k\in\mathbb{Z} \\
&(-1)^{\chi(\lambda,e_1)}\,, &&\text{if } s\in\left[(2k+1)\omega_{\lambda,e_1},2(k+1)\omega_{\lambda,e_1}\right),\, k\in\mathbb{Z} \end{split}\right.\,,
\end{equation*}
where $\omega_{\lambda,e_1}$ is the least period of $\mu_{\lambda,e_1}$ and $\chi:\mathcal{P}\longrightarrow\mathbb{Z}_2$ is the indicator function of $\mathcal{P}_*$, i.e., $\chi$ is zero everywhere but at the points $\left(\lambda,e_1\right)\in\mathcal{P}_*$ in which case $\chi(\lambda,e_1)=1$. We then define the \emph{angular function}
\begin{equation}\label{angular}
\theta_{\lambda,e_1}(s):=2\xi_{\lambda,e_1}\int_{\frac{\omega_{\lambda,e_1}}{2}}^s\frac{\mu_{\lambda,e_1}^2(t)\left(\mu_{\lambda,e_1}(t)+2\lambda\right)}{1-4\xi_{\lambda,e_1}^2\mu_{\lambda,e_1}^2(t)}\,dt
\end{equation}
and the \emph{radial} and \emph{height functions}, respectively,
\begin{eqnarray}
\rho_{\lambda,e_1}(s)&:=&\frac{\sigma_{\lambda,e_1}(s)\sqrt{4\xi_{\lambda,e_1}^2\mu_{\lambda,e_1}^2(s)-1}}{2\xi_{\lambda,e_1}\mu_{\lambda,e_1}(s)}\,,\label{radial}\\
h_{\lambda,e_1}(s)&:=&\frac{1}{2\xi_{\lambda,e_1}\mu_{\lambda,e_1}(s)}\,.\label{height}
\end{eqnarray}
From these definitions, some basic features of these functions can be deduced:
\begin{enumerate}
\item The height functions are periodic (with least period $\omega_{\lambda,e_1}$) and even. They have a minimum at $s=\omega_{\lambda,e_1}/2$ and a maximum at $s=\omega_{\lambda,e_1}$. (See Figure \ref{Fheight}.)
\begin{figure}[h!]
\centering
\includegraphics[width=0.35\linewidth]{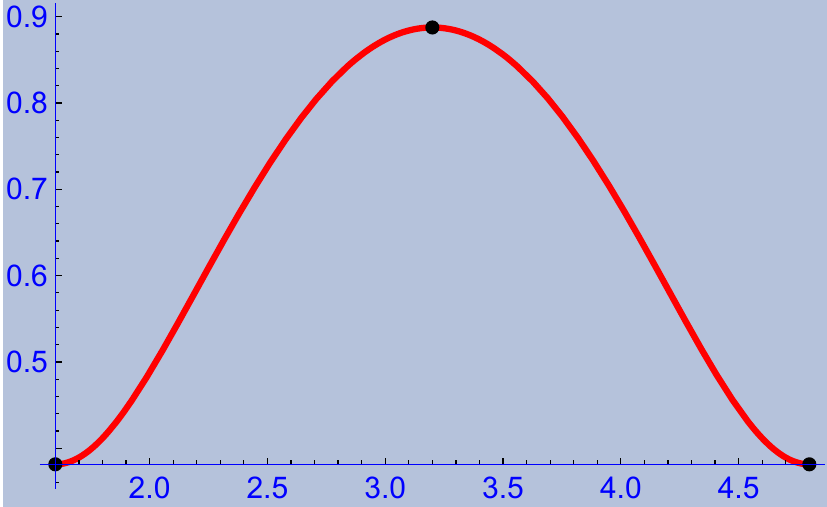}
\caption{The graph of the height function $h_{\lambda,e_1}$ for $\lambda=1.1$ and $e_1\simeq 1.46$ $(e_1<u_\lambda)$. In this picture we show the portion of the graph on the interval $[\omega_{\lambda,e_1}/2,\omega_{\lambda,e_1}/2+\omega_{\lambda,e_1}]$.}
\label{Fheight}
\end{figure}
\item The radial functions have two different qualitative behaviors depending on whether $e_1=u_\lambda$ or not:
\begin{enumerate}
\item If $e_1\neq u_\lambda$, then $\rho_{\lambda,e_1}$ is periodic (with least period $\omega_{\lambda,e_1}$) and positive. (See Figure \ref{Fradial}, Left.)
\item If $e_1=u_\lambda$, then $\rho_{\lambda,e_1}$ is periodic (with least period $2\omega_{\lambda,e_1}$) and $\rho_{\lambda,e_1}(s+\omega_{\lambda,e_1})=-\rho_{\lambda,e_1}(s)$. It has two zeros in the interval $[\omega_{\lambda,e_1}/2,\omega_{\lambda,e_1}/2+2\omega_{\lambda,e_1})$, precisely, at $s=\omega_{\lambda,e_1}$ and $s=2\omega_{\lambda,e_1}$. (See Figure \ref{Fradial}, Right.)
\end{enumerate}
\begin{figure}[h!]
\centering
\includegraphics[width=0.35\linewidth]{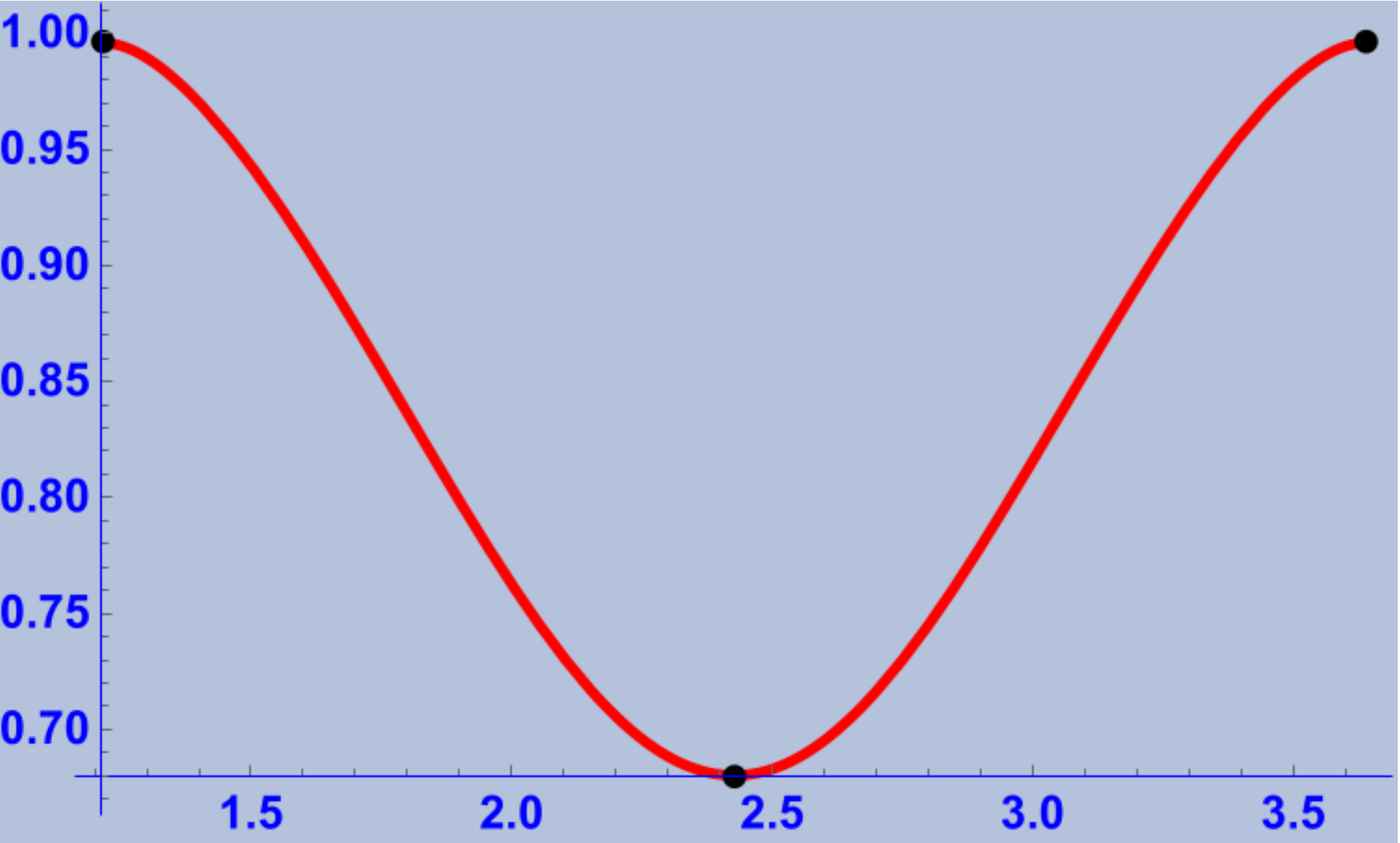}\quad\quad\quad\includegraphics[width=0.363\linewidth]{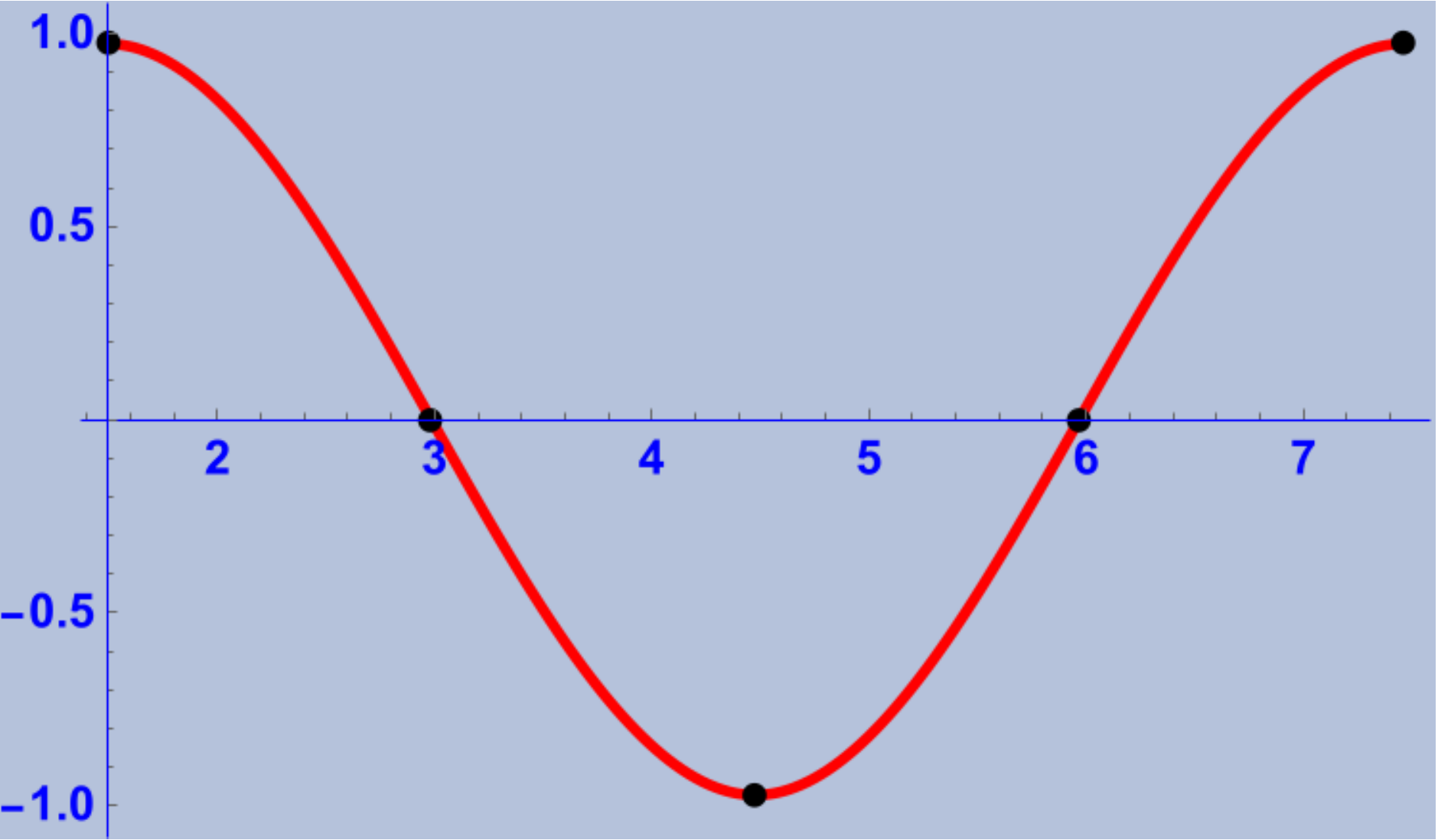}
\caption{The graphs of the radial functions $\rho_{\lambda,e_1}$ for: $\lambda=-1.1$ and $u_\lambda<e_1\simeq 4.59$ (left); and, $\lambda=-0.27$ and $e_1=u_\lambda\simeq 2.34$ (right). In the left figure we take the portion of the graph on the interval $[\omega_{\lambda,e_1}/2,\omega_{\lambda,e_1}/2+\omega_{\lambda,e_1}]$, while in the right one we are showing the graph on the interval $[\omega_{\lambda,e_1}/2,\omega_{\lambda,e_1}/2+2\omega_{\lambda,e_1}]$.}
\label{Fradial}
\end{figure}
\item The angular functions are arithmetic quasi-periodic (with quasi-period $\omega_{\lambda,e_1}$) and odd. In the interval $[\omega_{\lambda,e_1}/2,\omega_{\lambda,e_1}/2+\omega_{\lambda,e_1})$, $\theta_{\lambda,e_1}$ possesses an inflection point at $s=\omega_{\lambda,e_1}$ with $\theta_{\lambda,e_1}(\omega_{\lambda,e_1}/2)=\theta(\omega_{\lambda,e_1})/2$. However, their qualitative behavior depends on whether $e_1>u_\lambda$ or $e_1\leq u_\lambda$:
\begin{enumerate}
\item If $e_1>u_\lambda$, in the same interval, $\theta_{\lambda,e_1}$ has exactly two critical points, an absolute minimum somewhere between $s\in\left(\omega_{\lambda,e_1}/2,\omega_{\lambda,e_1}\right)$ and an absolute maximum in $(\omega_{\lambda,e_1},\omega_{\lambda,e_1}/2+\omega_{\lambda,e_1})$. These functions are increasing from the minimum to the maximum and they tend to $\infty$ as $s\to\infty$. (See Figure \ref{Fangular}, Left.)
\item If $e_1\leq u_\lambda$, then $\theta_{\lambda,e_1}$ is strictly decreasing and it tends to $-\infty$ as $s\to\infty$. (See Figure \ref{Fangular}, Right.)
\end{enumerate}
\begin{figure}[h!]
\centering
\includegraphics[width=0.35\linewidth]{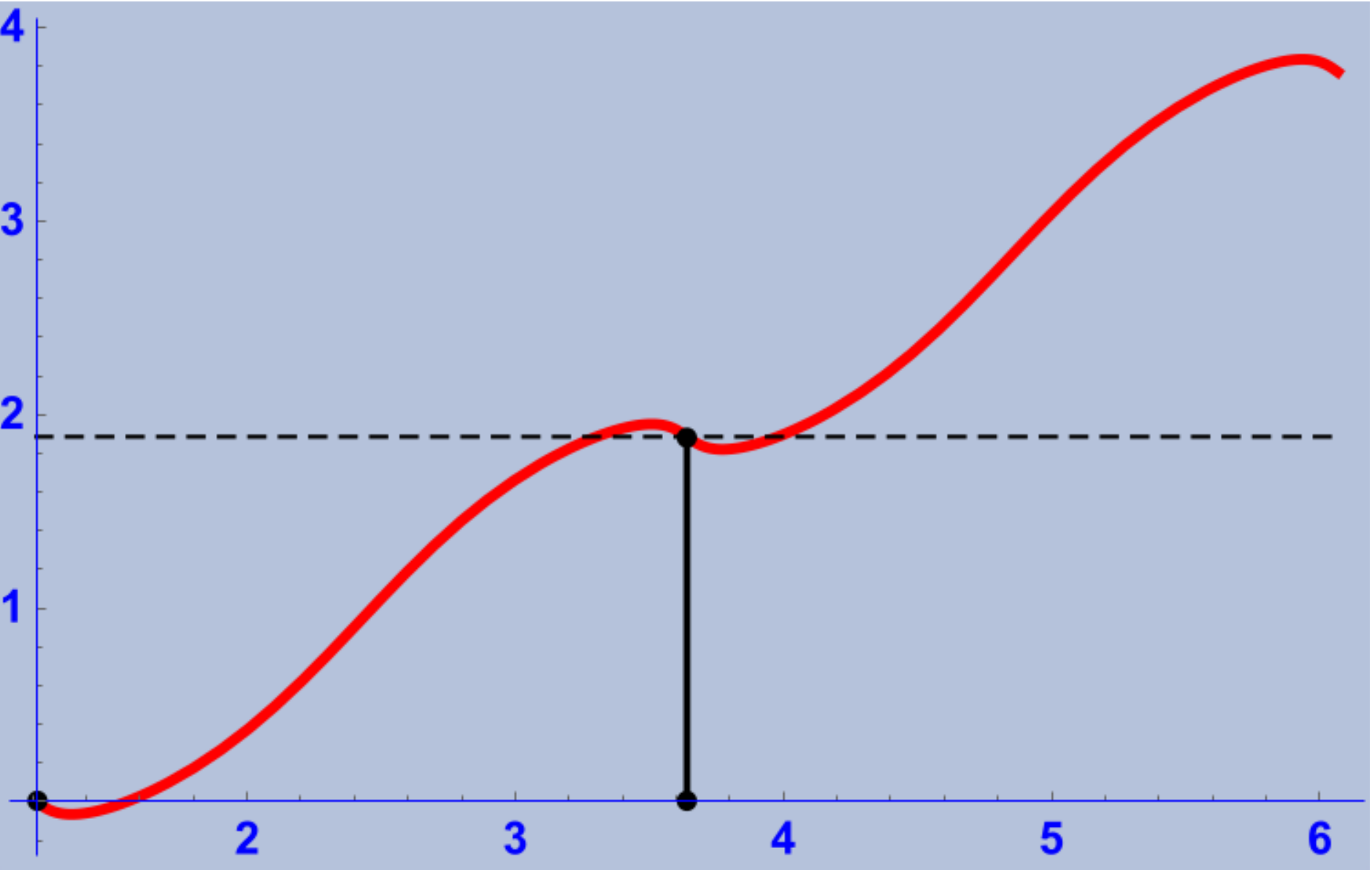}\quad\quad\quad\includegraphics[width=0.38\linewidth]{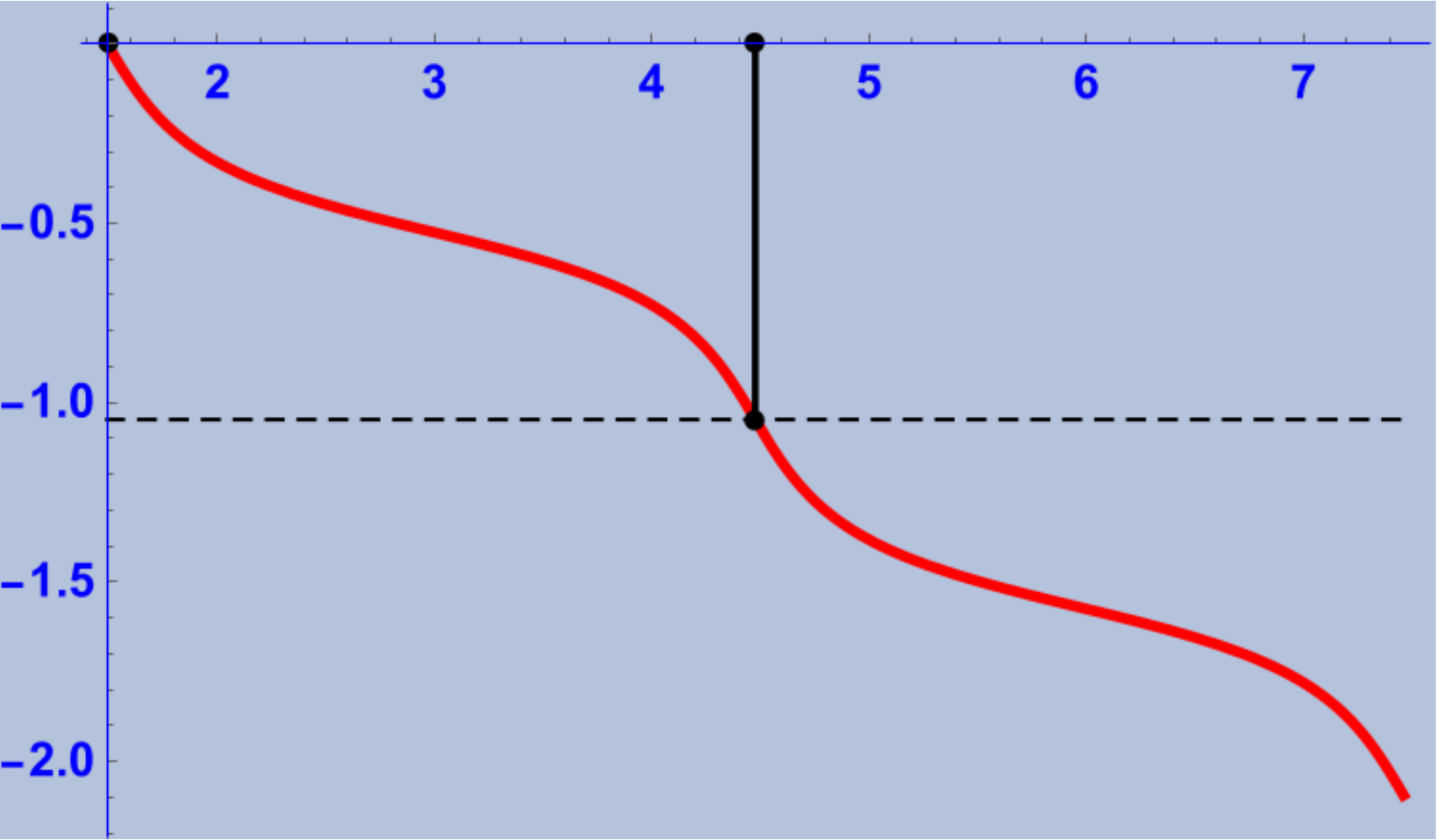}
\caption{The graphs of the angular functions $\theta_{\lambda,e_1}$ for: $\lambda=-1.1$ and $u_\lambda<e_1\simeq 4.59$ (left); and, $\lambda=-0.27$ and $e_1=u_\lambda\simeq 2.34$ (right). The black segment represents the jump of the quasi-periodic function. If the length of the segment is a rational multiple of $2\pi$, the critical curve is periodic (see Theorem \ref{Ccomplete}). In this pictures we show the portion of the graph on the interval $[\omega_{\lambda,e_1}/2,\omega_{\lambda,e_1}/2+2\omega_{\lambda,e_1}]$.}
\label{Fangular}
\end{figure}
\end{enumerate}

We have now all the necessary information to prove the parameterization of B-curves (the following result corresponds to Theorem \ref{A} in the Introduction).

\begin{thm}\label{Acomplete} Let $(\lambda,e_1)\in\mathcal{P}=\{(\lambda,e_1)\in\mathbb{R}^2\,\lvert\,e_1>\eta_\lambda\}$ where $\eta_\lambda$ is defined in Remark \ref{circles}, and $\mu$ be a solution of the Cauchy problem
\begin{equation*}
\left\{\begin{split} &\dot{\mu}^2=-\mu^2Q(\mu) \\
& \mu(\omega/2)=e_1 \end{split}\right.\,,
\end{equation*}
where $Q$ is the quartic polynomial defined in \eqref{polynomial} and $\omega$ is the least period of $\mu$. Then, $\mu\equiv\mu_{\lambda,e_1}$ is a positive periodic function. Moreover, define the curve
$$\gamma_{\lambda,e_1}(s)=\left(h_{\lambda,e_1}(s),-\rho_{\lambda,e_1}(s)\cos\theta_{\lambda,e_1}(s),\rho_{\lambda,e_1}(s)\sin\theta_{\lambda,e_1}(s)\right)$$
where $h_{\lambda,e_1}$, $\rho_{\lambda,e_1}$ and $\theta_{\lambda,e_1}$ are the height, radial and angular functions introduced above, \eqref{angular}-\eqref{height}. Then, $\gamma_{\lambda,e_1}$ is the B-curve with parameters $(\lambda,e_1)$, in its standard form. 
\end{thm}
\textit{Proof.} The first part of the statement was already proven in previous section. It remains to prove that $\gamma_{\lambda,e_1}$ is the B-curve with parameters $(\lambda,e_1)$. For simplicity, we will omit the subscripts $\lambda$ and $e_1$ throughout this proof.

We first note that for every $(\lambda,e_1)\in\mathcal{P}$ the function
$$s\longmapsto \frac{\mu(s)+2\lambda}{1-4\xi^2\mu^2(s)}$$
is real-analytic. In fact, if $(\lambda,e_1)$ is not exceptional, then $1-4\xi^2\mu^2(s)\leq1-4\xi^2e_2^2<0$. Note that $(\lambda,e_1)$ is exceptional if and only if $\xi=-1/(4\lambda)$ and $\lambda=-e_2/2<0$. Thus, in this case
$$\frac{\mu(s)+2\lambda}{1-4\xi^2\mu^2(s)}=4\lambda^2\frac{\mu(s)+2\lambda}{4\lambda^2-\mu^2(s)}=-4\lambda^2\frac{1}{\mu(s)-2\lambda}=\frac{-4\lambda^2}{\mu(s)+e_2}<0\,.$$
This implies that the angular function $\theta(s)$ is real-analytic. By construction, so is the height function $h(s)$.

We next distinguish between the cases where $(\lambda,e_1)\in\mathcal{P}$ is exceptional or not.

Suppose first that $(\lambda,e_1)\in\mathcal{P}$ is not exceptional. Then, as shown above $4\xi^2\mu^2-1>0$ holds and, hence, $\gamma(s)$ is a real-analytic spherical curve. Moreover, since $\mu(\omega/2)=e_1$, we have
$$\gamma\left(\frac{\omega}{2}\right)=\left(h\left(\frac{\omega}{2}\right),-\rho\left(\frac{\omega}{2}\right),0\right)=\frac{1}{2\xi e_1}\left(1,-\sqrt{4\xi^2 e_1^2-1},0\right).$$
Computing the derivative we get
$$\dot{\gamma}(s)\cdot\dot{\gamma}(s)=-\frac{4\lambda^2\mu^2+4\lambda\mu^5+\mu^6+\dot{\mu}^2}{\mu^2-4\xi^2\mu^4}=1\,,$$
and $\dot{\gamma}(\omega/2)=(0,0,-1)$. Here we have used the conservation law \eqref{ODE} to simplify above expression. Hence, $\gamma$ is an arc-length parameterized spherical curve. Then, computing the second order derivative of $\gamma$, and using the Euler-Lagrange equation \eqref{EL} and the conservation law \eqref{ODE} to simplify the expression, we obtain that
$$\ddot{\gamma}(s)\cdot\ddot{\gamma}(s)=\mu^4(s)+1\,,$$
from which we conclude that the geodesic curvature of $\gamma$ is either $\mu^2$ or $-\mu^2$. In what follows, we will discard the second case. From above computations we also deduce that
$$\gamma\left(\frac{\omega}{2}\right)\times\dot{\gamma}\left(\frac{\omega}{2}\right)=\left(\frac{e_1(e_1+2\lambda)}{\sqrt{1+(e_1[e_1+2\lambda])^2}},\frac{1}{\sqrt{1+(e_1[e_1+2\lambda])^2}},0\right).$$
This implies
$$\kappa\left(\frac{\omega}{2}\right)=\ddot{\gamma}\left(\frac{\omega}{2}\right)\cdot\left(\frac{e_1(e_1+2\lambda)}{\sqrt{1+(e_1[e_1+2\lambda])^2}},\frac{1}{\sqrt{1+(e_1[e_1+2\lambda])^2}},0\right).$$
A straightforward computation yields
$$\ddot{\gamma}\left(\frac{\omega}{2}\right)=\left(\frac{e_1^4+2\lambda e_1^3-1}{\sqrt{1+(e_1[e_1+2\lambda])^2}},\frac{2e_1(e_1+\lambda)}{\sqrt{1+(e_1[e_1+2\lambda])^2}},0\right),$$
and, hence, $\kappa(\omega/2)=e_1^2=\mu^2(\omega/2)$. That is $\kappa=\mu^2$ and $\gamma$ is an arc-length parameterized curve with curvature $\mu^2$ and the same initial conditions of the (standard) B-curve. Consequently, both of them coincide.

Finally, we treat the case where $(\lambda,e_1)\in\mathcal{P}$ is exceptional. In this case $4\xi^2\mu^2-1$ vanishes on $\omega\mathbb{Z}$ and is positive on $\mathbb{R}\setminus \omega\mathbb{Z}$. By definition,
\begin{equation*}
\gamma(s)=\left\{\begin{split} &\left(h(s),-\widehat{\rho}(s)\cos\theta(s),\widehat{\rho}(s)\sin\theta(s)\right), &&\text{if } s\in\left[2k\omega,(2k+1)\omega\right),\, k\in\mathbb{Z} \\
&\left(h(s),\widehat{\rho}(s)\cos\theta(s),-\widehat{\rho}(s)\sin\theta(s)\right), &&\text{if } s\in\left[(2k+1)\omega,2(k+1)\omega\right),\, k\in\mathbb{Z} \end{split}\right.\,,
\end{equation*}
where $\widehat{\rho}=\sqrt{4\xi^2\mu^2-1}/(2\xi\mu)$. 

It is then clear that $\gamma$ is continuous and real-analytic away from the discrete set $\omega\mathbb{Z}$. On the set $\mathbb{R}\setminus\omega\mathbb{Z}$ we can argue as in the non-exceptional case to prove that $\gamma$ is an arc-length parameterized spherical curve with $\kappa=\mu^2$ and such that
\begin{equation*}
\left\{\begin{split} &\gamma\left(\frac{\omega}{2}\right)=\frac{1}{2\xi e_1}\left(1,-\sqrt{4\xi^2e_1^2-1},0\right)\\
&\dot{\gamma}\left(\frac{\omega}{2}\right)=\left(0,0,-1\right)
\end{split}\right.\,.
\end{equation*}
We next prove that $\gamma$ is of class $\mathcal{C}^1$ on $\omega\mathbb{Z}$. By construction, $\gamma(s+k\omega)=R^k\cdot\gamma(s)$, for any $k\in\mathbb{Z}$, where $R$ is a rotation around the $Ox$-axis. Hence, it suffices to show that $\gamma$ is of class $\mathcal{C}^1$ at $s=\omega$. Note that the first component of $\gamma$ is everywhere real-analytic. Since $\theta$ is a real-analytic function, the second and third components of $\gamma$ are of class $\mathcal{C}^1$ at $s=\omega$ if and only if the function
\begin{equation*}
f(s)=\left\{\begin{split} &\sqrt{4\xi^2\mu^2(s)-1}\,, &&s\in[0,\omega)\\
&-\sqrt{4\xi^2\mu^2(s)-1}\,, &&s\in[\omega,2\omega)
\end{split}\right.\,,
\end{equation*}
is of class $\mathcal{C}^1$ at $s=\omega$. Recall that $\mu(\omega)=e_2$ and that, in the exceptional case $\xi=1/(2e_2)$. Thus, $f$ is continuous at $s=\omega$. Moreover, on the left of $\omega$ we have
$$\dot{\mu}=\mu\sqrt{-(\mu-e_1)(\mu-e_2)(\mu-e_3)(\mu-e_4)}\,,$$
where $e_1>e_2>0$, $e_3$ and $e_4$ are the roots of the polynomial $Q$, \eqref{polynomial}. Therefore,
$$\dot{f}=\frac{\mu\dot{\mu}}{e_2\sqrt{\mu^2-e_2^2}}=\frac{\mu^2\sqrt{-(\mu-e_1)(\mu-e_3)(\mu-e_4)}}{e_2\sqrt{\mu+e_2}}\,.$$
While, on the right of $\omega$,
$$\dot{\mu}=-\mu\sqrt{-(\mu-e_1)(\mu-e_2)(\mu-e_3)(\mu-e_4)}\,,$$
and so
$$\dot{f}=-\frac{\mu\dot{\mu}}{e_2\sqrt{\mu^2-e_2^2}}=\frac{\mu^2\sqrt{-(\mu-e_1)(\mu-e_3)(\mu-e_4)}}{e_2\sqrt{\mu+e_2}}\,.$$
This implies that
$$\lim_{s\to\omega^-}\dot{f}(s)=\lim_{s\to\omega^+}\dot{f}(s)=\frac{\sqrt{(e_1-e_2)(e_2-e_3)(e_2-e_4)}}{\sqrt{2}}\,.$$
Consequently, $\gamma$ is an arc-length parameterized spherical curve of class $\mathcal{C}^1$. Let $\widetilde{\gamma}$ be the standard B-curve with parameters $\lambda$ and $e_1$. Then, both $\gamma$ and $\widetilde{\gamma}$ are arc-length parameterized spherical curves with the same curvature and the same initial conditions, $\gamma(\omega/2)=\widetilde{\gamma}(\omega/2)$ and $\dot{\gamma}(\omega/2)=\dot{\widetilde{\gamma}}(\omega/2)$. Since $\widetilde{\gamma}$ is real-analytic and $\gamma$ is of class $\mathcal{C}^1$ and real-analytic on $\mathbb{R}\setminus\omega\mathbb{Z}$, we may conclude that $\gamma=\widetilde{\gamma}$.

More precisely, let $\mathcal{F}$ and $\widetilde{\mathcal{F}}$ be the (spherical) Frenet frame fields along $\gamma$ and $\widetilde{\gamma}$, respectively. By construction, $\mathcal{F}$ is continuous and real-analytic on $\mathbb{R}\setminus\omega\mathbb{Z}$, while $\widetilde{\mathcal{F}}$ is real-analytic everywhere. Then there exists a continuous map $A:\mathbb{R}\longrightarrow SO(3)$ such that $\widetilde{\mathcal{F}}=A\mathcal{F}$. On $\mathbb{R}\setminus\omega\mathbb{Z}$, $\mathcal{F}$ and $\widetilde{\mathcal{F}}$ are both solutions of the linear system
$$\dot{X}=X\cdot\begin{pmatrix} 0 & -1 & 0 \\ 1 & 0 & -\mu^2 \\ 0 & \mu^2 & 0 \end{pmatrix}.$$
Thus, $A$ is constant on $\mathbb{R}\setminus\omega\mathbb{Z}$. By continuity, one concludes that $A$ is constant everywhere. Moreover, since $\mathcal{F}(\omega/2)=\widetilde{\mathcal{F}}(\omega/2)$, we obtain that $A=Id_{3x3}$ is the identity. This proves that $\gamma=\widetilde{\gamma}$. \hfill$\square$

\begin{rem} A B-curve $\gamma_{\lambda,e_1}$ passes through the pole $(1,0,0)$ of the unit-sphere $\mathbb{S}^2$ if and only if the parameters $(\lambda,e_1)\in\mathcal{P}$ are exceptional, i.e., if $(\lambda,e_1)\in\mathcal{P}^*$ or, equivalently, if $e_1=u_\lambda$. Moreover, the B-curve passes through the pole, precisely, whenever $s\in\omega\mathbb{Z}$.
\end{rem}

The \emph{jump} of the B-curve $\gamma$ with parameters $(\lambda,e_1)\in\mathcal{P}$ is defined by
 \begin{equation}\label{jump}
\Psi(\lambda,e_1)\equiv \Psi_{\lambda}(e_1):=\theta_{\lambda,e_1}\left(\omega_{\lambda,e_1}\right)-\theta_{\lambda,e_1}(0)\,.
\end{equation}
The curve is periodic if and only if $\Psi_{\lambda}(e_1)$ is a rational multiple of $2\pi$. (See Figure \ref{Fangular}.)

The jump function can be seen as a period map in the following way. Let $SO(2)\subset SO(3)$ be the stabilizer of the vector ${\vec i}=(1,0,0)$. From Theorem \ref{Acomplete} it follows that the monodromy map $\mathfrak{m}$ is $SO(2)$-valued. Then, since $\mathcal{P}$ is contractible and $\mathfrak{m}$ is continuous, there exist continuous functions $\Phi:\mathcal{P}\longrightarrow \mathbb{R}$ such that
\begin{equation*}
\mathfrak{m}_{\lambda,e_1}=\begin{pmatrix} 1 & 0 & 0 \\ 0 & \cos\Phi_{\lambda,e_1} & -\sin\Phi_{\lambda,e_1} \\ 0 & \sin\Phi_{\lambda,e_1} & \cos\Phi_{\lambda,e_1} \end{pmatrix}.
\end{equation*}
We say that $\Phi$ is a \emph{period map}. The period map is unique ${\rm mod}\,\mathbb{Z}$ and $\Phi_{\lambda,e_1}\equiv \Psi_{\lambda}(e_1) \left({\rm mod}\,1\right)$. Keeping in mind that $\omega$ is real-analytic on the open set $\widehat{\mathcal{P}}$, that $\xi$ and the integrand are real-analytic on $\mathcal{P}$, we deduce that the period maps are real-analytic on $\widehat{\mathcal{P}}$.

We now introduce the strictly increasing real-analytic function $p:\mathbb{R}\longrightarrow \left(-1,-1/\sqrt{3}\,\right)$ given by
\begin{equation}\label{p}
p(\lambda):=-\sqrt{\frac{1+\eta_{\lambda}^4}{3+\eta_\lambda^4}}=-\sqrt{\frac{1-\lambda\eta_\lambda^3}{2-\lambda\eta_\lambda^3}}\,,
\end{equation}
where $\eta_\lambda$ is the (unique) positive solution of \eqref{cst}. This function will arise when analyzing the asymptotic behavior of $\Psi_{\lambda}(e_1)$.

We next prove the second main result of the paper, which shows the existence of periodic B-curves, i.e., B-strings (this result corresponds to Theorem \ref{C} in the Introduction).

\begin{thm}\label{Ccomplete} Let $\lambda\in\mathbb{R}$. Then, there exists an unbounded countable set $\Delta_\lambda\subset(\eta_\lambda,\infty)$ such that a B-curve $\gamma$ with parameters $\lambda$ and $e_1$ is periodic if and only if $e_1\in\Delta_\lambda$. In other words, for every $q\in\left(1+p(\lambda),1/2\right)\cap\mathbb{Q}$ there exists a B-string $\gamma$ with Lagrange multiplier $\lambda$ and monodromy
\begin{equation*}
\mathfrak{m}_q\equiv\mathfrak{m}_{\lambda,e_1}=\begin{pmatrix} 1 & 0 & 0 \\ 0 & \cos(2\pi q) & -\sin(2\pi q) \\ 0 & \sin(2\pi q) & \cos(2\pi q) \end{pmatrix}.
\end{equation*}
The rational number $q$ is the characteristic number of the B-string.
\end{thm}
\textit{Proof.} In order to show the existence of B-strings, we need to analyze the asymptotic behavior of the jump function (for simplicity, we will omit the subscript $\lambda$)
$$\Psi(e_1)=2\xi\int_0^\omega \frac{\mu^2\left(\mu+2\lambda\right)}{1-4\xi^2\mu^2}\,ds=4\xi\int_{e_2}^{e_1}\frac{\mu\left(\mu+2\lambda\right)}{(1-4\xi^2\mu^2)\sqrt{-(\mu-e_1)(\mu-e_2)(\mu-e_3)(\mu-e_4)}}\,d\mu\,.$$
This is a rather standard integral that can be solved in terms of complete elliptic integrals of the first and third kind (for details, see Appendix A). Moreover, by the properties of these elliptic integrals we obtain that $\Psi$ is real-analytic on the following sets:
$$\mathcal{P}_-=\{(\lambda,e_1)\in\mathcal{P}\,\lvert\,e_1<u_\lambda\}\,,\quad\quad\quad\mathcal{P}_+=\{(\lambda,e_1)\in\mathcal{P}\,\lvert\,e_1>u_\lambda\}\,,$$
and $\mathcal{P}_*$ (see Figure \ref{regions} for the plot of these regions). However, the jump function $\Psi$ is not real-analytic on $\mathcal{P}$ since it has a jump discontinuity on the exceptional locus $\mathcal{P}_*$ (see Remark \ref{discontinuity}).

\begin{figure}[h!]
\centering
\includegraphics[width=0.3\linewidth]{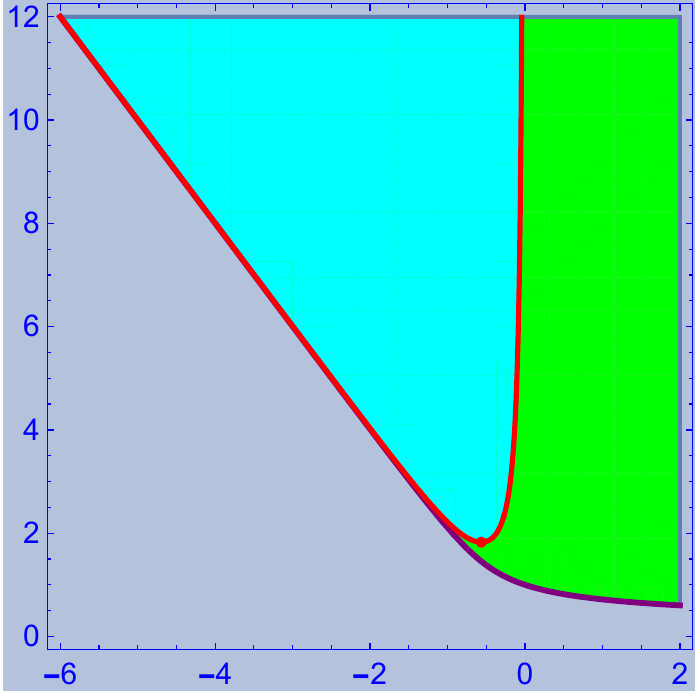}
\caption{The exceptional locus $\mathcal{P}_*$ (in red) and the regions $\mathcal{P}_-$ (in green) and $\mathcal{P}_+$ (in blue). The colored domain represents the region $\mathcal{P}=\mathcal{P}_-\cup\mathcal{P}_*\cup\mathcal{P}_+$ and its boundary is $\eta_\lambda$ (in purple).}
\label{regions}
\end{figure}

We next obtain the limits of $\Psi$ as the parameter $e_1$ approaches the boundaries of its domain, i.e., $e_1\in(\eta_\lambda,\infty)$. The proof of the following limits is just a technical computation involving the decomposition of $\Psi$ as the sum of three complete elliptic integrals. For the sake of clarity we simply state the limits here, while we postpone the proof to Appendix A.

We begin with the limit when $e_1\to\eta_\lambda^+$. In this case, for every $\lambda\in\mathbb{R}$ we have
$$\lim_{e_1\to\eta_\lambda^+} \Psi(e_1)=2\pi p(\lambda)\,,$$
where $p(\lambda)$ is, precisely, the function introduced in \eqref{p}. On the other hand, we have the following limits when $e_1\to\infty$, depending on the sign of $\lambda$,
$$\lim_{e_1\to\infty} \Psi(e_1)=-\pi\,,$$
if $\lambda\geq 0$, and
$$\lim_{e_1\to\infty} \Psi(e_1)=\pi\,,$$
if $\lambda<0$.

Finally, let $\Phi_{\pm}:\mathcal{P}\longrightarrow\mathbb{R}$ be the period maps such that $\Phi_-\lvert_{\mathcal{P}_-}=\Psi\lvert_{\mathcal{P}_-}$ and $\Phi_{+}\lvert_{\mathcal{P}_+}=\Psi\lvert_{\mathcal{P}_+}$. Then, there exists a $k\in\mathbb{Z}$ such that $\Phi_-=\Phi_++2\pi k$. From above limits, $k=-1$ and $\{ \Phi_-(\lambda,e_1)/2\pi\,\lvert\,e_1>\eta_\lambda\}$ contains the open interval $(-1/2,p(\lambda))\equiv (1+p(\lambda),1/2)\,({\rm mod}\,1)$. This concludes the proof. \hfill$\square$
\\

As mentioned in the previous proof, the function $\Psi_\lambda$ has a jump discontinuity on $\mathcal{P}_*$, i.e., when $e_1=u_\lambda$. Consequently, it is convenient to regularize this function in order to work with a continuous function. Let us define the \emph{regularized jump function} $\widehat{\Psi}:\mathcal{P}\longrightarrow\mathbb{R}$ by 
\begin{equation}\label{mjump}
\widehat{\Psi}(\lambda,e_1)\equiv\widehat{\Psi}_\lambda (e_1)=\left\{\begin{split} &\Psi_\lambda(e_1)\left({\rm mod}\,2\pi \right), &&e_1\neq u_\lambda\\
&\Psi_\lambda(e_1)+\pi\,, &&e_1=u_\lambda
\end{split}\right.\,.
\end{equation}
In terms of $\widehat{\Psi}_\lambda$, a B-curve $\gamma$ with parameters $(\lambda,e_1)\in\mathcal{P}$ is closed, i.e., a B-string, if and only if $\widehat{\Psi}_\lambda(e_1)=2\pi q$, where $q\in\mathbb{Q}$ is a rational number of the type $q=m/n$ for relatively prime natural numbers $m<n$. The rational number $q$ is the \emph{characteristic number} of the B-string, while the natural number $n$ is its \emph{wave number}.

From Theorem \ref{Ccomplete} we conclude that for any pair $(m,n)$ of relatively prime natural numbers and $\lambda\in\mathbb{R}$ satisfying
$$2\left(1+p(\lambda)\right)n<2m<n$$
there exists a B-string with multiplier $\lambda\in\mathbb{R}$ and characteristic number $q=m/n$.

\begin{rem}\label{discontinuity} Numerical experiments strongly support the ansatz that the regularized jump function $\widehat{\Psi}_\lambda$ is a real-analytic period map for the $1/2$-Bernoulli's bending variational problem. In addition, for every $\lambda\in\mathbb{R}$ the function $\widehat{\Psi}_\lambda:e_1\in(\eta_\lambda,\infty)\longmapsto(2\pi[1+p(\lambda)],\pi)$ is a strictly increasing diffeomorphism. This experimental fact would lead to the following conclusions (which are stronger than those of Theorem \ref{Ccomplete}):
\begin{enumerate}
\item For every $\lambda\in\mathbb{R}$ and $q\in (1+p(\lambda),1/2)\cap\mathbb{Q}$ there exists a unique B-string $\gamma$ with multiplier $\lambda$ such that $e_1=\widehat{\Psi}^{-1}_\lambda(2\pi q)$. In particular, B-strings $\gamma$ with multiplier $\lambda$ are in one-to-one correspondence with the countable set $(1+p(\lambda),1/2)\cap\mathbb{Q}$.
\item The wave number of $\gamma\equiv\gamma_q$ is $n\geq 3$, and hence B-strings are not embedded (see Lemma \ref{lem} for the number of self-intersection points).
\end{enumerate}
\end{rem}

In the last part of this section we will prove Theorem \ref{B}. We first show a technical lemma regarding the number of points of self-intersection.

\begin{lem}\label{lem} Let $\gamma$ be a B-string with multiplier $\lambda$ and characteristic number $q=m/n$. Assume also that $e_1\neq u_\lambda$. Then, 
\begin{enumerate}
\item The points $\gamma(\omega/2+k\omega)$ and $\gamma(\omega+k\omega)$, $k=0,...,n-1$ are simple.
\item If $e_1<u_\lambda$, then $\gamma$ possesses, exactly, $n\left(n-m-1\right)$ ordinary double points.
\item If $e_1>u_\lambda$, then $\gamma$ possesses, at least, $n m$ points of self-intersection. The angular function has a unique absolute maximum at $s_*\in[\omega/2,\omega]$ such that $\theta(s_*)>\theta(\omega)=\pi m /n$, and we have three different scenarios:
\begin{enumerate}
\item Case $\theta(s_*)<\pi(m+1)/n$. In this case the B-string has exactly $n m$ ordinary double points.
\item Case $\pi(m+k)/n<\theta(s_*)<\pi(m+k+1)/n$ for some $k\in\mathbb{N}$. In this case the B-string has exactly $n(m+2k)$ points of self-intersection.
\item Case $\pi(m+k)/n<\theta(s_*)=\pi(m+k+1)/n$ for some $k\in\mathbb{N}\cup\{0\}$. In this case the B-string has exactly $n(m+2k+1)$ points of self-intersection.
\end{enumerate}
\end{enumerate}
\end{lem}
\textit{Proof.} We begin by proving that the points $\gamma(\omega/2+k\omega)$ and $\gamma(\omega+k\omega)$, $k=0,...,n-1$ are simple. Since
$\gamma(s+\omega)=R_{2\pi m/n}\cdot \gamma(s)$, where $R$ is a rotation around the $Ox$-axis (see previous Theorems \ref{Acomplete} and \ref{Ccomplete}) it suffices to show that $\gamma(\omega/2)$ and $\gamma(\omega)$ are simple. Consider $\gamma(\omega/2)$. We exhibit that, if $s_*\in [\omega/2,\omega/2+n\omega)$ satisfies $\gamma(s_*)=\gamma(\omega/2)$, then $s_*=\omega/2$. The radial function reaches its maximum at $\omega/2+h\omega$, $h\in {\mathbb Z}$.  Hence $s_*=\omega/2+p\omega$, for some integer $p=0,\dots  n-1$. In addition, $\theta(s_*)\equiv \theta(\omega/2)=0 ({\rm mod}\, 2\pi)$. Then, $\theta(s_*)=2k\pi$, $k\in {\mathbb Z}$. On the other hand, $\theta(s_*)=\theta(\omega/2+p\omega)=2\pi mp/n$. This implies $mp=kn$. Since $m$ and $n$ are relatively prime, we have $p=hn$, $h\in {\mathbb Z}$. Therefore, $s_*=\omega/2+hn\omega$.  But $\omega/2\le s_*<\omega/2+n\omega$. Thus, $h=0$ and $s_*=\omega/2$. Next, consider $\gamma(\omega)$. The radial function reaches its minimum at $h\omega$, $h\in {\mathbb Z}$. Thus, if $s_*\in [\omega/2,\omega/2+n\omega)$ satisfies $\gamma(s_*)=\gamma(\omega)$, then $s_*=p\omega$,  for some integer $p=0,\dots n$. In addition, $\theta(s_*)\equiv \theta(\omega) = 2\pi m/n$, mod $2\pi{\mathbb Z}$, ie $\theta(s_*)=2\pi(m/n+k)$, $k\in {\mathbb Z}$. On the other hand, $\theta(s_*)=\theta(p\omega) =2\pi pm/n$. Then, $m(p-1)=nk$. Thus, $p-1=hn$, $h\in {\mathbb Z}$. Since $1\le p \le n$, the only option is $p=1$. This proves that $s_*=\omega$.

To prove the other assertions, let $k_1$ and $k_2$ be the integers such that
\begin{equation}\label{help}
k_1n+k_2m=1\,,
\end{equation}
whose existence is guaranteed because $m$ and $n$ are relatively prime.

Consider first that $e_1<u_\lambda$. The angular function is strictly decreasing on $(\omega/2,\omega)$ and $\theta(\omega/2)=0$ while $\theta(\omega)=-2\pi (n-m)/n$. Thus, for every $j=1,...,n-m-1$ there exists a unique $s_j\in(\omega/2,\omega)$ such that $\theta(s_j)=-\pi j/n$. We are going to prove that, for every $k=0,...,n-1$, $ \gamma(s_j-k\omega)$ is a multiple point. For this purpose, consider
$$\widehat{s}_j=-s_j-\left(k+jk_2\right)\omega\,.$$
The height and radial functions are even and periodic, with period $\omega$. Then, $\rho(s_j)=\rho(\widehat{s}_j)$ and $h(s_j)=h(\widehat{s}_j)$. Taking into account that $\theta$ is an odd function and using \eqref{help} we have
\begin{eqnarray*}
\theta(\widehat{s}_j)&=&-\theta(s_j)-2\pi\left(k+jk_2\right))\frac{m}{n}=\pi\frac{j}{n}-2\pi j\frac{1-k_1n}{n}-2\pi k \frac{m}{n}\\
&=&-\pi\frac{j}{n}-2\pi k \frac{m}{n}+2\pi j k_1=\theta(s_j-k\omega)+2\pi j k_1\,.
\end{eqnarray*}
This implies that $\gamma(s_j)=\gamma(\widehat{s}_j)$. It remains to prove that $s_j\not\equiv \widehat{s}_j\left({\rm mod}\,n\omega\right)$. By contradiction, suppose $s_j\equiv \widehat{s}_j\left({\rm mod}\,n\omega\right)$, then there exists a $p\in\mathbb{Z}$ such that $2s_j=p\omega$, but this is impossible because $s_j\in (\omega/2,\omega)$. We next prove that the points $\gamma(s_j-k\omega)$, $1\leq j\leq n-m-1$ and $0\leq k\leq n-1$ are distinct. Let $0\leq \widehat{k}<k\leq n-1$,  then
$$\theta(s_j-k\omega)=-\pi\frac{j}{n}-2\pi k\frac{m}{n}=\theta(s_j-\widehat{k}\omega)-2\pi\frac{m(k-\widehat{k})}{n}\,.$$
Since $m$ and $n$ are relatively prime, $m(k-\widehat{k})/n$ is an integer number if and only if $k-\widehat{k}$ is an integer multiple of $n$. On the other hand, it is clear that $1\leq k-\widehat{k}<n-1$.  Hence, $m(k-\widehat{k})/n\not\in\mathbb{Z}$ and so $\theta(s_j-k\omega)\not\equiv\theta(s_j-\widehat{k}\omega)\left({\rm mod}\,2\pi\right)$. This implies that $\gamma(s_j-k\omega)\neq\gamma(s_j-\widehat{k}\omega)$.

Suppose now that $1\leq j<\widehat{j}\leq n-m-1$, then $\theta(s_j)=-\pi j/n>-\pi\widehat{j}/n=\theta(s_{\widehat{j}})$. Since $\theta$ is strictly decreasing, we have $\omega/2<s_j<s_{\widehat{j}}<\omega$. On the other hand, $\rho$ is also strictly decreasing on $(\omega/2,\omega)$, so $\rho(s_j-k\omega)=\rho(s_j)\neq\rho(s_{\widehat{j}})=\rho(s_{\widehat{j}}-\widehat{k}\omega)$, for every $k,\widehat{k}=0,...,n-1$.

In conclusion, $\gamma(s_j-k\omega)\neq\gamma(s_{\widehat{j}}-\widehat{k}\omega)$, for every $k,\widehat{k}=0,...,n-1$, and it then follows that the points $\gamma(s_j-k\omega)$, $1\leq j\leq n-m-1$, $0\leq k\leq n-1$ are distinct. This proves that $\gamma$ possesses, at least, $n(n-m-1)$ points of self-intersection.

We finally conclude from the first assertion and from the fact that the angular function is strictly decreasing that the points constructed above are the only multiple points of $\gamma$.

In what follows, we analyze the case $e_1>u_\lambda$. Denote by $\widehat{s},\widetilde{s}\in(\omega/2,\omega)$ the absolute minimum and absolute maximum, respectively, of $\theta$ on the closed interval $[\omega/2,\omega]$. Then, for every $j=1,...,m$ there exists a unique $s_j\in\left(\widehat{s},\widetilde{s}\right)$ such that $\theta(s_j)=\pi j/n$. Then, arguing as in the case $e_1<u_\lambda$, we prove that $\gamma(s_j-k\omega)$, $j=1,...,m$ and $k=0,...,n-1$ are distinct multiple points of $\gamma$. However, contrary to the case $e_1<u_\lambda$, these may not be the only multiple points of $\gamma$ since the angular function has a unique absolute maximum at $s_*\in[\omega/2,\omega]$ such that $\theta(s_*)>\theta(\omega)=\pi m /n$, and we have three different scenarios:
\begin{enumerate}
\item[(a)] Case $\theta(s_*)<\pi(m+1)/n$. In this case the B-string has exactly $n m$ points of self-intersection. These points are the ones constructed above. (See Figure \ref{Example5-1} of Subsection 5.4.)
\item[(b)] Case $\pi(m+k)/n<\theta(s_*)<\pi(m+k+1)/n$. In this case the B-string has exactly $n(m+2k)$ points of self-intersection. All of them are ordinary double points. The first $n m$ multiple points are the ones constructed in the proof, while the others arise as the two solutions of the equations $\theta(s)=\pi(m+p)/n$ with $p=1,...,k$. (See Figure \ref{Example7-1} of Subsection 5.6.)
\item[(c)] Case $\pi(m+k)/n<\theta(s_*)=\pi(m+k+1)/n$. In this case the B-string has exactly $n(m+2k+1)$ points of self-intersection. The first $n(m+2k)$ are ordinary double points (as in previous case), while the remaining $n$ points are points of tangential self-intersection. More precisely, these points are, exactly, $\gamma(s_*+p\omega)$, $p=0,...,n-1$. (See Figure \ref{Example6-1} of Subsection 5.5.)
\end{enumerate}
This finishes the proof. \hfill$\square$
\\

A similar argument may be used to obtain the analogue result for the case $e_1=u_\lambda$, which we state in the following remark.

\begin{rem}\label{remarknew} Let $\gamma$ be a B-string with multiplier $\lambda$ and characteristic number $q=m/n$. Assume that $e_1=u_\lambda$ and let $k$ be the largest natural number relatively prime with $n$ and such that $2k<n$. Then, $\gamma$ has a multiple point of multiplicity $n$ and $n(k-m)$ ordinary double points (see Figures \ref{Example2-1} and \ref{Example3-1}  of Subsections 5.2 and 5.3, respectively).
\end{rem}

We are now in the right position to prove Theorem \ref{B}, which we state here in a more technical form.

\begin{thm}\label{Bcomplete} Let $\gamma$ be a B-string with parameters $(\lambda,e_1)\in\mathcal{P}$ and characteristic number $q=m/n$. Then, the following conclusions hold true:
\begin{enumerate}
\item The trajectory of $\gamma$ is invariant by the group generated by the rotation $2\pi/n$ around the $Ox$-axis and $\gamma$ is contained in the spherical region bounded by the planes $x=1/(2\xi e_1)$ and $x=1/(2\xi e_2)$.
\item If $e_1\neq u_\lambda$, then $\gamma$ does not intersect the $Ox$-axis. Moreover:
\begin{enumerate}
\item If $e_1<u_\lambda$, $n-m$ is the linking number with the $Ox$-axis (equipped with the upward orientation) and $\gamma$ possesses, exactly, $n\left(n-m-1\right)$ ordinary double points.
\item If $e_1>u_\lambda$, $-m$ is the linking number with the $Ox$-axis (equipped with the upward orientation) and $\gamma$ possesses, at least, $n m$ points of self-intersection.
\end{enumerate}
\item If $e_1=u_\lambda$, then $\gamma$ intersects the $Ox$-axis $n$ times and the moving point $\gamma(s)$ travels counter-clockwise around the $Ox$-axis (equipped with the upward orientation). Moreover, $n-m$ is the turning number of the plane projection of $\gamma$ to the plane $x=0$.
\end{enumerate}
\end{thm}
\textit{Proof.} The first assertion is trivial. In fact, the stabilizer of the trajectory of a B-string is the subgroup generated by the monodromy, that is the group generated by the rotation of an angle $2\pi/n$ around the $Ox$-axis. Moreover, from Theorem \ref{Acomplete}, it follows that
$$\gamma(s)=\left(h(s),-\rho(s)\cos\theta(s),\rho(s)\sin\theta(s)\right),$$
where $h(s)$ is defined in \eqref{height}. As customary in our proofs, we are avoiding to write the subscripts $\lambda$ and $e_1$, for simplicity in the expressions. Since $h(s)>0$, the trajectory of the curve is clearly contained in $\mathbb{S}^2_+=\{(x,y,z)\in\mathbb{S}^2\,\lvert\,x>0\}$. Similarly, since $e_2\leq\mu(s)\leq e_1$, from the definition of $h(s)$, one also concludes that the trajectory of the curve lies in the spherical region bounded by the planes of the statement.

Suppose now that $e_1\neq u_\lambda$. In this case, the radial function $\rho$ defined in \eqref{radial} is strictly positive and so the linking number ${\rm Lk}(\gamma,Ox)$ is the winding number of $\widehat{\gamma}=(-\rho\cos\theta,\rho\sin\theta)$ around the origin, that is, the degree of the circle map
$$s\in\mathbb{R}\setminus n\omega\mathbb{Z}\cong\mathbb{S}^1\longmapsto (-\cos\theta,\sin\theta)\in\mathbb{S}^1\,,$$
which is equal to
$$-\frac{n}{2\pi}\left(\theta\left[\frac{\omega}{2}+\omega\right]-\theta\left[\frac{\omega}{2}\right]\right)=-m\,.$$
If $e_1>u_\lambda$, we have $\theta(\omega/2+\omega)>\theta(\omega/2)$ and, if $e_1<u_\lambda$, $\theta(\omega/2+\omega)<\theta(\omega/2)$. This proves the second statement (the number of points of self-intersection follows from Lemma \ref{lem}). 

Finally, suppose that $e_1=u_\lambda$, then $\gamma(s)$ passes through the pole $(1,0,0)$ if and only if $s$ is a zero of the radial function $\rho$. In other words, if and only if, $s\equiv 0({\rm mod}\,\omega)$. This implies that $\gamma$ crosses the $Ox$-axis when $s=0,...,(n-1)\omega$, i.e., $n$ times. 

Recall that in this exceptional case
$$\Psi(e_1)=\theta(\omega)-\theta(0)=2\pi\left(q-\frac{1}{2}\right)=2\pi\left(\frac{m}{n}-\frac{1}{2}\right),$$
where $q=m/n$ is the characteristic number of $\gamma$, and we denote by $\widehat{m}/\widehat{n}=q-1/2$. The plane projection of $\gamma$ is the ``polar" curve
$$\widetilde{\gamma}(s)=\rho(s)\left(-\cos\theta(s),\sin\theta(s)\right).$$
The radial and angular functions satisfy $\theta(s+\omega)-\theta(s)=2\pi\widehat{m}/\widehat{n}$ and $\rho(s+\omega)=-\rho(s)$. Thus, if $\widehat{n}=2k$ and $k$ is an odd integer, the least period of $\widetilde{\gamma}$ is $k\omega$, while if $\widehat{n}=2k$ and $k$ is even, the least period of $\widetilde{\gamma}$ is $\widehat{n}\omega$. Similarly, if $\widehat{n}$ is odd, $\omega_{\widetilde{\gamma}}=2\widehat{n}\omega$. Let $\widetilde{T}$ and $\widetilde{N}$ be the unit tangent and unit normal vector fields along $\widetilde{\gamma}$. Then
$$\widetilde{N}\cdot d\widetilde{T}=-\frac{2\dot{\rho}^2\dot{\theta}+\rho^2\dot{\theta}^3-\rho\ddot{\rho}\,\dot{\theta}+\rho\dot{\rho}\,\ddot{\theta}}{\dot{\rho}^2+\rho^2\dot{\theta}^2}\,ds\,.$$
This implies that
\begin{equation}\label{tc}
\frac{1}{2\pi}\int \widetilde{N}\cdot d\widetilde{T}=\frac{1}{2\pi}\left(f(s)-\theta(s)\right)+c\,,
\end{equation}
where $f$ is the continuous determination of $-\arctan(\rho\,\dot{\theta}/\dot{\rho})$ such that $f(0)=0$. Considering the properties of the radial and angular functions, $\arctan(\rho\,\dot{\theta}/\dot{\rho})$ possesses jump discontinuities at the points $p_k=\omega/2+k\omega$, $k\in\mathbb{Z}$, and is real-analytic elsewhere. At the points of discontinuity we have
\begin{equation*}
\left\{\begin{split} 
&\lim_{s\to p_k^+}\arctan\left(\frac{\rho\,\dot{\theta}}{\dot{\rho}}\right)=\frac{\pi}{2}\\
&\lim_{s\to p_k^-}\arctan\left(\frac{\rho\,\dot{\theta}}{\dot{\rho}}\right)=-\frac{\pi}{2}\\
\end{split}\right.\,.
\end{equation*}
Consequently, $f$ is a quasi-periodic function with quasi-period $\omega$ such that $f(\omega)-f(0)=\pi$. Then, we deduce the following properties from \eqref{tc}:
\begin{enumerate}
\item If $\widehat{n}=2k$ and $k$ is an odd integer, then $2m=\widehat{m}+k$, $n=k$ and ($\omega_{\widetilde{\gamma}}=k\omega$)
$$\frac{1}{2\pi}\int_0^{k\omega} \widetilde{N}\cdot d\widetilde{T}=-\frac{\widehat{m}}{2}+\frac{k}{2}=n-m\,.$$
\item If $\widehat{n}=2k$ and $k$ is an even integer, then $m=\widehat{m}+k$, $n=\widehat{n}$ and ($\omega_{\widetilde{\gamma}}=\widehat{n}\omega$)
$$\frac{1}{2\pi}\int_0^{\widehat{n}\omega}\widetilde{N}\cdot d\widetilde{T} =\frac{\widehat{n}}{2}-\widehat{m}=n-m\,.$$
\item If $\widehat{n}$ is an odd integer, $m=2\widehat{m}+\widehat{n}$, $n=2\widehat{n}$ and ($\omega_{\widetilde{\gamma}}=2\widehat{n}\omega$)
$$\frac{1}{2\pi}\int_0^{2\widehat{n}\omega} \widetilde{N}\cdot d\widetilde{T}=\widehat{n}-2\widehat{m}=n-m\,.$$
\end{enumerate}
This proves that the total curvature, i.e., the turning number of $\widetilde{\gamma}$ is $n-m$, as claimed. \hfill$\square$

\begin{rem} Theorem \ref{Bcomplete} implies that the geometry of a B-string is encoded by its projection onto an annular region (disc-type region if the B-string intersects the $Ox$-axis) of the oriented plane through the origin and orthogonal to the symmetry axis. The multiple points of the B-string are projected onto the multiple points of the projected plane curve and, vice-versa, each multiple point of the projection gives rise to a multiple point of the string. 

If $e_1\leq u_\lambda$, the projection is counter-clockwise oriented, while if $e_1>u_\lambda$ is clockwise oriented. Its symmetry group is the same than the one of the B-string. If $e_1\neq u_\lambda$, the linking number of the string is the homotopy class of the projection, viewed as a plane curve of the annular region.
\end{rem}

\section{Theoretical Aspects}

In this section we comment on the theoretical aspects behind Theorems \ref{Acomplete}, \ref{Ccomplete} and \ref{Bcomplete} and on the rich ``hidden" geometry surrounding the $1/2$-Bernoulli's bending variational problem, which is typical of variational problems related to non-commutative completely integrable Hamiltonian contact systems (\cite{FT,GS1,GM,GS2,J,OR}) and to Liouville integrable geometric variational problems (\cite{LS1}).

\subsection{The Phase Space}

Let 
\begin{equation*}
\Omega:=\begin{pmatrix} 0 & -\varpi_0^1 & -\varpi_0^2 \\ \varpi_0^1 & 0 & -\varpi_1^2 \\ \varpi_0^2 & \varpi_1^2 & 0 \end{pmatrix},
\end{equation*}
be the \emph{Maurer-Cartan form} of $SO(3)$. Then, $\{\varpi_0^1,\varpi_0^2,\varpi_1^2\}$ is a basis for the space of left-invariant $1$-forms. Using the algorithmic procedure illustrated in \cite{GM,G}, the \emph{momentum space} of the functional $\mathcal{B}_\lambda$ is the $5$-dimensional submanifold $M$ of $T^*[SO(3)]$ defined by the embedding
$$\left(A,\mu,\dot{\mu}\right)\in SO(3)\times\mathbb{R}^+\times\mathbb{R}\longmapsto \left(\frac{\mu}{2}+\lambda\right)\varpi_0^1\lvert_A-\frac{\dot{\mu}^2}{\mu^2}\varpi_0^2\lvert_A+\frac{1}{2\mu}\varpi_1^2\lvert_A\in T^*[SO(3)]\,.$$
The restriction of the Liouville form of $T^*[SO(3)]$ to $M$ is the $1$-form
$$\zeta=\left(\frac{\mu}{2}+\lambda\right)\varpi_0^1-\frac{\dot{\mu}^2}{\mu^2}\varpi_0^2+\frac{1}{2\mu}\varpi_1^2\,.$$
The $2$-form $d\zeta$ has maximal rank and its characteristic line bundle is generated by the vector field
$$X=\partial_{\varpi_0^1}+\mu^2\partial_{\varpi_1^2}+\dot{\mu}\partial_\mu+\left(\mu-2\lambda\mu^4-\mu^5+2\frac{\dot{\mu}^2}{\mu}\right)\partial_{\dot{\mu}}\,,$$
where $\left(\partial_{\varpi_0^1},\partial_{\varpi_0^2},\partial_{\varpi_1^2},\partial_\mu,\partial_{\dot{\mu}}\right)$ is the trivialization of $TM$ dual to the coframe $\left(\varpi_0^1,\varpi_0^2,\varpi_1^2,d\mu,d\dot{\mu}\right)$. The integral curves of $X$ are the \emph{canonical lifts} of the B-curves (not necessarily in their standard form) with multiplier $\lambda$, that is curves of the following type
$$\Gamma:s\in\mathbb{R}\longmapsto \left(\left(\gamma,\dot{\gamma},\gamma\times\dot{\gamma}\right),\mu,\dot{\mu}\right)\in M\,,$$
where $\gamma$ is an arc-length parameterized B-curve with multiplier $\lambda$, curvature $\mu^2$ and where $\dot{\mu}$ is the derivative of $\mu$ with respect to the arc-length parameter. Thus, the problem has been reduced to the integration of the characteristic vector field of the $SO(3)$-invariant $2$-form $d\zeta$.

\subsection{The Moment Map}

The \emph{moment map} $\mathfrak{M}:M\longrightarrow so(3)^*$ for the $SO(3)$-action on $(M,d\zeta)$ is the restriction of the moment map for the $SO(3)$-action on $T^*[SO(3)]$ equipped with its standard symplectic form. So, we get
$$\mathfrak{M}:\left(A,\mu,\dot{\mu}\right)\longmapsto {\rm Ad}^*(A)\cdot H(\mu,\dot{\mu})\in so(3)^*\,,$$
where 
$$H(\mu,\dot{\mu})=\left(\frac{\mu}{2}+\lambda\right)\varpi_0^1-\frac{\dot{\mu}^2}{\mu^2}\varpi_0^2+\frac{1}{2\mu}\varpi_1^2\in so(3)^*\,.$$
Using the Killing form we identify $so(3)^*$ with $so(3)$, which is isomorphic to $\mathbb{R}^3$ equipped with the Lie algebra structure defined by the usual vector cross product. Modulo these identifications, the moment map can be written as
$$\mathfrak{M}:\left(A,\mu,\dot{\mu}\right)\in M\longmapsto A\cdot\left(\frac{1}{2\mu},-\frac{\dot{\mu}}{2\mu^2},\frac{\mu}{2}+\lambda\right)^T\in\mathbb{R}^3\,.$$
By construction, $\mathfrak{M}$ is constant along the integral curves of $X$. This implies that if $\gamma$ is a B-curve with multiplier $\lambda$ and parameter $e_1$, then
$$\frac{1}{2\mu}\gamma-\frac{\dot{\mu}}{2\mu^2}\dot{\gamma}+\left(\frac{\mu}{2}+\lambda\right)\gamma\times\dot{\gamma}=\mathcal{J}\,,$$
is constant. Note that $\lVert\mathcal{J}\rVert=\xi$, where $\xi$ is the constant of integration defined in \eqref{xi}. The oriented line passing through the origin and parallel to $\mathcal{J}$ is the axis of symmetry of $\gamma$. The element $\mathcal{J}\in \mathbb{R}^3$ is the \emph{momentum} of the B-curve and so the constant of integration $\xi$ is the length of the momentum. Identifying $\mathbb{R}^3$ with $so(3)$ and letting $\mathcal{J}^*$ be the corresponding  fundamental vector field of $\mathbb{S}^2$, then $\mathcal{J}^*$ is the sum of the (adapted) Killing vector fields along $\gamma$ which arise in the Lagrangian approach (\cite{LS2}). 

The momentum map has maximal rank at each point of $M$ and its image is the open set
$$\Upsilon_\lambda:=\{\mathcal{J}\in\mathbb{R}^3\,\lvert\,\lVert\mathcal{J}\rVert>\frac{1}{2\eta_\lambda}\sqrt{1+\left(\eta_\lambda\left[\eta_\lambda+2\lambda\right]\right)^2}\,\}\subset\mathbb{R}^3\,,$$
where $\eta_\lambda$ is defined in Remark \ref{circles}.

If $\mathcal{J}\in\Upsilon_\lambda$, the quartic polynomial
$$Q_\mathcal{J}(t)=t^4+4\lambda t^3+4\left(\lambda-\lVert \mathcal{J}\rVert^2\right)t^2+1\,,$$
possesses exactly two distinct positive real roots, $e_1>e_2>0$. The integral curves of $X$ with momentum $\mathcal{J}$ are, precisely, the canonical lifts of B-curves with multiplier $\lambda$ and parameter $e_1$.

The \emph{(reduced) phase curve} $\mathcal{C}_\mathcal{J}^*$ of $\mathcal{J}$ is the connected component of the singular algebraic curve
$$C_\mathcal{J}:=\{\left(x,y\right)\in\mathbb{R}^2\,\lvert\,y^2=-x^2Q_\mathcal{J}(x)\}\subset\mathbb{R}^2\,,$$
contained in the half-plane $x>0$. It turns out that $\mathcal{C}_\mathcal{J}^*$ is a smooth curve and $C_\mathcal{J}\otimes\mathbb{C}$ is the affine part of a singular algebraic curve of $\mathbb{CP}^2$. Such a curve is, in general, elliptic and it is rational if and only if either $\lambda=0$ or $\lambda<0$ and $(e_1+e_2)^2=4e_1^3e_2^3$.

\subsection{The Marsden-Weinstein Reduction and the Arnold Connection}

Let $\mathcal{J}\in\Upsilon_\lambda$. The (Marsden-Weinstein) \emph{reduced space} of $\mathcal{J}$ is the $2$-dimensional torus $\mathcal{O}_\mathcal{J}:=\mathfrak{M}^{-1}(\mathcal{J})\subset M$. This torus is invariant by the action of the stabilizer $\mathcal{S}_\mathcal{J}:=\{A\in SO(3)\,\lvert\,A\cdot\mathcal{J}=\mathcal{J}\}\cong SO(2)$ and the vector field $X$ is tangent to $\mathcal{O}_\mathcal{J}$.

The map $\pi_\mathcal{J}:(A,\mu,\dot{\mu})\in\mathcal{O}_\mathcal{J}\longmapsto (\mu,\dot{\mu})\in\mathcal{C}_\mathcal{J}^*$ is $\mathcal{S}_\mathcal{J}$-invariant and gives, on $\mathcal{O}_\mathcal{J}$, the structure of a principal circle bundle. The line bundle of $T\mathcal{O}_\mathcal{J}$ spanned by $X_\mathcal{J}:=X\lvert_{\mathcal{O}_\mathcal{J}}$ defines a connection on the principal bundle $\mathcal{O}_\mathcal{J}\longmapsto\mathcal{C}_\mathcal{J}^*$.

\subsection{Closure Conditions and Integrability by Quadratures}

Let $e_1$ be the largest positive root of $Q_\mathcal{J}(x)$. Then, $(e_1,0)\in\mathcal{C}_\mathcal{J}^*$. Let $\Delta_\mathcal{J}\subset\mathcal{S}_\mathcal{J}$ be the (discrete) holonomy group of the connection, with reference point $(e_1,0)$. We have that $\Delta_\mathcal{J}$ is isomorphic to the monodromy of a B-curve with multiplier $\lambda$ and parameter $e_1$. From a theoretical point of view, the closure condition for a critical curve can be rephrased as follows.

\begin{rem} A B-curve with momentum $\mathcal{J}$ is periodic if and only if $\Delta_\mathcal{J}$ is finite.
\end{rem}

Let $R=(R_1,R_2,R_3)\in SO(3)$ be a positively oriented orthogonal basis such that
$$\mathcal{J}=\frac{1}{2e_1}R_1+\left(\frac{e_1}{2}+\lambda\right)R_3\,,$$
and $\widetilde{\mathcal{O}}_\mathcal{J}\subset\mathcal{O}_\mathcal{J}$ be the holonomy bundle of the connection passing through $(R,e_1,0)$. Then:
\begin{enumerate}
\item The map $\widetilde{\mathcal{O}}_\mathcal{J}\longmapsto\mathcal{O}_\mathcal{J}$ is a covering map with deck transformation group $\Delta_\mathcal{J}$.
\item If $\Delta_\mathcal{J}$ is finite, $\widetilde{\mathcal{O}}_\mathcal{J}\cong\mathbb{S}^1$. Otherwise, $\widetilde{\mathcal{O}}_\mathcal{J}\cong\mathbb{R}$.
\item The set $\widetilde{\mathcal{O}}_\mathcal{J}$ is an integral curve of $X$ and, hence, $(A,\mu,\dot{\mu})\in\widetilde{\mathcal{O}}_\mathcal{J}\longmapsto A_1\in\mathbb{S}^2$ is a B-curve with multiplier $\lambda$ and momentum $\mathcal{J}$.
\end{enumerate}

We next explain how to find the parameterization of a B-curve given in Theorem \ref{Acomplete}. Let $(\lambda,e_1)\in\mathcal{P}$ and
$$\mathcal{J}\equiv\mathcal{J}_{\lambda,e_1}=\xi {\vec i}=\frac{1}{2e_1}\sqrt{1+\left(e_1\left[e_1+2\lambda\right]\right)^2}\,{\vec i}\,.$$
Assume that $\mu:\mathbb{R}\longrightarrow\mathbb{R}^+$ is a (periodic) solution to the Cauchy problem
\begin{equation*}
\left\{\begin{split} &\dot{\mu}^2=-\mu^2Q_\mathcal{J}(\mu) \\
& \mu(\varpi/2)=e_1 \end{split}\right.\,.
\end{equation*}
We denote by $\varphi=(\mu,\dot{\mu}):\mathbb{R}\longrightarrow\mathcal{C}_\mathcal{J}^*$ the parameterization of the phase curve. Choose any cross section of $\varphi^*(\mathcal{O}_\mathcal{J})\longmapsto\mathbb{R}$. This amounts to find $A:\mathbb{R}\longrightarrow SO(3)$ satisfying
$$A\cdot\left(\frac{1}{2\mu},-\frac{\dot{\mu}}{2\mu^2},\frac{\mu}{2}+\lambda\right)^T=(\xi,0,0).$$
Such a map can be easily found by elementary linear algebra. For instance,
\begin{equation*}
A=\begin{pmatrix} \frac{1}{2\xi\mu} & -\frac{\dot{\mu}}{2\xi\mu^2} & \frac{2\lambda+\mu}{2\xi} \\ -\frac{\sqrt{4\xi^2\mu^2-1}}{2\xi\mu} & -\frac{\dot{\mu}}{2\xi\mu^2\sqrt{4\xi^2\mu^2-1}} & \frac{2\lambda+\mu}{2\xi\sqrt{4\xi^2\mu^2-1}} \\ 0 & -\frac{\mu(2\lambda+\mu)}{\sqrt{4\xi^2\mu^2-1}} & -\frac{\dot{\mu}}{\mu\sqrt{4\xi^2\mu^2-1}} \end{pmatrix}
\end{equation*}
satisfies above equation. Clearly, any other cross section is of the form
\begin{equation*}
B=\begin{pmatrix} 1 & 0 & 0 \\ 0 & \cos\Psi & -\sin\Psi \\ 0 & \sin\Psi & \cos\Psi \end{pmatrix}\cdot A\,.
\end{equation*}
Parallel sections of the connection are solutions of
\begin{equation*}
\dot{B}=B\cdot\begin{pmatrix} 0 & -1 & 0 \\ 1 & 0 & -\mu^2 \\ 0 & \mu^2 & 0 \end{pmatrix}\,.
\end{equation*}
It follows from the expression of $A$ and $B$, that $B$ is, precisely, a parallel section if and only if
$$\Psi=2\xi\int\frac{\mu^2(\mu+\lambda)}{1-4\xi^2\mu^2}\,ds\,.$$
Hence, the map $s\longmapsto (B(s),\mu(s),\dot{\mu}(s))$ is a parallel cross section and $s\longmapsto B_1(s)$ is the standard configuration of a B-curve with parameters $(\lambda,e_1)$ exhibited in Theorem \ref{Acomplete}.

We finish this section by reformulating Theorems \ref{Ccomplete} and \ref{Bcomplete} in terms of the momentum.

Let $\lambda\in\mathbb{R}$ and $\Sigma_\lambda:={\rm Im}(\widehat{\Psi}_\lambda)\cap 2\pi\mathbb{Q}$ (which is a countable set containing $2\pi\left((1+p(\lambda),1/2)\cap\mathbb{Q}\right)$). The function $\xi$ is a strictly increasing real-analytic diffeomorphism of $(\eta_\lambda,\infty)$ onto $(\widehat{\eta}_\lambda,\infty)$ where
$$\widehat{\eta}_\lambda:=\frac{1}{2\eta_\lambda}\sqrt{1+\left(\eta_\lambda\left[\eta_\lambda+2\lambda\right]\right)^2}\,.$$
Then, we may choose $\xi$ as a fundamental parameter and express $e_1$ and $\widehat{\Psi}_\lambda$ as functions of $\xi\in (\widehat{\eta}_\lambda,\infty)$. Let $\xi_\lambda^*\in(\widehat{\eta}_\lambda,\infty]$ be defined by $u_\lambda=e_1(\xi_\lambda^*)$ (if $\lambda<0$, then $\xi_\lambda^*=-1/(4\lambda)$, while if $\lambda\geq 0$, we set $\xi_\lambda^*=\infty$). Note that there exist countably many $\lambda<0$ such that $\widehat{\Psi}_\lambda(\xi_\lambda^*)\in\Sigma_\lambda$.

\begin{rem} With this notation the results of Theorems \ref{Ccomplete} and \ref{Bcomplete} can be reformulated as follows:
\begin{enumerate}
\item An arc-length parameterized curve $\gamma$ is a B-string with multiplier $\lambda$ if and only if
$$\frac{1}{2\mu}\gamma-\frac{\dot{\mu}}{2\mu^2}\dot{\gamma}+\left(\frac{\mu}{2}+\lambda\right)\gamma\times\dot{\gamma}=\mathcal{J}\,,$$
is constant and $\widehat{\Psi}_\lambda(\lVert\mathcal{J}\rVert)\in\Sigma_\lambda$.
\item For every $\mathcal{J}$ such that $\widehat{\Psi}_\lambda(\lVert\mathcal{J}\rVert)\in\Sigma_\lambda$, there exists a B-string with momentum $\mathcal{J}$. Two B-strings are equivalent if and only if their momenta have the same length.
\item Let $\widehat{\Psi}_\lambda(\lVert\mathcal{J}\rVert)=2\pi m/n\in\Sigma_\lambda$. Then, the stabilizer of a B-string with momentum $\mathcal{J}$ is generated by the rotation around $\mathcal{J}$ of an angle $2\pi/n$. Moreover:
\begin{enumerate}
\item If $\lVert\mathcal{J}\rVert<\xi_\lambda^*$, $n-m$ is the linking number of the B-string with the oriented axis $\mathcal{A}_\mathcal{J}:=\{O+t\mathcal{J}\,\lvert\,t\in\mathbb{R}\}$. In this case, the B-string possesses, exactly, $n(n-m-1)$ ordinary double points.
\item If $\lVert\mathcal{J}\rVert=\xi_\lambda^*$, the B-string turns counter-clockwise around $\mathcal{A}_\mathcal{J}$ and intersects this axis $n$ times. The turning number of the plane projection of the B-string to the oriented plane through the origin and orthogonal to $\mathcal{A}_\mathcal{J}$ is $n-m$.
\item If $\lVert\mathcal{J}\rVert>\xi_\lambda^*$, $-m$ is the linking number of the B-string with the axis $\mathcal{A}_\mathcal{J}$. In this case, the B-string possesses, at least, $nm$ points of self-intersection. (See Lemma \ref{lem} for more details about the intersection points.)
\end{enumerate}
\end{enumerate}
\end{rem}

\begin{rem} Assuming the ansatz that $\widehat{\Psi}_\lambda$ is strictly increasing, then $\Sigma_\lambda=2\pi\left((1+p(\lambda),1/2)\cap\mathbb{Q}\right)$ and, for every $2\pi m/n\in\Sigma_\lambda$ there exists a unique equivalence class of B-strings with multiplier $\lambda$ such that $\widehat{\Psi}_\lambda(\lVert\mathcal{J}\rVert)=2\pi m/n$.
\end{rem}

\section{Examples}

In this section we will consider several examples which illustrate all the theoretical findings of previous sections.

\subsection{Case $\eta_\lambda<e_1<u_\lambda$:}

In this example we consider a B-string with multiplier $\lambda=1.1$ and characteristic number $q=2/5$ ($n=5$ and $m=2$). Since $\lambda>0$ the B-string is of \emph{negative type} (i.e., $e_1<u_\lambda=\infty$). 

\begin{figure}[h!]
\centering
\begin{subfigure}[b]{0.3\linewidth}
\includegraphics[height=5cm]{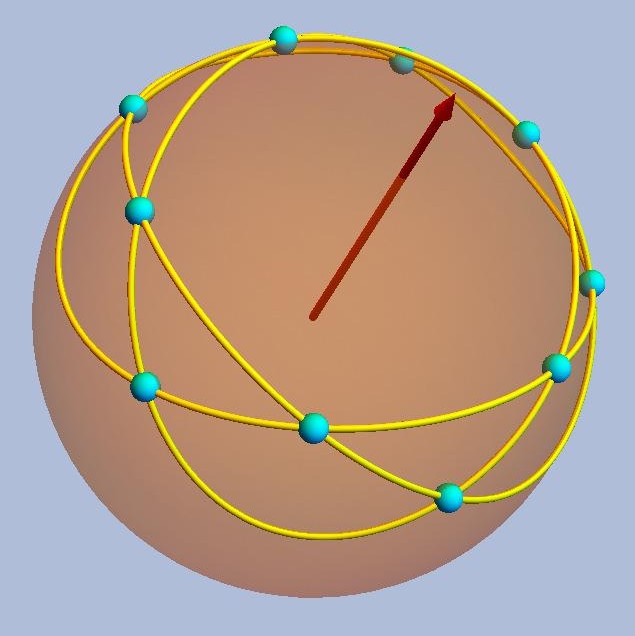}
\end{subfigure}
\quad
\begin{subfigure}[b]{0.3\linewidth}
\includegraphics[height=5cm]{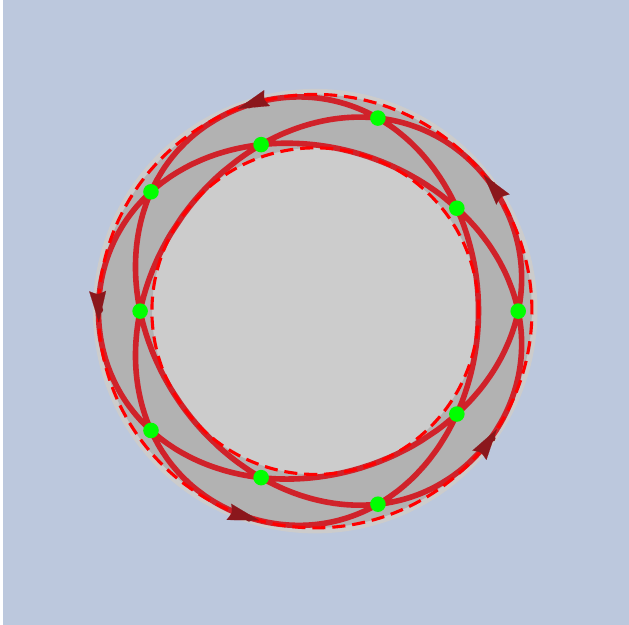}
\end{subfigure}
\quad
\begin{subfigure}[b]{0.3\linewidth}
\includegraphics[height=5cm]{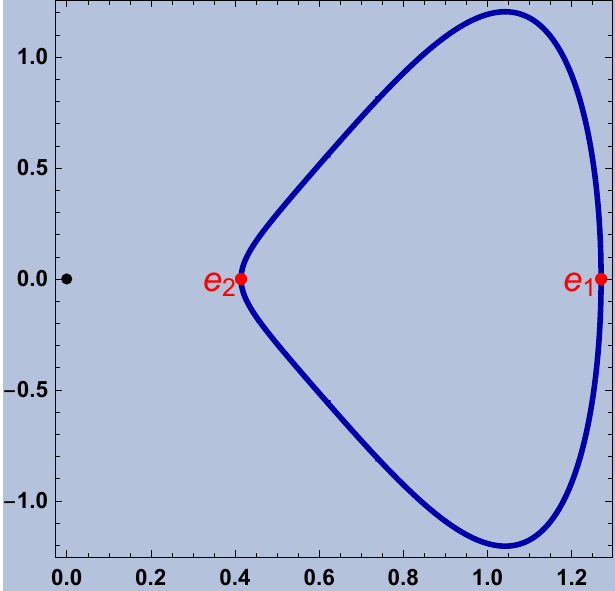}
\end{subfigure}
\caption{A B-string of negative type ($e_1<u_\lambda$) with multiplier $\lambda=1.1$ and characteristic number $q=2/5$, together with its plane projection. On the right: the phase curve, which consists of the isolated singular point (the origin, in black) and a smooth connected component contained in the half-plane $x>0$.}
\label{Example1-1}
\end{figure}

According to Theorems \ref{Ccomplete} and \ref{Bcomplete}, the B-string has a counter-clockwise five-fold symmetry; its linking number with the (upward oriented) $Ox$-axis is $n-m=3$ (which coincides with the winding number and the turning number of the plane projection), and it possesses $n(n-m-1)=10$ ordinary double points. In Figure \ref{Example1-1}, we show the corresponding B-string, its plane projection and its associated phase curve, where we illustrate these properties.

\subsection{Case $e_1=u_\lambda$ (without ordinary double points):}

In this example we consider an \emph{exceptional} B-string with multiplier $\lambda\simeq -0.11$ and characteristic number $q=4/9$ ($n=9$ and $m=4$).

\begin{figure}[h!]
\centering
\begin{subfigure}[b]{0.3\linewidth}
\includegraphics[height=5cm]{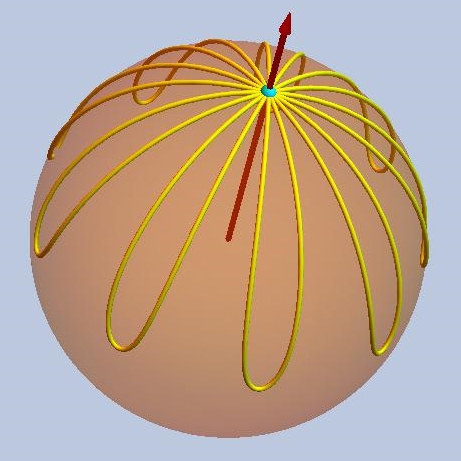}
\end{subfigure}
\quad
\begin{subfigure}[b]{0.3\linewidth}
\includegraphics[height=5cm]{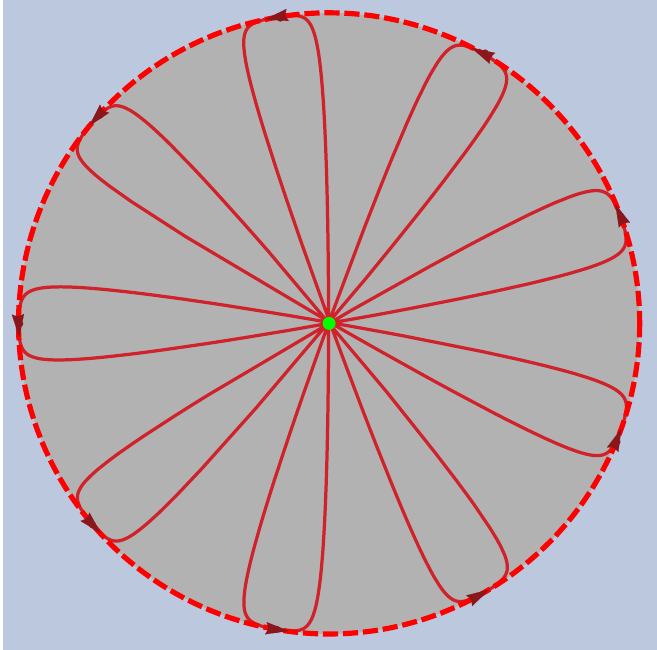}
\end{subfigure}
\quad
\begin{subfigure}[b]{0.3\linewidth}
\includegraphics[height=5cm]{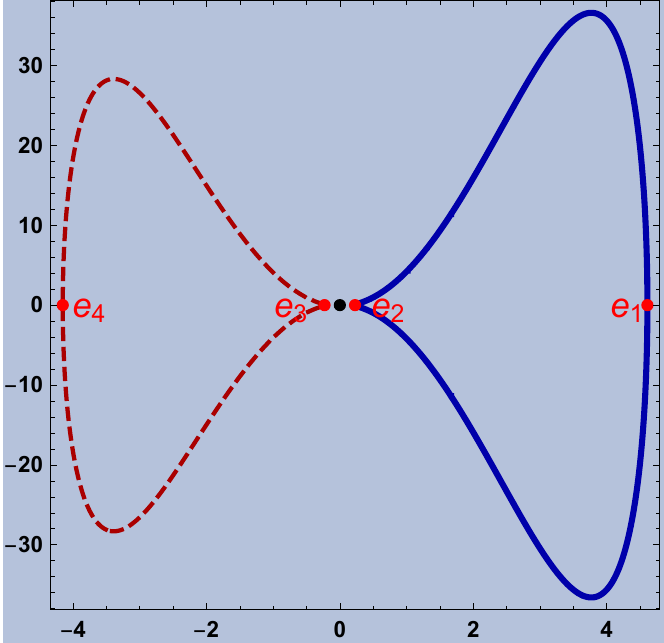}
\end{subfigure}
\caption{An exceptional ($e_1=u_\lambda$) B-string with multiplier $\lambda\simeq -0.11$ and characteristic number $q=4/9$, together with its plane projection. On the right: the phase curve, which consists of the isolated singular point (the origin, in black) and two smooth reduced phase curves. One is contained in the half-plane $x>0$ (colored in blue) and the other one (the dashed dark-red curve) is contained in the half-plane $x<0$.}
\label{Example2-1}
\end{figure}

Observe in Figure \ref{Example2-1} that, as stated in Theorems \ref{Ccomplete} and \ref{Bcomplete}, the B-string possesses a counter-clockwise nine-fold symmetry, the turning number of the plane projection is $n-m=5$ and that it possesses a multiple point (at the pole $(1,0,0)$) of multiplicity $n=9$. Observe that the largest natural number $k$ relatively prime with $n=9$ and such that $2k<n$ is $k=4$ and so, we conclude from Remark \ref{remarknew}, that this B-string has $n(k-m)=0$ ordinary double points.

\subsection{Case $e_1=u_\lambda$ (with $n(k-m)$ ordinary double points):}

In this example we consider an \emph{exceptional} B-string with multiplier $\lambda\simeq -0.45$ and characteristic number $q=2/9$ ($n=9$ and $m=2$).

\begin{figure}[h!]
\centering
\begin{subfigure}[b]{0.3\linewidth}
\includegraphics[height=5cm]{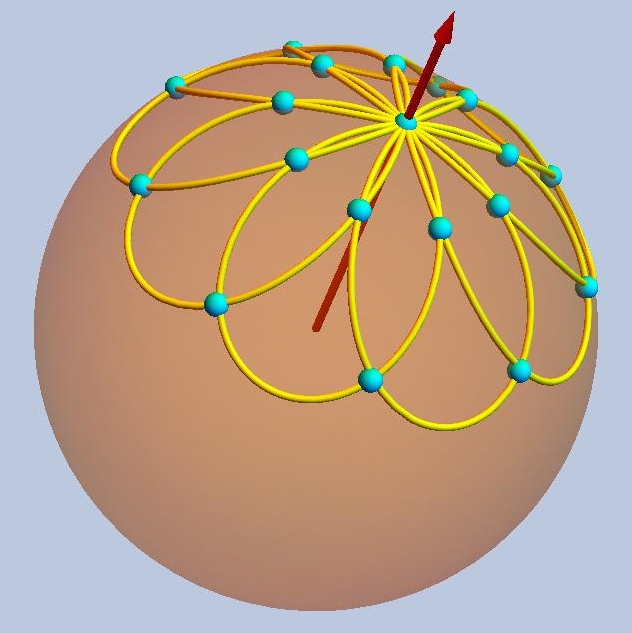}
\end{subfigure}
\quad
\begin{subfigure}[b]{0.3\linewidth}
\includegraphics[height=5cm]{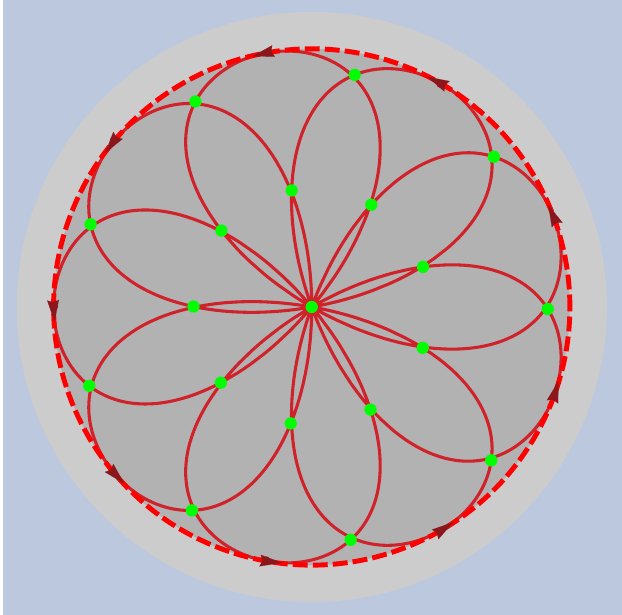}
\end{subfigure}
\quad
\begin{subfigure}[b]{0.3\linewidth}
\includegraphics[height=5cm]{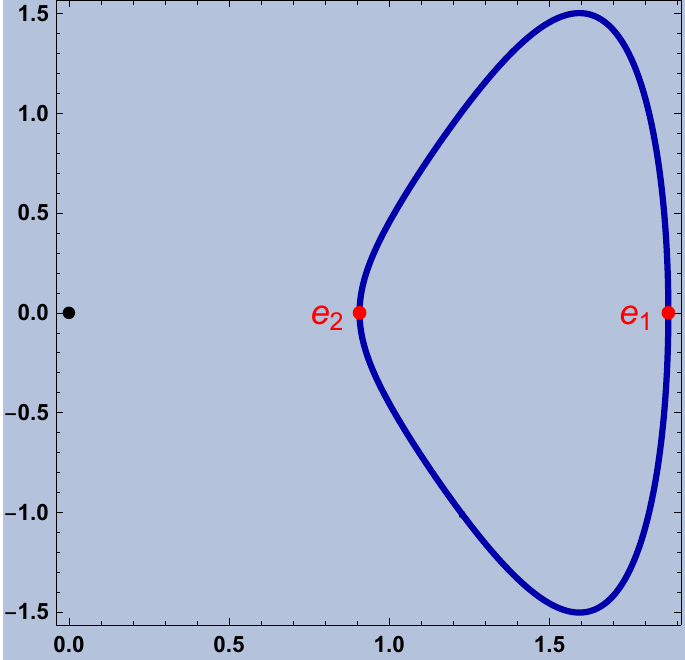}
\end{subfigure}
\caption{An exceptional ($e_1=u_\lambda$) B-string with multiplier $\lambda\simeq -0.45$ and characteristic number $q=2/9$, together with its plane projection. On the right: the phase curve, which consists of the isolated singular point (the origin, in black) and a smooth reduced phase curve contained in the half-plane $x>0$.}
\label{Example3-1}
\end{figure}

In Figure \ref{Example3-1} we can see that, as stated in Theorems \ref{Ccomplete} and \ref{Bcomplete}, the B-string possesses a counter-clockwise nine-fold symmetry, the turning number of the plane projection is $n-m=7$ and that it possesses a multiple point (at the pole $(1,0,0)$) of multiplicity $n=9$. Moreover, the B-string also possesses $n(k-m)=18$ ordinary double points, since in this case the largest natural number $k$ relatively prime with $n=9$ and such that $2k<n$ is $k=4$ (see Remark \ref{remarknew}).

\subsection{Case $e_1>u_\lambda$ (with $nm$ ordinary double points):}

In this example we consider a B-string of \emph{positive type} with multiplier $\lambda= -0.5$ and characteristic number $q=3/8$ ($n=8$ and $m=3$). By positive type, we mean that $e_1>u_\lambda$ holds.

\begin{figure}[h!]
\centering
\begin{subfigure}[b]{0.3\linewidth}
\includegraphics[height=5cm]{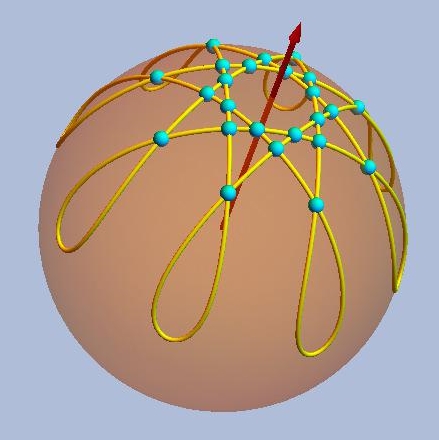}
\end{subfigure}
\quad
\begin{subfigure}[b]{0.3\linewidth}
\includegraphics[height=5cm]{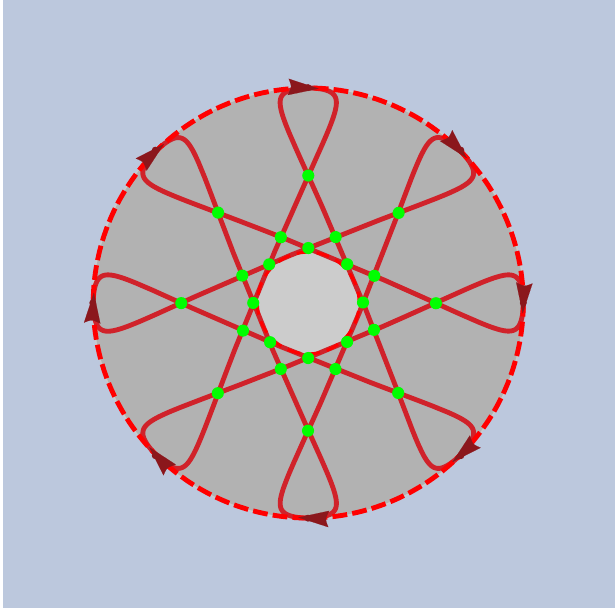}
\end{subfigure}
\quad
\begin{subfigure}[b]{0.3\linewidth}
\includegraphics[height=5cm]{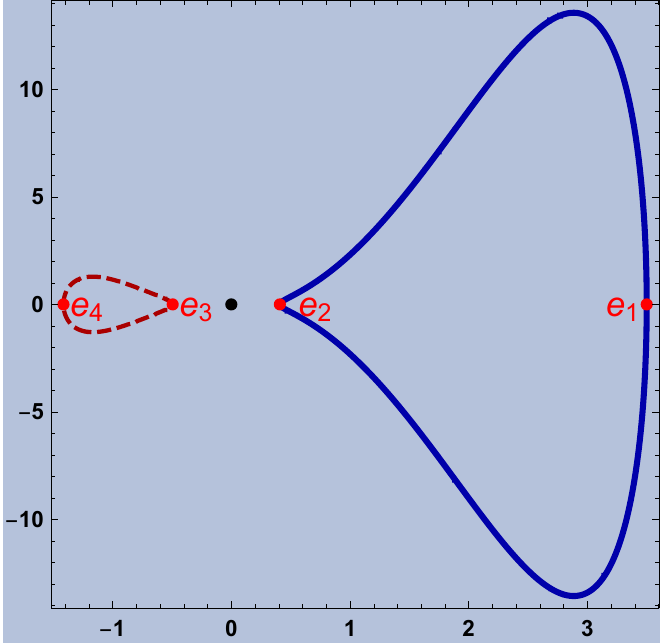}
\end{subfigure}
\caption{A B-string of positive type ($e_1>u_\lambda$) with multiplier $\lambda= -0.5$ and characteristic number $q=3/8$, together with its plane projection. On the right: the phase curve, which consists of the isolated singular point (the origin, in black) and two smooth reduced phase curves. One is contained in the half-plane $x>0$ (colored in blue) and the other one (the dashed dark-red curve) is contained in the half-plane $x<0$.}
\label{Example5-1}
\end{figure}

\begin{figure}[h!]
\centering
\begin{subfigure}[b]{0.3\linewidth}
\includegraphics[height=5cm]{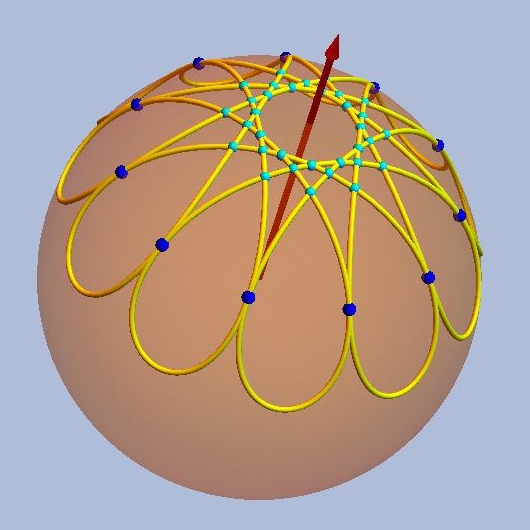}
\end{subfigure}
\quad
\begin{subfigure}[b]{0.3\linewidth}
\includegraphics[height=5cm]{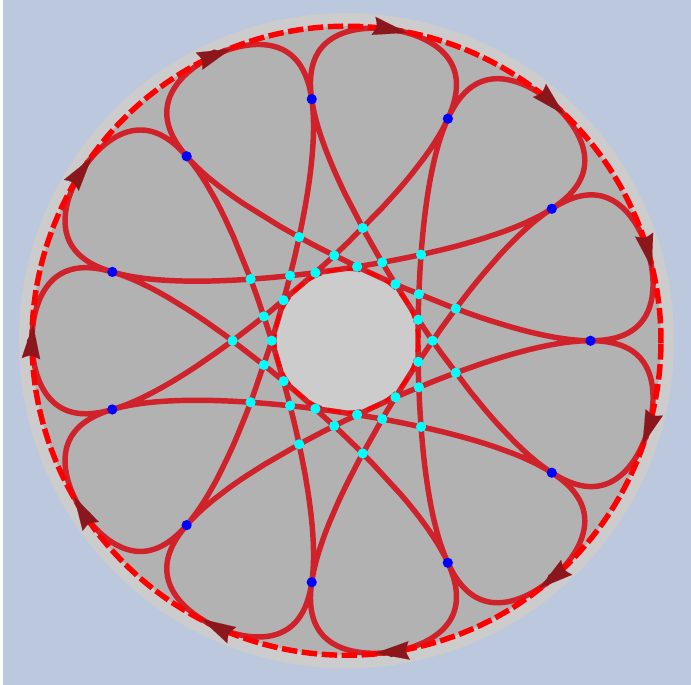}
\end{subfigure}
\quad
\begin{subfigure}[b]{0.3\linewidth}
\includegraphics[height=5cm]{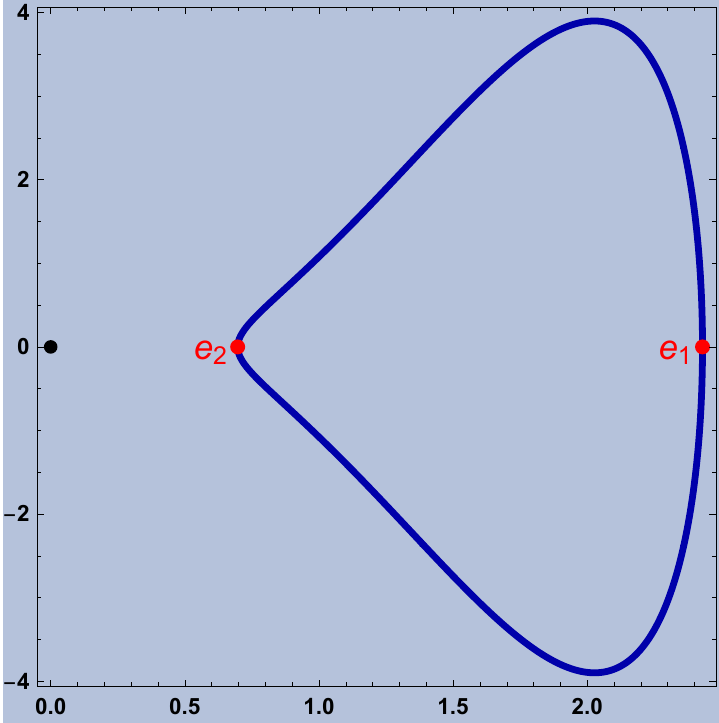}
\end{subfigure}
\caption{A B-string of positive type ($e_1>u_\lambda$) with multiplier $\lambda=-0.5$ and characteristic number $q=3/11$, together with its plane projection. On the right: the phase curve, which consists of the isolated singular point (the origin, in black) and a smooth connected component contained in the half-plane $x>0$.}
\label{Example6-1}
\end{figure}

Observe in Figure \ref{Example5-1} that, as stated in Theorems \ref{Ccomplete} and \ref{Bcomplete}, the B-string possesses a clockwise eight-fold symmetry, the winding number (also, the turning number) of the plane projection is $-m=-3$ and it possesses $nm=24$ ordinary double points.

\subsection{Case $e_1>u_\lambda$ (with tangential double points):}

In this example we consider a B-string of \emph{positive type} with multiplier $\lambda= -0.5$ and characteristic number $q=3/11$ ($n=11$ and $m=3$).

This B-string possesses a clockwise eleven-fold symmetry, the winding number (also, the turning number) of the plane projection is $-m=-3$ and it possesses $nm=33$ ordinary double points and $n=11$ points of tangential self-intersection.

\subsection{Case $e_1>u_\lambda$ (with $n(m+2k)$ ordinary double points):}

In this example we consider a B-string of \emph{positive type} with multiplier $\lambda=-0.5$ and characteristic number $q=2/9$ ($n=9$ and $m=2$).

\begin{figure}[h!]
\centering
\begin{subfigure}[b]{0.3\linewidth}
\includegraphics[height=5cm]{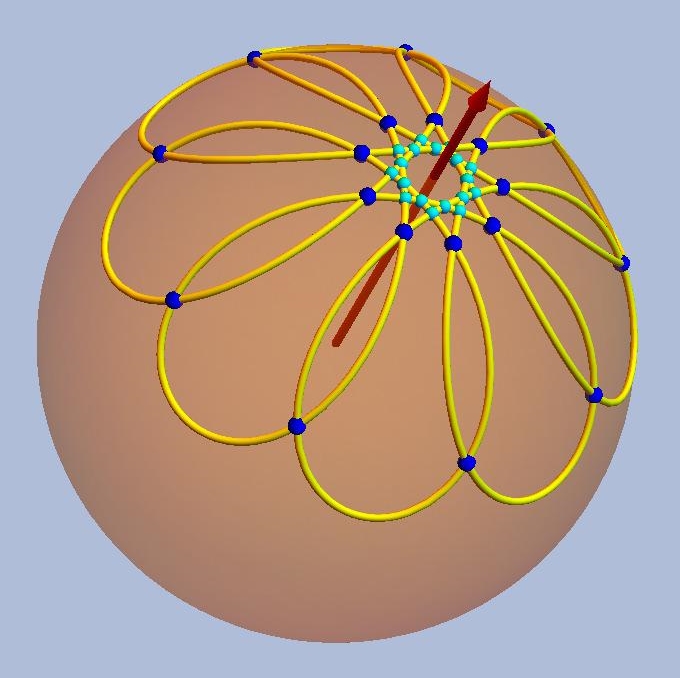}
\end{subfigure}
\quad
\begin{subfigure}[b]{0.3\linewidth}
\includegraphics[height=5cm]{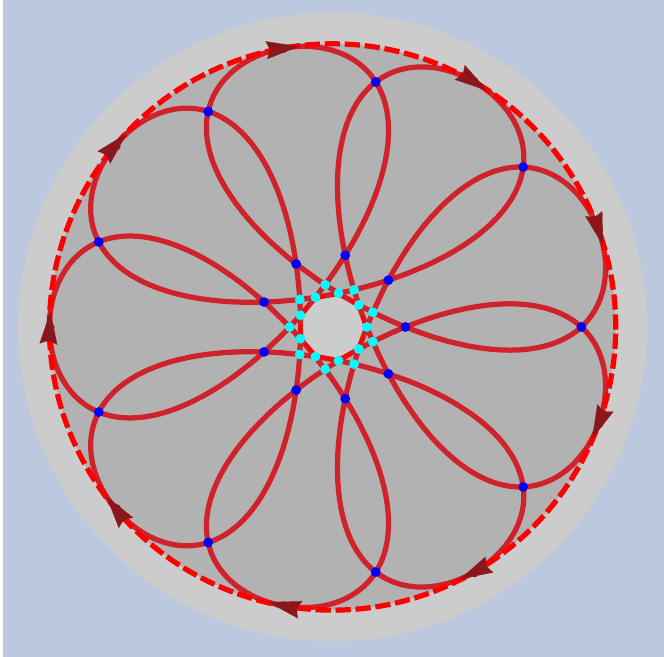}
\end{subfigure}
\quad
\begin{subfigure}[b]{0.3\linewidth}
\includegraphics[height=5cm]{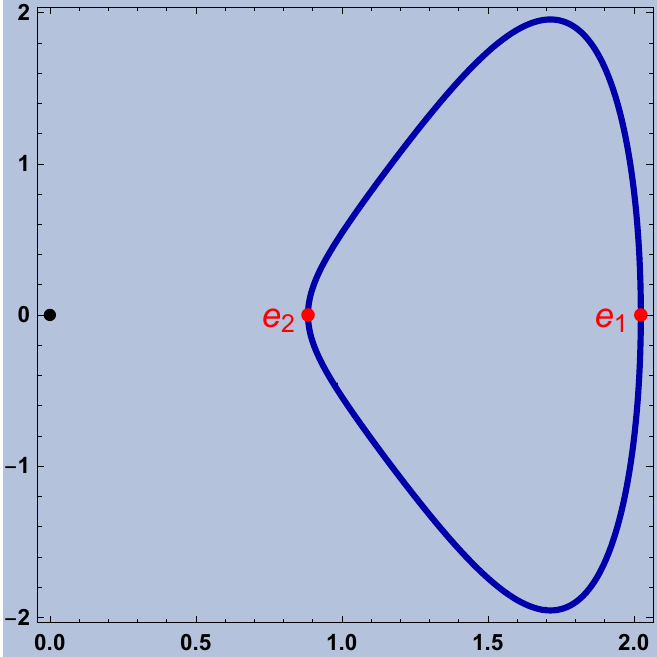}
\end{subfigure}
\caption{A B-string of positive type ($e_1>u_\lambda$) with multiplier $\lambda=-0.5$ and characteristic number $q=2/9$, together with its plane projection. On the right: the phase curve, which consists of the isolated singular point (the origin, in black) and a smooth connected component contained in the half-plane $x>0$.}
\label{Example7-1}
\end{figure}

Observe in Figure \ref{Example7-1} that, as stated in Theorems \ref{Ccomplete} and \ref{Bcomplete}, the B-string possesses a clockwise nine-fold symmetry, the winding number (also, the turning number) of the plane projection is $-m=-2$ and that it possesses $n(m+2k)=36$ ordinary double points, since $k=1$.

\section*{Appendix A. The curvature of the Extrema and the Complete Elliptic Integral $\Psi_\lambda$} 

This appendix has two parts. In the first one we will show how to build the $\mu$-invariant from incomplete elliptic integrals of the third kind and how to compute its least period in terms of complete elliptic integrals. In the second part we will decompose the integral $\Psi_\lambda$ and compute its limits as $e_1$ approaches the boundaries of its domain. 

\noindent{\bf Part I: The Curvature of the Extrema} 

\begin{figure}[h]
\begin{center}
\includegraphics[height=5cm,width=5cm]{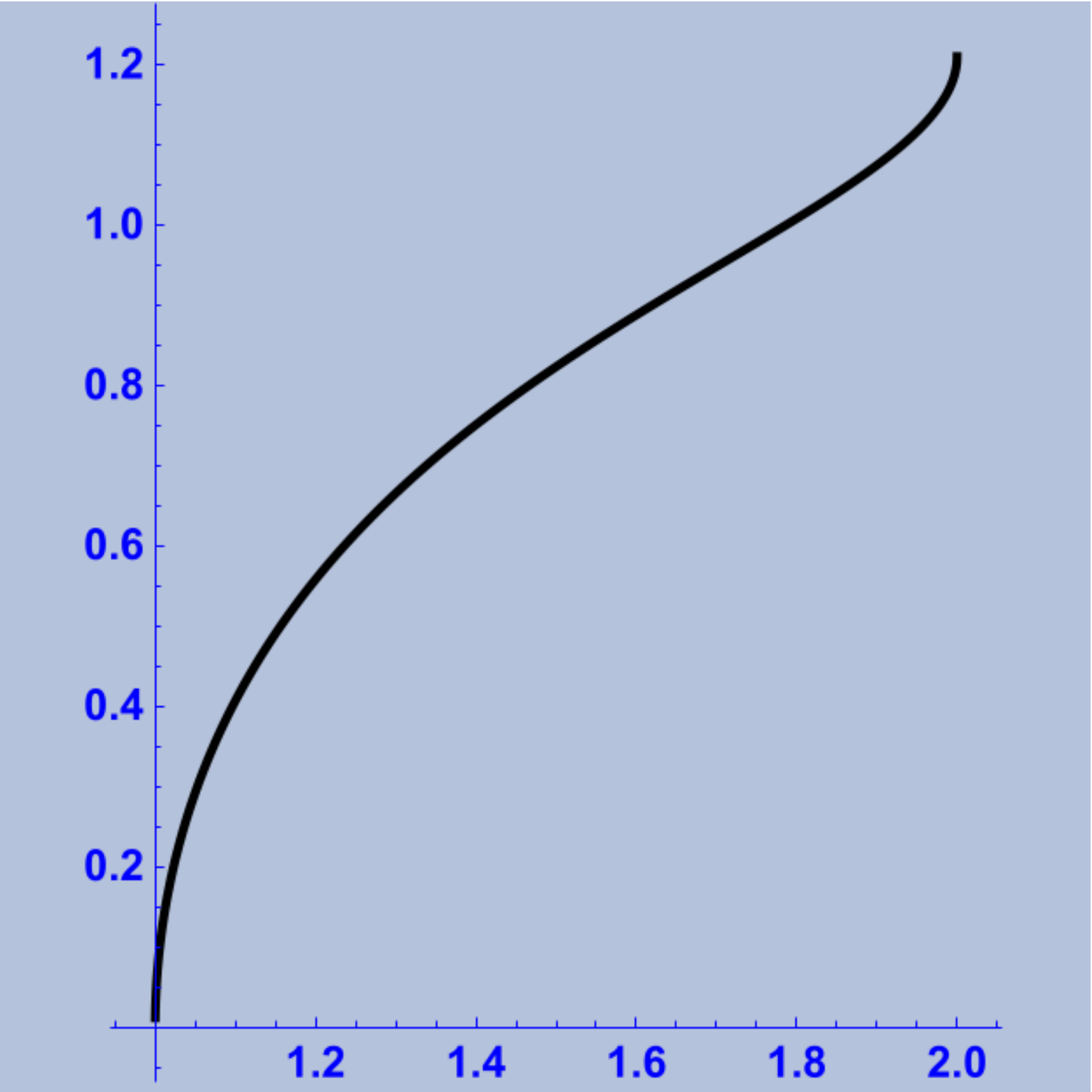}\quad\quad\quad
\includegraphics[height=5cm,width=5cm]{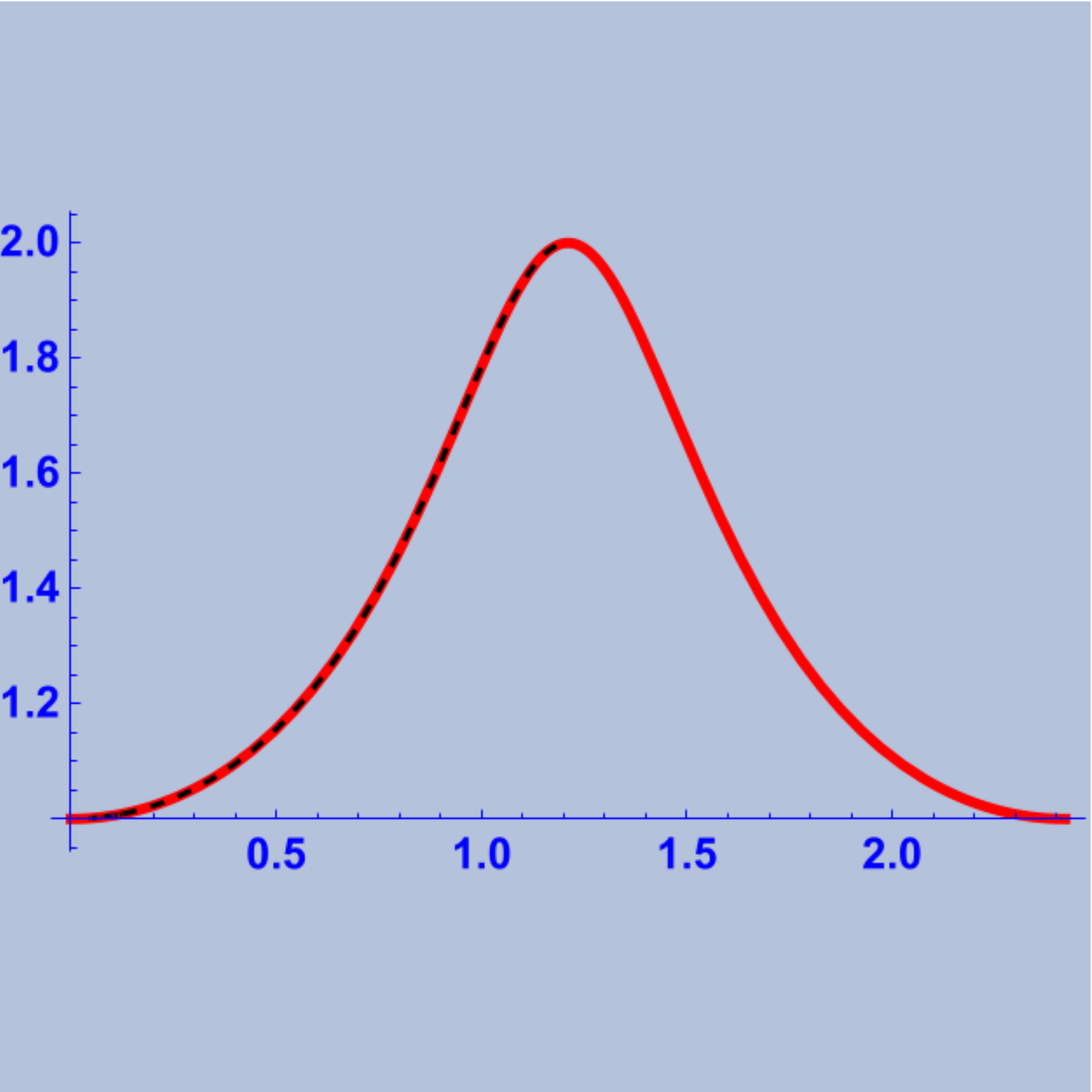}
\caption{\small{On the left: the graph of the $h$-function for $e_1=2$ and $e_2=1$. On the right: the graph of the $\mu$ function for the same values of $e_1$ and $e_2$.
}}\label{FIGmu}
\end{center}
\end{figure}

\noindent Let $K(\phi,\delta)$ and $\Pi(\zeta,\phi,\delta)$ be the Legendre's incomplete elliptic integrals of the first and third kind, defined as 
$$K(\phi,\delta)=\int_0^{\phi} \frac{1}{\sqrt{1-\delta \sin^2 (\theta)}}\,d\theta\,,\quad\quad\Pi(\zeta,\phi,\delta)= \int_0^{\phi} \frac{1}{(1-\zeta \sin^2(\theta))\sqrt{1-\delta\sin^2(\theta)}}\,d\theta\,,$$
and $K(\delta)=K(\pi/2,\delta)$, $\Pi(\zeta,\delta)=\Pi(\zeta, \pi/2, \delta)$ be the corresponding complete elliptic integrals. Let ${\rm am}(u,\delta)$ be the Jacobi's amplitude with parameter $\delta$ and ${\rm sn}(u,\delta) =\sin({\rm am}(u,\delta))$ the associated Jacobi's elliptic function. For simplicity, we denote by
$$\alpha\equiv\alpha(\lambda,e_1)=\frac{e_2-e_1}{e_2-e_4}\,,\quad
\beta\equiv\beta(\lambda,e_1)=\frac{2}{\sqrt{(e_1-e_3)(e_2-e_4)}}\,,\quad
\delta\equiv\delta(\lambda,e_1)=\frac{(e_1-e_2)(e_3-e_4)}{(e_1-e_3)(e_2-e_4)}\,,$$
and 
$$\zeta\equiv \zeta(\lambda,e_1)=\frac{e_4(e_2-e_1)}{e_1(e_2-e_4)},
$$
where $e_1>e_2>0$ and $e_3,e_4$ are the roots of the polynomial $Q$. Recall that $e_2, e_3$ and $e_4$ are functions of the fundamental parameters $\lambda$ and $e_1$.  Let $\mu$ be the solution of (\ref{ODE}) with $\mu(0)=e_2$ and $\omega>0$ be its least period. By construction, $\mu$ is strictly increasing on $[0,\omega/2]$ and $\mu(\omega/2)=e_1$.  
Let $h:[e_2,e_1]\to [0,\omega/2]$ be defined by
$$h(y)= \int_{e_2}^{y}\frac{1}{x\sqrt{-Q(x)}}\,dx =  \frac{\omega}{2}-\int_y^{e_1}\frac{1}{x\sqrt{-Q(x)}}\,dx.$$
Then, from (\ref{ODE}) it follows that $\mu|_{[0,\omega/2]}=h^{-1}$.  Since $\mu$ is even, this is enough to reconstruct $\mu$ on the whole real axis.  Using $257.12$ and $340.04$ of \cite{BF}, we obtain
$$\int_y^{e_1}\frac{1}{x\sqrt{-Q(x)}}\,dx=\frac{\beta}{e_1}\left(\frac{\alpha}{\zeta}u(y)-\frac{\alpha-\zeta}{\zeta}\Pi(\zeta,{\rm am}(u(y),\delta), \delta)\right),$$
where
$$u(y)={\rm sn}^{-1}\left(\sqrt{\frac{(e_2-e_4)(e_1-y)}{(e_1-e_2)(y-e_4)}},\delta \right).
$$
Putting $y=e_2$, we see that the least period of $\mu$ is
$$\omega = 2\frac{\beta}{e_1}\left(\frac{\alpha}{\zeta}K(\delta)-\frac{\alpha-\zeta}{\zeta}\Pi(\zeta,\delta) \right).
$$
Figure \ref{FIGmu} reproduces the graphs of the $h$-function and of the $\mu$-function on the intervals $[e_2,e_1]$ and $[0,\omega]$,  when $e_1=2$, $e_2=1$, $e_3=(-3+i\sqrt{23})/8$, $e_4=\overline{e_3}$.  The $h$-function is evaluated via the Mathematica library of elliptic functions while the $\mu$ function is evaluated  solving numerically (\ref{EL}) with initial conditions $\mu(0)=e_2$ and $\dot{\mu}(0)=0$. The black-dashed portion of the graph of $\mu$ on $[0,\omega/2]$ is obtained by symmetrizing the graph of $h$ with respect to the bisector of the first quadrant, showing that the two methods are in agreement with each other.

The $1/2$-Bernoulli's bending energy of a B-string can also be evaluated in terms of the wave number and complete elliptic integrals of the first kind as
$${\mathcal B}_{\lambda}(\gamma_{m, n})=2n\int_{e_2}^{e_1}\frac{1}{\sqrt{-Q(y)}}\,dy+n\lambda\omega=n\left(2\beta K(\delta)+\lambda\omega\right).$$

\noindent{\bf Part II: The Complete Elliptic Integral $\Psi_\lambda$}  

\noindent Using \eqref{ODE} to make a change of variable in the definition of $\Psi_\lambda$, \eqref{jump}, we have
 $$\Psi_\lambda(e_1)=2\xi\int_0^\omega \frac{\mu^2\left(\mu+2\lambda\right)}{1-4\xi^2\mu^2}\,ds=4\xi\int_{e_2}^{e_1}\frac{\mu\left(\mu+2\lambda\right)}{(1-4\xi^2\mu^2)\sqrt{-(\mu-e_1)(\mu-e_2)(\mu-e_3)(\mu-e_4)}}\,d\mu.$$
This integral can be solved in terms of complete elliptic integrals of the first and third kind. 

For simplicity, we denote by
$$\zeta_+\equiv\zeta_+(\lambda,e_1)=-\frac{(e_1-e_2)(-1+2e_4\xi)}{(e_2-e_4)(-1+2e_1\xi)}\,,\quad\quad
\zeta_-\equiv\zeta_-(\lambda,e_1)=-\frac{(e_1-e_2)(1+2e_4\xi)}{(e_2-e_4)(1+2e_1\xi)}\,.$$

We then have
$$\Psi_\lambda=2\pi\left(I+[1-\chi] II+III\right),$$
where $\chi$ is the indicator function of the exceptional locus $\mathcal{P}_*$ and,
\begin{eqnarray}
I&=&\frac{\beta}{4\pi\xi}\left(-2+\alpha\left[\frac{-1-4\lambda\xi}{\zeta_+(-1+2e_1\xi)}+\frac{1-4\lambda\xi}{\zeta_-(1+2e_1\xi)}\right]\right)K(\delta)\,,\label{I1}\\
II&=&\frac{\beta(\alpha-\zeta_+)(1+4\lambda\xi)}{4\pi\zeta_+\xi(-1+2e_1\xi)}\,\Pi(\zeta_+,\delta)\,,\label{II2}\\
III&=&\frac{\beta(\alpha-\zeta_-)(-1+4\lambda\xi)}{4\pi\zeta_-\xi(1+2e_1\xi})\,\Pi(\zeta_-,\delta)\,.\label{III3}
\end{eqnarray}
These formulas follow from three standard elliptic integrals. The first one (cf. 340.01 and 341.03 of \cite{BF}) is
$$\int_0^{K(\delta)}\frac{1-a\, {\rm sn}^2(u,\delta)}{1-b\, {\rm sn^2}(u,\delta)}\,du=\frac{a}{b}K(\delta)-\frac{a-b}{b}\Pi(b,\delta).$$
The second elliptic integral (cf. 257 and 259 of \cite{BF}) is
$$\int_{e_2}^{e_1}\frac{d\mu}{\sqrt{-(\mu-e_1)(\mu-e_2)(\mu-e_3)(\mu-e_4)}}=\beta K(\delta)\,,$$
where $\beta$ and $\delta$ are as above. The third relevant elliptic integral (cf. 257.39 and 259.04 of \cite{BF}) is
\begin{eqnarray*}\int_{e_2}^{e_1}\frac{d\mu}{(p-\mu)\sqrt{-(\mu-e_1)(\mu-e_2)(\mu-e_3)(\mu-e_4)}}&=&\frac{\beta}{p-e_1}\int_0^{K(\delta)}\frac{1-\alpha\, {\rm sn}^2(u,\delta)}{1-\widetilde{\zeta }\,{\rm sn^2}(u,\delta)}\,du\\&=&\frac{\beta}{p-e_1}\left(\frac{\alpha}{\widetilde{\zeta}}K(\delta)-\frac{\alpha-\widetilde{\zeta}}{\widetilde{\zeta}}\Pi(\widetilde{\zeta},\delta)\right),
\end{eqnarray*}
where $p\neq e_1$, $\alpha$, $\beta$ and $\delta$ are as above, and
$$\widetilde{\zeta}=\frac{(p-e_4)(e_1-e_2)}{(e_1-p)(e_2-e_4)}\,.$$

We begin proving that $\Psi_\lambda\to 2\pi p(\lambda)$ when $e_1$ approaches $\eta_\lambda$ from the right. By construction we have that $e_2\to\eta_\lambda$ too, and, hence, it follows that
\begin{equation*}
\left\{\begin{split} 
&\lim_{e_1\to\eta_\lambda^+} e_3(\lambda,e_1)=\frac{-1+\sqrt{1-\eta_\lambda^4}}{\eta_\lambda^3} \\
&\lim_{e_1\to\eta_\lambda^+} e_4(\lambda,e_1)=\frac{-1-\sqrt{1-\eta_\lambda^4}}{\eta_\lambda^3}
\end{split}\right.\,,
\end{equation*}
and
$$\lambda=\frac{1-\eta_\lambda^4}{2\eta_\lambda^3}\,.$$
We now see that the coefficients in \eqref{I1}-\eqref{III3} tend to, respectively,
$$\lim_{e_1\to\eta_\lambda^+} \alpha(\lambda,e_1)=0\,,\quad\quad\quad\lim_{e_1\to\eta_\lambda^+} \beta(\lambda,e_1)=\frac{2\eta_\lambda}{\sqrt{3+\eta_\lambda^4}}\,,\quad\quad\quad
\lim_{e_1\to\eta_\lambda^+} \delta(\lambda,e_1)=0\,,$$
while
$$\lim_{e_1\to\eta_\lambda^+} \xi(\lambda,e_1)=\frac{\sqrt{1+\eta_\lambda^4}}{2\eta_\lambda^3}\,,\quad\quad\quad\lim_{e_1\to\eta_\lambda^+} \zeta_+(\lambda,e_1)=0\,,\quad\quad\quad\lim_{e_1\to\eta_\lambda^+} \zeta_-(\lambda,e_1)=0\,,$$
and
$$\lim_{e_1\to\eta_\lambda^+} \frac{\alpha}{\zeta_+}=\frac{\eta_\lambda^4\left(\eta_\lambda^2-\sqrt{1+\eta_\lambda^4}\right)}{\sqrt{1+\eta_\lambda^4}+\sqrt{1-\eta_\lambda^8}+\eta_\lambda^6},\,\,\quad\quad \lim_{e_1\to\eta_\lambda^+} \frac{\alpha}{\zeta_-}=-\frac{\eta_\lambda^4\left(\eta_\lambda^2+\sqrt{1+\eta_\lambda^4}\right)}{\sqrt{1+\eta_\lambda^4}+\sqrt{1-\eta_\lambda^8}-\eta_\lambda^6}\,.$$
Finally, recalling that $K(0)=\Pi(0,0)=\pi/2$, using the above limits and (\ref{p}) we conclude that
\begin{equation}\label{fenkelimits}
\lim_{e_1\to\eta_\lambda^+} \Psi_\lambda =-2\pi\sqrt{\frac{1+\eta_\lambda^4}{3+\eta_\lambda^4}}=2\pi p(\lambda).
\end{equation}

In what follows, we prove the limit when $e_1\to\infty$. This limit will depend on the sign of $\lambda$. More precisely, we will see that
$$\lim_{e_1\to\infty} I(\lambda,e_1)=\lim_{e_1\to\infty} III(\lambda,e_1)=0\,,$$
and
\begin{equation*}
\left\{\begin{split} 
&\lim_{e_1\to\infty} II(\lambda,e_1)=-\frac{1}{2}\,,\quad\quad&\text{if}\,\lambda\geq 0\\
&\lim_{e_1\to\infty} II(\lambda,e_1)=\frac{1}{2}\,,\quad\quad&\text{if}\,\lambda<0
\end{split}\right.\,.
\end{equation*}
In order to prove these limits we observe that, as $e_1\to\infty$, the following asymptotic estimates hold true:
\begin{eqnarray*}
&e_2\sim1/e_1\,,\quad e_3\sim -1/e_1\,,\quad e_4\sim-e_1\,,\quad \xi\sim e_1/2\,,\quad \beta\sim 2/e_1\,,\\ 
&\delta\sim 1-4/e_1^2\,,\quad \alpha\sim (1-e_1^2)/(1+e_1^2)\,,\quad \zeta_-\sim 1-4/e_1^2\,,
\end{eqnarray*}
and 
\begin{equation}\label{otherproblems}
\left\{\begin{split} 
&\zeta_+\sim 1-\frac{4\lambda^2}{e_1^4}\,,\quad\quad&\text{if}\,\lambda\neq 0\\
&\zeta_+\sim 1-\frac{1}{e_1^6}\,,\quad\quad&\text{if}\,\lambda=0
\end{split}\right.\,.
\end{equation}
Moreover, recall that $K(\delta)\sim-\log(1-\delta)/2$ as $\delta\to 1^-$ and so, in our case, we have
$$K\left(\delta(\lambda,e_1)\right)\sim\frac{-1}{2}\log\frac{4}{e_1^2}\,,$$
as $e_1\to\infty$. Combining this and the above estimates we conclude that
$$I(\lambda,e_1)\sim \frac{1}{\pi e_1^2}\log\frac{4}{e_1^2}\,,$$ 
as $e_1\to\infty$. This proves the first limit.

For the other limits we need some basic properties of the complete elliptic integral of the third kind $\Pi(\zeta,\delta)$. Let $\Lambda=\{(\zeta,\delta)\in[0,1)\times[0,1)\,\lvert\,\zeta\geq\delta\}$ and consider the function
$$f:(\zeta,\delta)\in\Lambda\longmapsto\frac{2}{\pi}\sqrt{1-\zeta}\sqrt{1-\delta}\,\Pi(\zeta,\delta)\in\mathbb{R}\,.$$
This function $f$ is bounded below by $2/\pi$ and above by $1$. In addition, $f(\zeta,0)=1$ for every $\zeta\in[0,1)$ and $f(\zeta,\delta)\to 1$ when $\zeta\to 1^-$, for every value of $\delta\in[0,1)$. Moreover, for every $\zeta\in[0,1)$,
$$f(\zeta,\zeta)=\frac{2}{\pi}\sqrt{1-\zeta}\,E\left(\frac{\zeta}{\zeta-1}\right)$$
where $E$ is the complete elliptic integral of the second kind. From this we infer that
$$\lim_{\zeta\to 1^-} f(\zeta,\zeta)=\frac{2}{\pi}\,.$$
From these properties we deduce the following facts:
\begin{enumerate}
\item If $\gamma:(a,\infty)\longrightarrow\Lambda$ is a smooth curve such that $\gamma(t)\to(1,1)$ when $t\to\infty$ and $\zeta=\delta$ is an asymptote of $\gamma$ as $t\to\infty$, then
$$\lim_{t\to\infty}f(\gamma(t))=\frac{2}{\pi}\,.$$
\item If $\widetilde{\gamma}:(a,\infty)\longrightarrow\Lambda$ is a smooth curve such that $\widetilde{\gamma}(t)\to(1,1)$ when $t\to\infty$ and $\zeta=1$ is an asymptote of $\widetilde{\gamma}$ as $t\to\infty$, then
$$\lim_{t\to\infty}f(\widetilde{\gamma}(t))=1\,.$$
\end{enumerate}
Combining both things, it follows that as $t\to\infty$,
\begin{equation*}
\left\{\begin{split} 
&\Pi(\gamma(t))\sim \frac{1}{\sqrt{1-\gamma_1(t)}\sqrt{1-\gamma_2(t)}}\\
&\Pi(\widetilde{\gamma}(t))\sim \frac{\pi}{2\sqrt{1-\widetilde{\gamma}_1(t)}\sqrt{1-\widetilde{\gamma}_2(t)}}
\end{split}\right.\,.
\end{equation*}

In view of these properties, fix $\tau_\lambda$ sufficiently large and consider the curves
\begin{equation*}
\left\{\begin{split} 
&\gamma_\lambda:e_1\in(\tau_\lambda,\infty)\longmapsto(\zeta_-(\lambda,e_1),\delta(\lambda,e_1))\in\Lambda\\
&\widetilde{\gamma}_\lambda:e_1\in(\tau_\lambda,\infty)\longmapsto(\zeta_+(\lambda,e_1),\delta(\lambda,e_1))\in\Lambda
\end{split}\right.\,.
\end{equation*}
From above estimates when $e_1\to\infty$, we have the following asymptotic behavior for $\gamma_\lambda$,
$$\gamma_\lambda\sim\left(1-\frac{4}{e_1^2},1-\frac{4}{e_1^2}\right)\,,$$
while for $\widetilde{\gamma}_\lambda$, it depends on the value of $\lambda$,
\begin{equation*}
\left\{\begin{split} 
&\widetilde{\gamma}_\lambda\sim \left(1-\frac{4\lambda^2}{e_1^4},1-\frac{4}{e_1^2}\right),\quad\quad&\text{if}\,\lambda\neq 0\\
&\widetilde{\gamma}_0\sim \left(1-\frac{1}{e_1^6},1-\frac{4}{e_1^2}\right),\quad\quad&\text{if}\,\lambda=0
\end{split}\right.\,.
\end{equation*}
Thus, $\gamma_\lambda\to(1,1)$ and $\zeta=\delta$ is an asymptote of $\gamma_\lambda$ as $e_1\to\infty$. Similarly, $\widetilde{\gamma}_\lambda\to(1,1)$ as $e_1\to\infty$. Hence, it follows from above facts that
\begin{equation*}
\left\{\begin{split} 
&\Pi\left(\zeta_-(\lambda,e_1),\delta(\lambda,e_1)\right)\sim\frac{e_1^2}{e_1^2-4}\\
&\Pi\left(\zeta_+(\lambda\neq 0,e_1),\delta(\lambda\neq 0,e_1)\right)\sim\frac{\pi e_1^3}{8\lvert\lambda\rvert}\\
&\Pi\left(\zeta_+(\lambda=0,e_1),\delta(\lambda=0,e_1)\right)\sim\frac{\pi e_1^4}{4}\\
\end{split}\right.\,,
\end{equation*}
when $e_1\to\infty$.

It is then clear, combining this and above estimates, that
\begin{equation*}
\left\{\begin{split} 
&\lim_{e_1\to\infty}III(\lambda,e_1)=0\\
&\lim_{e_1\to\infty}II(\lambda\geq 0,e_1)=-\frac{1}{2}\\
&\lim_{e_1\to\infty}II(\lambda<0,e_1)=\frac{1}{2}\\
\end{split}\right.\,.
\end{equation*}
This completes the proof about the claimed limits for $\Psi_\lambda$ when $e_1\to\infty$.

\section*{Appendix B. Closed $1/2$-Elasticae in the Plane} 

We briefly comment about $1/2$-elasticae in $\r^2$ in order to clarify some assertions made in the Introduction. We begin with the non-existence of closed convex $1/2$-elasticae other than circles. In $\r^2$ the phase curves for convex $1/2$-elasticae are the singular rational curves (see the picture on the left of Figure \ref{FIGmEuclidean})
$${\mathcal C}_{e_1,e_2}\,:\, y^2+x^4\left(x^2-[e_1+e_2]x + e_1e_2\right)=0\,,$$ 
where $e_1>e_2>0$. Then, following the general argument of the Introduction, $\mu$ is a solution of $\dot{\mu}^2+\mu^4(\mu^2-[e_1+e_2]\mu+e_1e_2)=0$ and $(\mu,\dot{\mu})$ is a periodic, regular parameterization of the smooth component of the phase curve lying in the positive half-plane ${\mathbb H}^2=\{(x,y)\,\lvert\, x>0\}$.

\begin{figure}[h]
\begin{center}
\includegraphics[height=5cm,width=5cm]{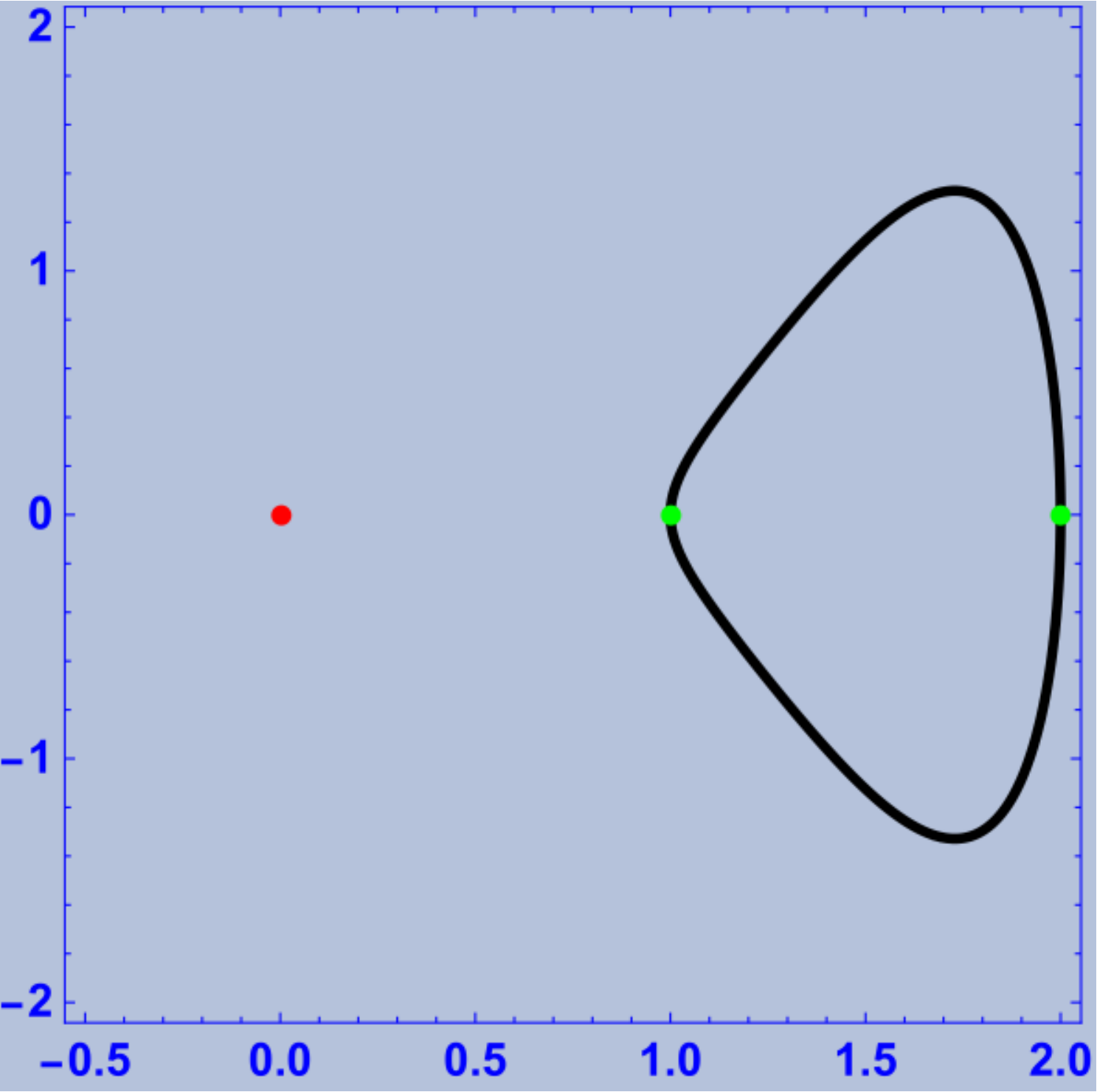}\quad\quad\quad
\includegraphics[height=5cm,width=5cm]{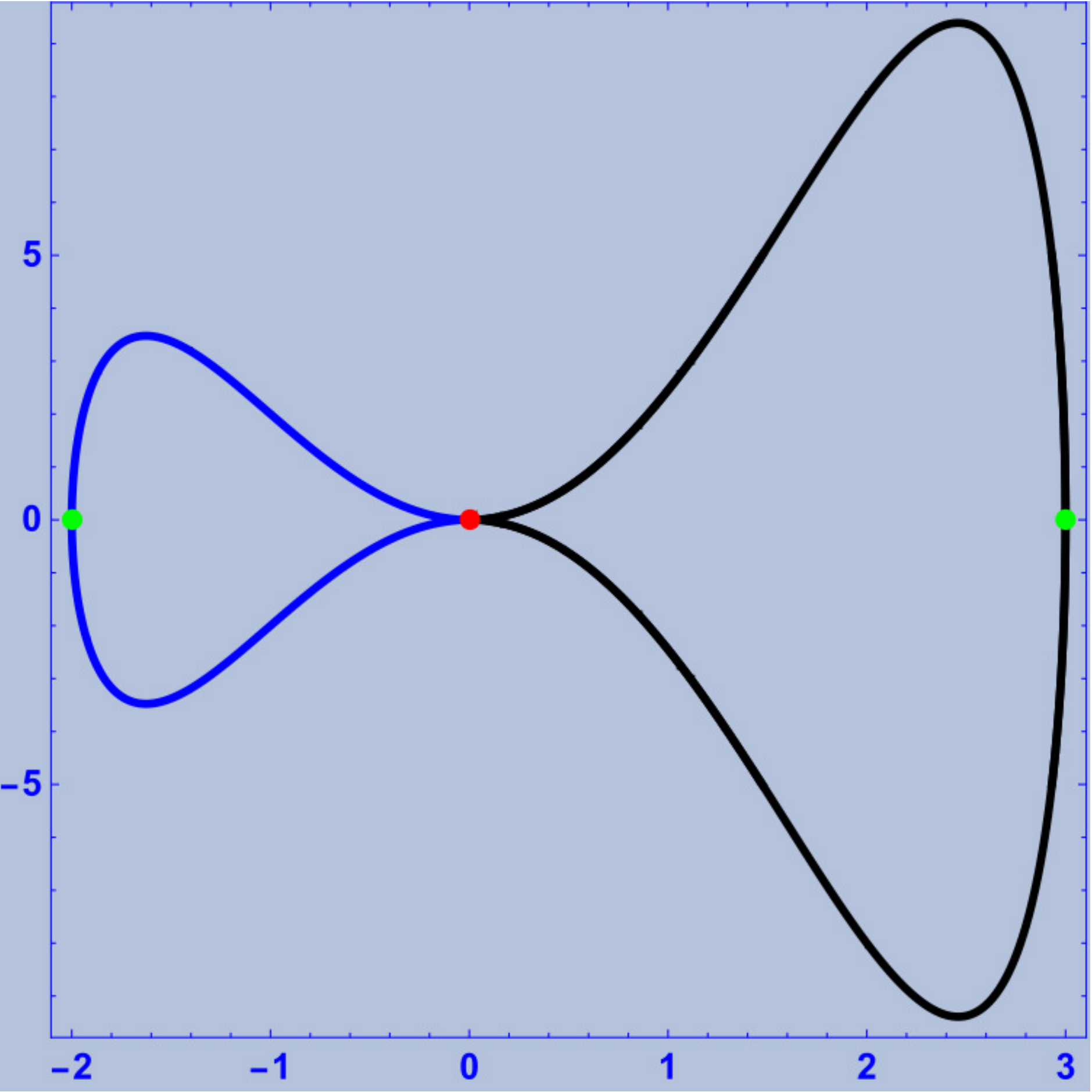}
\caption{Phase curves of convex planar $1/2$-elasticae. On the left, $e_1=2$ and $e_2=1$; while, on the right, $e_1=3$ and $e_2=-2$.}
\label{FIGportrEuclidea}
\end{center}
\end{figure}

Integrating by quadratures, the arc-length parameterization of a critical curve, up to rigid motions, is
$$\gamma_{\lambda,d}(s)=\sqrt{\frac{1}{d}}\left(\lambda s +\frac{1}{2}\int \mu(s)ds,-\frac{1}{2\mu(s)}\right).$$
Let $\omega$ be the least period of $\mu$. Then,
$$\gamma(\omega)-\gamma(0)=\sqrt{\frac{1}{d}}\left(\lambda\omega+\frac{1}{2}\int_0^{\omega}\mu(s)ds,0\right).$$
On the other hand,
$$\omega = 2\int_{e_2}^{e_1}\frac{d\mu}{\mu^2\sqrt{-(\mu-e_1)(\mu-e_2)}}=\pi\frac{e_1+e_2}{(e_1e_2)^{3/2}}\,,$$
and
$$\frac{1}{2}\int_0^{\omega}\mu(s)ds = \int_{e_2}^{e_1}\frac{d\mu}{\mu\sqrt{-(\mu-e_1)(\mu-e_2)}}=\frac{\pi}{\sqrt{e_1e_2}}\,.$$
Hence
$$\gamma(\omega)-\gamma(0)=-\frac{\pi(e_1-e_2)}{(e_1e_2)^{3/2}}(1,0).$$
This implies that  the trajectory of  $\gamma_{\lambda,d}$  is invariant by the subgroup generated by a non-trivial translation along the $Ox$-axis.  In particular, it is unbounded.

\begin{figure}[h]
\begin{center}
\includegraphics[height=5cm,width=5cm]{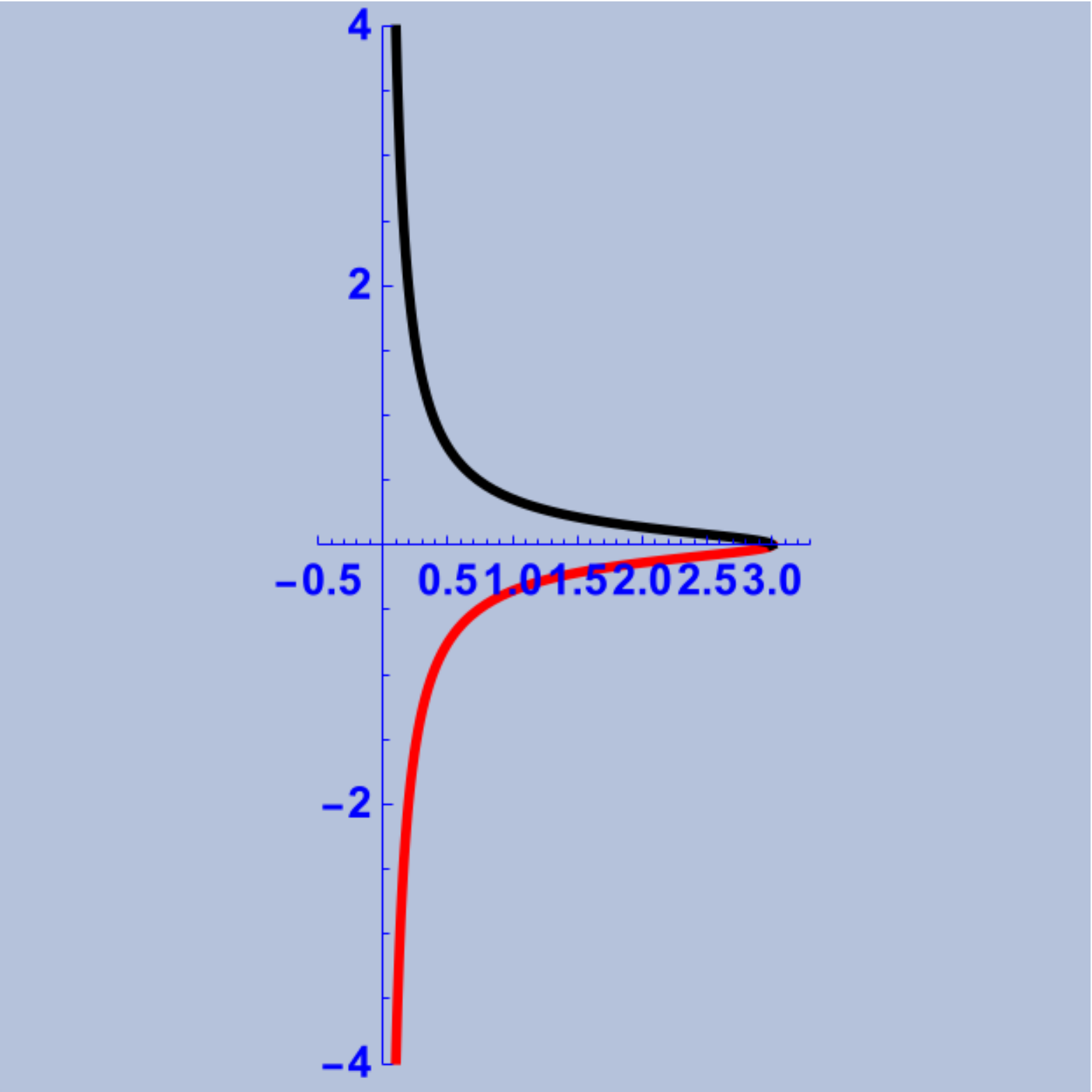}\quad\quad\quad
\includegraphics[height=5cm,width=5cm]{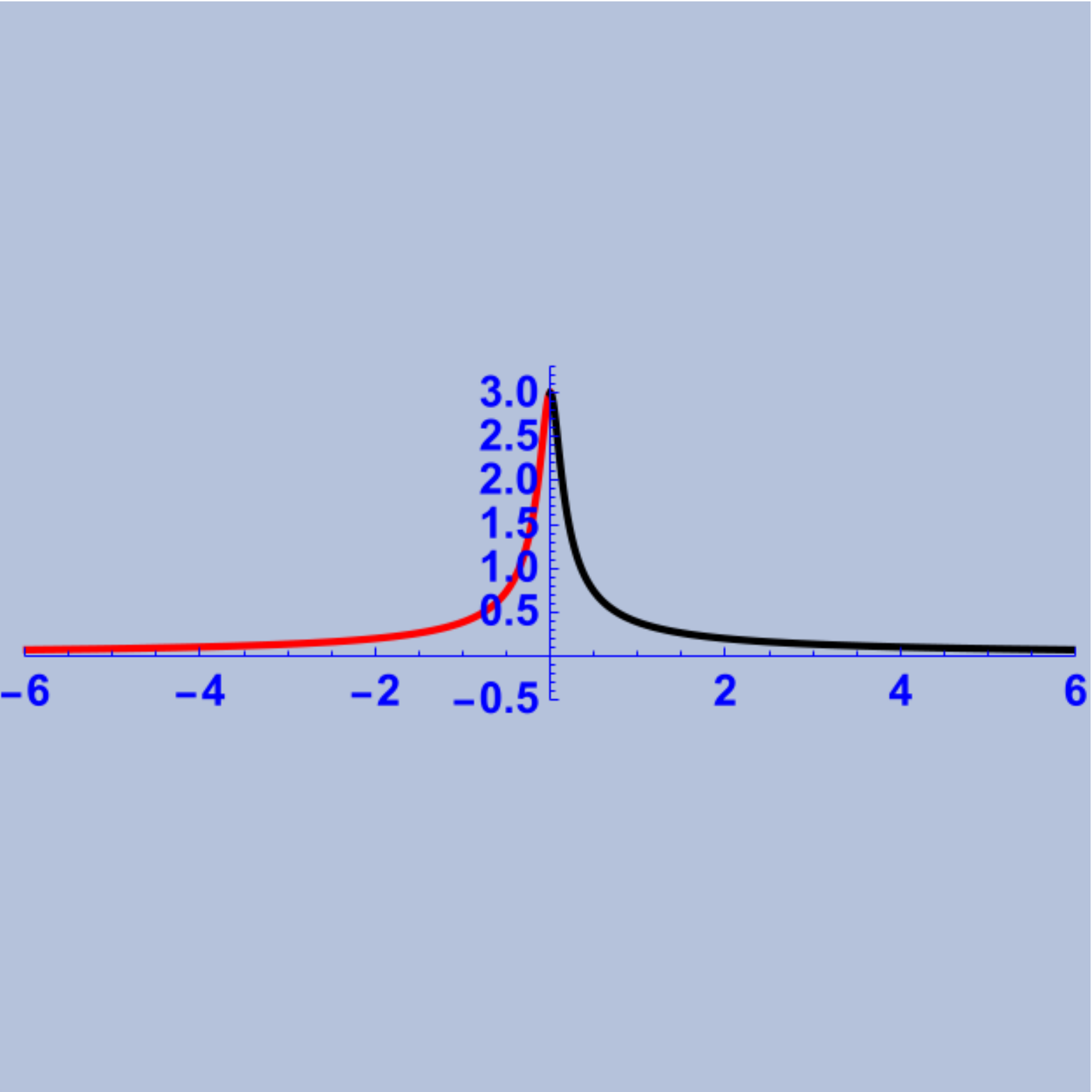}
\caption{On the left, the functions $h_+$ (black) and $h_-$ (red). On the right, the graph of the $\mu$-invariant of a planar $1/2$-elastic curve with $e_1=3$ and $e_2=-2$.}
\label{FIGmEuclidean}
\end{center}
\end{figure}

Next we focus on the non-existence of non-convex critical curves with periodic curvature.  In this case $e_1>0>e_2$.  By contradiction,  suppose that $\kappa$ is periodic and non-constant.  Without loss of generality $\kappa(0)={\rm max}(\kappa)>0$.  Let $J=(a,b)$,  $a<0<b$ be the connected component of $\{s\in {\mathbb R}\,\lvert\, \kappa(s)>0\}$ containing the origin.  Since $\kappa$ is not strictly positive, at least one among $a$ or $b$ is finite. Put $\mu =\sqrt{\kappa\lvert_{J}}$. Then, $m:s\in J\to (\mu,\mu')\in {\mathbb H}^2$ is a parameterization of an open arc contained in ${\mathcal C}^+_{e_1,e_2}:= {\mathcal C}_{e_1,e_2}\cap {\mathbb H}^2$ (see the picture on the right of Figure \ref{FIGportrEuclidea}; ${\mathcal C}^+_{e_1,e_2}$ is represented in the black part). By construction, $\mu(0)=e_1$ and, in addition, there exist $\epsilon>0$ such that $\mu'>0$ on $(-\epsilon,0)$ and $\mu'<0$ on $(0,\epsilon)$. Taking into account that $\mu$ is a solution of the Euler-Lagrange equation, $m$ is an integral curve of the vector field
$$\vec{X}\lvert_{(x,y)}=y\,\partial_{y}+\frac{1}{x}\left(2y^2+\frac{1}{2}[e_1+e_2]x^5-x^6\right)\partial_{y}\in {\mathfrak X}({\mathbb H}^2)\,.$$
Thus, $m(s)$ cannot invert his motion along ${\mathcal C}^+_{e_1,e_2}$. Since ${\mathcal C}^+_{e_1,e_2}$ is homeomorphic to ${\mathbb R}$,  this implies that $m$ is a homeomorphism onto its image. Therefore, $m((a,b))$ intersects the $Ox$-axis at one point, namely at $m(0)=(e_1,0)$. Hence, $\mu'>0$ on $(a,0)$ and $\mu'<0$ on $(0,b)$.  Let $h_{+}$  and $h_{-}$ be the inverses of $\mu\lvert_{(a,0)}$ and $\mu\lvert_{(0,b)}$ respectively (see the picture on the left of Figure \ref{FIGmEuclidean}; the graph of $h_+$ is shown in black while the graph of $h_-$ in red). Then,
\[\begin{split}h_{\pm}(m)=&\int_{e1}^{m}\frac{1}{\mu^2\sqrt{-\mu^2+(e_1+e_2)\mu-e_1e_2}}d\mu =\\
=& \mp\frac{1}{(e_1^2e_2^2)^{3/4}m}\left( \sqrt{e_1|e_2|(e_1-m)(m-e_2)}  -(e_1+e_2)m\,{\rm arctanh} \sqrt{\frac{|e_2|(e_1-m)}{e_1(m-e_2)}} \right).
\end{split}\]
Consequently, $\mu$ can be extended to a function $\hat{\mu}:{\mathbb R}\to (0,e_1]$ attaining its maximum at $s=0$, strictly increasing on $(-\infty,0)$ and strictly decreasing on $(0,+\infty)$ (see the picture on the right of Figure \ref{FIGmEuclidean}) such that
$$\lim_{s \to\-\infty}\hat{\mu}(s)=  \lim_{s \to +\infty}\hat{\mu}(s)=0\,.$$
This implies $a=-\infty$ and $b=+\infty$, which is a contradiction.

\section*{Acknowledgements}

The first author is partially supported by PRIN 2017 ``Real and Complex Manifolds: Topology, Geometry and Holomorphic Dynamics" (protocollo 2017JZ2SW5-004) and by the GNSAGA of INDAM. The present research was also partially supported by MIUR grant ``Dipartimenti di Eccellenza" 20182022, CUP: E11G18000350001, DISMA, Politecnico di Torino.

The authors would like to thank the referees for their valuable comments which have helped to improve the manuscript.

\begin{flushleft}
Emilio M{\footnotesize USSO}\\ 
Dipartimento di Scienze Matematiche (Department of Mathematical Sciences) - Politecnico di Torino, Corso Duca degli Abruzzi 24, I-10129 Torino, Italy\\
E-mail: emilio.musso@polito.it
\end{flushleft}

\begin{flushleft}
\'Alvaro P{\footnotesize \'AMPANO}\\
Department of Mathematics and Statistics, Texas Tech University, Lubbock, TX, 79409, USA\\
E-mail: alvaro.pampano@ttu.edu
\end{flushleft}

\end{document}